\documentclass[11pt,twoside,leqno,euscript]{aomamlt2e} 
  \pageno{825}
\received{May 13, 2002}
  
 \newcounter{notes}
\newenvironment{mitemize}
{\begin{list}{{\rm \arabic{notes}.}}{\usecounter{notes}%
\setlength\parsep{.5ex plus.1ex minus.1ex}
\setlength\itemsep{.5ex  plus.1ex minus.1ex}
\setlength\topsep{6pt  minus2pt}}}
{\end{list}}

 \renewcommand{\thesection}{\arabic{section}}
\makeatletter
\renewcommand{\@seccntformat}[1]{\csname
the#1\endcsname.\hspace{0.5em}\setcounter{Subsec}{0}\setcounter{Subsubsec}{0}}\makeatother

 \newcommand{\mathscr}{\EuScript}
\newtheorem{lemma}{Lemma}[section]
\newtheorem{prop}[lemma]{Proposition}

\newtheorem{cor}[lemma]{Corollary}

\newtheorem{thm}[lemma]{Theorem}

\theorembodyfont{\upshape}
\newtheorem{rem}[lemma]{{\it Remark}}

\newtheorem{notation}[lemma]{{\it Notation}}

\newtheorem{defi}[lemma]{{\it Definition}}

\newtheorem{exo}[lemma]{{\it Example}}
\newtheorem{exos}[lemma]{{\it Examples}}

\def\hb{{\rm H}_{\rm b}}

\def\hbc{{\rm H}_{\rm cb}}

\def\linfty{L^\infty_{\rm w*}}
\def\linftya{L^\infty_{\mathrm{w*,alt}}}
\def\la{L^\infty_{\mathrm{alt}}}

\def\aut{\mathrm{Aut}}
\def\inn{\mathrm{Inn}}
\def\out{\mathrm{Out}}

\def\H{\mathcal{H}}

\def\R{\mathscr{R}}

\def\cc{\mathcal C}
\def\creg{\mathcal C_\mathrm{reg}}
\def\bone{\mathbf{1\kern-1.6mm 1}}
\def\bu{\bullet}
\def\weak{weak-* }
\def\cat#1{{CAT(#1)}}
\def\id{{\rm Id}}
\def\ro{\rho}   
\def\fhi{\varphi}
\def\teta{\vartheta}

\def\ti{-}
\def\lra{\longrightarrow}
\def\Mind#1#2#3{{#1\mathbf{I}}_{#2}^{#3}}
\def\MMind#1#2#3{{#1\mathbf{i}}_{#2}^{#3}}
\def\No{N\raise4pt\hbox{\tiny o}\kern+.2em}
\def\no{n\raise4pt\hbox{\tiny o}\kern+.2em}
\def\bsl{\backslash}
\begin{document}
\currannalsline{164}{2006} 

 \title{Orbit equivalence rigidity\\ and bounded cohomology}

 \acknowledgements{}
\twoauthors{Nicolas Monod}{Yehuda Shalom}

 \institution{University of Chicago,   Chicago, IL\\
\email{monod@math.uchicago.edu}\\
\vglue-9pt
School of Mathematical Sciences, Tel-Aviv University, Tel-Aviv, ISRAEL\\
\email{yeshalom@post.tau.ac.il}}

\centerline{\bf Abstract}
\vglue12pt
We establish new results and introduce new methods in the theory of
measurable orbit equivalence, using bounded cohomology of group
representations. Our rigidity statements hold for a wide
(uncountable) class of groups arising from negative curvature geometry.
Amongst our applications are (a) measurable Mostow-type rigidity
theorems
for products of negatively curved groups; (b) prime factorization
results
for measure equivalence; (c) superrigidity for orbit equivalence; (d)
the
first examples of continua of type $II_1$ equivalence relations with
trivial outer automorphism group that are mutually not stably
isomorphic.

  \vglue6pt
\centerline{\bf Contents}
\def\sni#1{\vglue-1pt\noindent{#1}.\hskip5pt}
\def\ssni#1{\vglue-1pt\noindent\hskip18pt{#1}.}
\sni{1} Introduction
\sni{2} Discussion and applications of the main results
\sni{3} Background in bounded cohomology
\sni{4} Cohomological induction through couplings
\sni{5} Strong rigidity
\sni{6} Superrigidity
\sni{7} Groups in the class $\cc$ and ME invariants
\vglue-1pt\noindent References

\vglue-22pt
\phantom{up}
\section{Introduction}
\label{sec_intro}%
\vglue-4pt

In this paper, a companion to~\cite{Monod-Shalom1}, we continue
our attempts to widen the scope of rigidity theory, using new
techniques made available by the bounded cohomology methods
recently developed by Burger and
Monod~\cite{Burger-Monod3}, \cite{Monod}. In the present paper, we
focus our attention on rigidity of measurable orbit equivalence,
an area which has seen remarkable achievements by R.~Zimmer during
the 80's, and in the last few years has flourished again with the
striking work of A.~Furman~\cite{Furman1}, \cite{Furman2}, \cite{Furman3} and
D.~Gaboriau~\cite{Gaboriau00}, \cite{GaboriauCRAS}, \cite{GaboriauL2}. Our main
purpose is to establish new rigidity phenomena, some reminiscent
of those known in the case of higher rank lattices, for a large
(uncountable) class of groups arising geometrically in the general framework of ``negative curvature'':

\begin{exos}
\label{exos_c}%
Consider the collection of all countable groups $\Gamma$ which admit either:
(i)~A nonelementary simplicial action on some simplicial tree, proper on the set of edges; or
(ii)~A nonelementary proper isometric action on some proper \cat{-1} space; or
(iii)~A nonelementary proper isometric action on some Gromov-hyperbolic graph of bounded valency.
\end{exos}

\emph{Non-Abelian free groups} are outstanding examples of groups in this
class; indeed, the main rigidity results below are already interesting
in that case. Notice that since {\it any nontrivial free product of
two countable groups} is in the list above (unless they are finite of order~$2$),
 this class is uncountable; it also contains the uncountable class
of nonelementary subgroups of Gromov-hyperbolic groups. In particular, 
this collection of groups includes the fundamental group of any closed
manifold of negative sectional curvature.

The Examples~\ref{exos_c} are given as a matter of convenience to make this introduction
more concrete; it is in fact only a certain cohomological property of
these groups which plays a role in our approach. Indeed, we introduce the following:

\begin{notation}
Denote by~$\creg$ {\it the class of countable groups} $\Gamma$ with\break
$\hb^2(\Gamma,\ell^2(\Gamma))\neq 0$.
\end{notation}

This definition refers to the bounded cohomology of $\Gamma$ with coefficients in the regular representation; see 
Sections~\ref{sec_background_hb} and~\ref{sec_in_c} for the relevant background. When stating our results in 
Section~\ref{sec_examples} in full generality, we use a possibly larger class~$\cc$. For the time being, however, 
suffice it to indicate that indeed~$\creg$ is strongly related to the geometric notion of negative curvature,
as the following indicates:

\begin{thm}
\label{thm_exos_c}%
All the groups of Examples~{\rm \ref{exos_c}} belong to~$\creg$.
\end{thm}

This statement can be seen as a cohomological property of negative curvature and relies on the results of~\cite{Monod-Shalom1} complemented with~\cite{Mineyev-Monod-Shalom}. However, we shall offer in Section~\ref{sec_shortcut} a short independent proof that many examples, including free groups, belong to the class~$\creg$.

Before recalling the notion of measurable orbit equivalence, let
us fix the following convention: For a discrete group $\Gamma$ we
say that a standard measure space $(X, \mu)$ is a
\emph{probability $\Gamma$-space} if $\mu(X)=1$ and $\Gamma$
acts measurably on $X$, preserving $\mu$. In this paper, all such
actions are  assumed {\it essentially free\/}; i.e., the stabiliser
of almost every point is trivial.

\begin{defi}
\label{defi_OE}%
Let $\Gamma$ and $\Lambda$  be countable groups and $(X,\mu)$,
$(Y,\nu)$ be probability $\Gamma$\ti\ and $\Lambda$-spaces
respectively. A measurable isomorphism $F:X \to Y$ is said to be
an \emph{Orbit Equivalence} of the actions if for a.e.\ $x \in X$:
$F(\Gamma x)=\Lambda F(x)$, i.e., if $F$ takes almost every
$\Gamma$-orbit bijectively onto a $\Lambda$-orbit. In that
case, the two actions are called \emph{Orbit Equivalent} (OE), and
we say that a (possibly different) isomorphism $\widetilde F:X \to
Y$ induces this orbit equivalence if\break\vskip-12pt\noindent $\widetilde F(\Gamma
x)=F(\Gamma x)$ for a.e.\ $x \in X$.
\end{defi}

The starting point of orbit equivalence rigidity theory lies in
the remarkable lack-of-rigidity phenomenon established by
Ornstein-Weiss~\cite{Ornstein-Weiss} (generalised by Connes-Feldman-Weiss~\cite{Connes-Feldman-Weiss}), 
following H.~Dye~\cite{Dye}, for the class of amenable
groups:\itshape\ Any two ergodic probability measure-preserving actions of
 countable amenable groups are \upshape OE. (Shortly we shall mention another
different motivation for OE rigidity theory, related to geometric
group theory.) To put our main results in a better perspective, we
observe first that this absence of rigidity can be extended also
to some nonamenable groups (see Theorem~\ref{thm_free_continuum}):

\medskip

 \itshape
Any given probability measure-preserving action of a countable free group is orbit equivalent to actions of uncountably many different groups.\upshape

\medskip

Of course, a similar lack of rigidity follows for product actions of
products of free groups. The main point of several of our results
is this: For such product groups, a surprisingly rigid behaviour occurs
if we rule out product actions by the following ergodicity property.

\begin{defi}
\label{defi_Irreducible}%
Let $\Gamma=\Gamma_1\times \Gamma_2$ be a product of countable groups. A $\Gamma$-space
$(X,\mu)$ is called \emph{irreducible} if both $\Gamma_i$ act ergodically on $X$.
\end{defi}

For clarity of the exposition we shall formulate here some of our main results for two factors only,
and in partial generality; Section~\ref{sec_reformulation} contains the general statements.

Observe that irreducibility depends on the given product structure on
$\Gamma$, rather than on $\Gamma$ alone. Among the many natural examples of
irreducible actions, we mention here those we shall make explicit
use of: Bernoulli actions (see below), products of unbounded real
linear groups acting on homogeneous spaces (see Section 2.5 below),
and left-right multiplication actions of products of groups which
are both embedded densely in one compact group (see the proof of
Theorem~\ref{thm_out} below).

\begin{thm}[OE Strong Rigidity -- Products]
\label{thm_strong_rigidity1} Let $\Gamma_1${\rm ,} $\Gamma_2$ be
torsion-free groups in~$\creg${\rm ,} $\Gamma = \Gamma_1 \times \Gamma_2${\rm ,} and let
$(X,\mu)$ be an irreducible probability\break $\Gamma$-space. Let
$(Y, \nu)$ be any other probability $\Gamma$-space \/{\rm (}\/not
necessarily irreducible\/{\rm ).}\/
If the $\Gamma$-actions on $X$ and $Y$ are {\rm OE,} then
they are isomorphic with respect to an automorphism of $\Gamma$.
 More precisely{\rm ,} there is $f \in \aut(\Gamma)$ such that the
orbit equivalence is induced by a Borel isomorphism $F:X \to Y$
with $F(\gamma x) = f(\gamma)F(x)$ for all $\gamma \in \Gamma$ and
a.e.\ $x$.
\end{thm}

Notice that composing an action with a group automorphism yields
an orbit equivalent action, but in general one which is not isomorphic.  
Unlike
the case of higher rank lattices, for some groups covered by the
theorem (such as products of free groups), there is an abundance
of such automorphisms which should be ``detected''. As observed in
Section~\ref{sec_reformulation} below, Theorem~\ref{thm_strong_rigidity1}
is not valid in general if the groups are not in the
class~$\creg$.

Using  Theorem~\ref{thm_strong_rigidity1} we are able to produce
the first examples of finitely generated groups outside the distinguished family of
higher rank lattices in semi-simple Lie groups, possessing
infinitely many nonorbit equivalent actions (see also the ``exotic'' 
infinitely generated groups in~\cite{BG81}). In fact we show
more:

\begin{thm}[Many groups with many actions]

\label{thm_continuum}%
There exists a continuum $\,2^{\textstyle\aleph_0}$ of finitely
generated torsion-free groups{\rm ,} each admitting a continuum of
measure-preserving free actions on standard probability spaces{\rm ,}
such that no two actions in this whole collection are orbit
equivalent.
\end{thm}

Although we are able to include products of
(non-Abelian) free groups in this family, it is still an open
problem to produce infinitely many mutually nonorbit equivalent
actions of one free group.

(Added in proof: D.~Gaboriau and S.~Popa have since obtained a continuum of 
non-OE actions
of a free group~\cite{Gaboriau-Popa}, while G.~Hjorth established that all infinite Kazhdan groups 
share this property~\cite{Hjorth_Dye}.)

To proceed one step further, we recall the following notion:

\begin{defi}
\label{defi_mild_mixing}%
A measure-preserving  action of a group $\Lambda$ on a measure
space $(Y,\nu)$ is called \emph{mildly mixing} if there are no
nontrivial recurrent sets, i.e., if for any measurable
$A\subseteq X$ and any sequence $\lambda_i\to \infty$ in
$\Lambda$, one has $\nu(\lambda_i A\triangle A)\to 0$ only when
$A$ is null or co-null.
\end{defi}

Here is now a \emph{superrigidity}-type result:

\begin{thm}[OE superrigidity for products
--  torsion free case] \ Let\break $\Gamma =\Gamma_1 \times \Gamma_2$ and
\label{thm_superrigidity0}%
$(X,\mu)$ be as in Theorem~{\rm \ref{thm_strong_rigidity1}.} Let
$\Lambda$ be {\rm any torsion-free} countable group and let
$(Y,\nu)$ be any mildly mixing probability $\Lambda$-space.

\nobreak If the $\Gamma$- and $\Lambda$-actions are {\rm OE}  then
$\Lambda$ is isomorphic to $\Gamma${\rm ,} and the actions on $X,Y$ are
isomorphic \/{\rm (}\/with respect to an isomorphism $\Gamma\cong \Lambda${\rm ).}
\end{thm}

Actually we prove a more general statement, dropping the
torsion-freeness assumption on $\Lambda$, thereby allowing
``commensurable situations''. We state here the following result,
which is generalised further in Section 2:

\begin{thm}[OE superrigidity --
product]
\label{thm_superrigidity1} \hskip-6pt
 Let $\Gamma =\Gamma_1 \times \Gamma_2$ and $(X,\mu)$ be
as in Theorem~{\rm \ref{thm_strong_rigidity1}.} Let $\Lambda$ be {\rm
any} countable group and let $(Y,\nu)$ be any mildly mixing
probability $\Lambda$-space.
  If the $\Gamma$- and $\Lambda$-actions are {\rm OE} then both
the groups $\Gamma$ and $\Lambda${\rm ,} as well as the actions{\rm ,} are
commensurable. More precisely\/{\rm :}\/

\begin{itemize}
\ritem{(i)} There exist  a finite index subgroup $\Gamma_0 < \Gamma$\/{\rm ,}\/ whose projections to both
factors
$\Gamma_i$ are onto\/{\rm ,}\/
 a finite normal subgroup $N \lhd \Lambda$ with $|N|=[\Gamma : \Gamma_0]$\/{\rm ,}\/ and a
short exact sequence
 
\hfil $1 \to N \to \Lambda \to \Gamma_0 \to 1$ \hfill

\vglue-9pt
\noindent 
such that\/{\rm :}\/

  \ritem{(ii)} The $\Gamma$-action induced from the
$\,\Lambda/N \cong \Gamma_0$-action on $(N\bsl Y, \nu)$ is {\rm
isomorphic} to its action on $(X, \mu)$ \/{\rm (}\/with respect to an
automorphism of $\Gamma$\/{\rm ).}\/
\end{itemize}

In particular{\rm ,} if either the $\Gamma$-action on $X$ is
aperiodic \/{\rm (}\/i.e.{\rm ,} remains ergodic under any finite index subgroup\/{\rm ),}\/
or $\Lambda$ is torsion-free{\rm ,} then $\Lambda$ is isomorphic to
$\Gamma$ and the actions on $X,Y$ are isomorphic \/{\rm (}\/with respect to
an isomorphism $\Gamma\cong \Lambda${\rm )}.
\end{thm}

This theorem is optimal in the sense that   any $\Lambda$
satisfying~(i) above admits an action which is OE to an
irreducible action of $\Gamma$. A crucial ingredient in the proof
of this theorem is a remarkable idea of A.~Furman
from~\cite{Furman1} in the framework of simple Lie groups, which
we adapt here for our purposes. In Example~\ref{exo_mild_mixing}
below we show by means of a counter-example why the mild mixing
condition is natural in our context, and how
Theorem~\ref{thm_superrigidity1} may fail for actions which are
weakly mixing, and ``close to being'' mildly mixing. Of course, the
simplest examples of mildly mixing actions are (strongly) mixing
actions, and those exist for any group, as in the following
standard construction: For a countable group $\Gamma$ and any
probability distribution $\mu$ (different from Dirac) on the interval $[0,1]$, call the
natural shift $\Gamma$-action on the product space
$([0,1]^\Gamma, \mu^\Gamma)$ a \textit{Bernoulli} $\Gamma$-action. Any such action can easily be seen
to be mixing, and this takes care at the same time of irreducibility and aperiodicity. We
therefore have:

\begin{cor}
\label{cor_Bernoulli}%
Let $\Gamma = \Gamma_1 \times \Gamma_2$ where each $\Gamma_i$ is a
torsion-free group in~$\creg$. If a Bernoulli
$\Gamma$-action is orbit equivalent to a Bernoulli $\Lambda$-action for some {\rm arbitrary} group
$\Lambda${\rm ,} then $\Gamma$ and $\Lambda$ are isomorphic and
 with respect to some isomorphism $\Gamma\cong\Lambda$ the actions are
 isomorphic by a Borel isomorphism which induces the given orbit
 equivalence.\hfill\qed
\end{cor}

As shown by the result of Ornstein and Weiss cited above, amenable
groups share a sharp lack of rigidity in the measurable orbit
equivalence theory. Our next two results are analogous to two of
the  theorems above, only that here we replace the setting of
products by one involving amenable radicals. We show a similar
rigid behaviour modulo the intrinsic lack of rigidity caused by
the presence of such radicals.

\begin{thm}[OE Strong Rigidity -- Radicals]
\label{thm_strong_rigidity2}%
 Let $\Gamma$ be a group and $M \lhd \Gamma$ a normal
amenable subgroup such that the quotient $\bar \Gamma = \Gamma /M$
is torsion-free and in~$\creg$. Let $(X,\mu)$,
$(Y,\nu)$ be  probability $\Gamma$-spaces on which $M$ acts
ergodically.
 If the two $\Gamma$-actions are {\rm OE} then there is a
Borel isomorphism $F:X \to Y$ such that for all $\gamma \in
\Gamma$ and a.e.\ $x \in X$: $F(\gamma M x)= f (\bar \gamma) M
F(x)${\rm ,} where $\bar \gamma=\gamma M$ and $f$ is some automorphism
of $\bar \Gamma$.
\end{thm}

Here is the superrigidity-type version:

\begin{thm}[OE Superrigidity -- Radicals]
\label{thm_superrigidity2}%
Let $\Gamma$ and $(X,\mu)$ be as in
Theorem~{\rm \ref{thm_strong_rigidity2}.} Let $\Lambda$ be {\rm any
countable group} and let $(Y,\nu)$ be any mildly mixing probability
$\Lambda$-space.
 If the $\Gamma$\ti\ and $\Lambda$-actions are {\rm OE} then there
exists an infinite normal amenable subgroup $N \lhd \Lambda$ such
that $\Lambda / N$ is isomorphic to $\Gamma /M$.
 Moreover{\rm ,} there is an isomorphism $f: \Gamma /M \to \Lambda /N$
such that the {\rm OE} is induced by a Borel isomorphism $F:X \to Y$
satisfying $F(\gamma M x)= f (\bar \gamma) N F(x)$.
\end{thm}

In a different direction, we can apply
Theorem~\ref{thm_strong_rigidity1} to study countable ergodic
relations of type $\mathrm{II}_1$. We first recall some
terminology (see also~\cite{Furman2}, \cite{Furman3}).

Let $\Gamma$ be a countable group and $(X,\mu)$ be an ergodic
probability $\Gamma$-space. Let $\R = \R_{\Gamma,X}\subseteq
X\times X$ denote the (type $\mathrm{II}_1$) equivalence relation
on $X$ defined by that action, i.e.\ $(x,y)\in\R$ if and only if $\Gamma x
=\Gamma y$. Two such relations are isomorphic if and only if the two actions
are OE. Further, the group of automorphisms $\aut(\R)$ of the
relation $\R$ is the group of measure-preserving isomorphisms $F:
X\to X$ such that $F(\Gamma x) = \Gamma F(x)$ for a.e.\ $x\in X$.
Moreover, one defines the inner and outer automorphism groups by
$$\inn(\R) = \big\{F\in\aut(\R) \ :\  F(x)\in\Gamma x\ \ \mu\mathrm{-a.e.}\big\}, \kern.5cm \out(\R) = \aut(\R)/\inn(\R).$$
While $\inn(\R)$ (the so-called \emph{full group}) is always very
large (e.g.\ it acts essentially transitively on the collection of
all measurable subsets of a given measure), it is of interest to
find relations~-- or group actions~-- for which $\out(\R)$ is
small, or even trivial. The first construction of some
$\R_{\Gamma,X}$ with trivial outer automorphism group is due to
S.~Gefter~\cite{Gefter93}, \cite{Gefter96}. Recently
A.~Furman~\cite{Furman3} has produced more examples within a
comprehensive study of the problem in the setting of higher rank
lattices (these are used, along with Zimmer's cocycle
superrigidity, by both authors).  Furman constructs a continuum of
mutually nonisomorphic type $\mathrm{II}_1$ relations with
trivial outer automorphism group which are all \emph{weakly
isomorphic} (see (i) in Definition~\ref{defi_WOE} below), being
obtained by restricting one fixed relation $\R_{\Gamma,X}$ to
subsets of different measure. We show the following:

\begin{thm}[Many Relations with Trivial Out]
\label{thm_out}%
There exists a continuum of mutually non weakly isomorphic relations of type $\mathrm{II}_1$  with trivial outer automorphism group.
\end{thm}

As mentioned earlier, the study  of orbit equivalence can be
motivated also from an entirely different point of view, being a
measurable counterpart to geometric (or quasi-isometric)
equivalence of groups. This analogy, as well as the following
notion, were suggested by M.~Gromov~\cite[0.5.E]{Gromov}:

\begin{defi}
\label{defi_ME}%
Two countable groups $\Gamma,\Lambda$ are called \emph{Measure
Equivalent} (ME) if there is a standard (infinite-) measure space
$(\Sigma,m)$ with commuting measure-preserving $\Gamma$\ti\ and
$\Lambda$-actions, such that each one of the actions admits a
finite measure fundamental domain. (In particular, both actions
are free, even though not necessarily the product action -- see
also Remark~\ref{rem_essentially_free} below.) The space
$(\Sigma,m)$ endowed with these actions is called an ME
\emph{coupling} of $\Gamma$ and $\Lambda$.
\end{defi}

The analogy with geometric  group theory can be seen as follows:
Replacement of  $\Sigma$  in  in Definition~\ref{defi_ME} by a locally
compact space on which $\Gamma$ and $\Lambda$ act properly,
continuously and co-compactly, in a commuting way, results in a
notion strictly equivalent to $\Gamma$ being quasi-isometric to
$\Lambda$, see~\cite[0.2.C]{Gromov}.

On the other side, ME relates back to OE because of the following
fact, observed by Zimmer and Furman (see Section~\ref{sec_ME_WOE}
below): For two discrete groups $\Gamma$ and $\Lambda$,  admitting
some OE actions is equivalent to having an ME coupling where the
two groups have the same co-volume. (The case of arbitrary
co-volumes corresponds to \emph{weak orbit equivalence} which we
actually cover in all of our results, but preferred not to discuss
in the introduction~-- see Section~\ref{sec_examples} below.)
Thus, results concerning orbit and measure equivalence can be
transformed one to the other (a fact we shall take advantage of,
following Furman's approach), and may both come under the title
``measurable group theory''~-- a counterpart to geometric group
theory.

\begin{thm}[ME Rigidity -- Factors]
\label{thm_factors}%
 Let $\Gamma = \Gamma_1 \times
\cdots \times \Gamma_n$ and $\Lambda = \Lambda_1 \times \cdots
\times \Lambda_{n'}$ be products of
torsion-free countable groups. Assume that all the
$\Gamma_i$\/{\rm '}\/s are in~$\creg$. If $\Gamma$ is {\rm ME} to $\Lambda$\/{\rm ,}\/ then
$n \ge n'$\/{\rm ,}\/ and if equality holds then{\rm ,} after permutation of the indices{\rm ,}
$\Gamma_i$ is {\rm ME} to $\Lambda_i$ for all $i$.
\end{thm}

This may be viewed as a far reaching extension of the phenomena
established by R.~Zimmer~\cite{Zimmer83} and
S.~Adams~\cite{Adams94a} to the effect that the orbit relation
generated by ``negatively curved'' groups is not a product
relation. Illustrating the analogy with geometric group theory, we
point out that the arguments of
Eskin-Farb~\cite{Eskin-Farb1}, \cite{Eskin-Farb2} or
Kleiner-Leeb~\cite{Kleiner-Leeb} can be used to show that if two
products of nonelementary hyperbolic groups are quasi-isometric,
then so are the factors (after permuting indices).

For amenable radicals we have the following analogue:

\begin{thm}[ME Rigidity -- Quotients by Radicals]
\label{thm_kernel_intro}%
\hskip-8pt
Let
$\Gamma,\Lambda$ be countable groups and let $M\lhd \Gamma, N\lhd
\Lambda$ be amenable normal subgroups such that
$\bar\Gamma=\Gamma/M$ and $\bar\Lambda=\Lambda/N$ are in~$\creg$
and are torsion-free. If $\,\Gamma$ is {\rm ME} to $\Lambda${\rm ,} then $\bar\Gamma$ is {\rm ME}
to $\bar\Lambda$.
\end{thm}

\phantom{up}
\vglue-20pt

As mentioned earlier, our new approach to orbit equivalence
rigidity uses notably the new approach to bounded cohomology recently
developed by\break Burger-Monod~\cite{Burger-Monod3}, \cite{Monod}.
The latter provides both results as well as ``working tools'' which
turn out to be very effective in the setting of measurable orbit
equivalence. Aiming the paper at the broader audience interested in orbit
equivalence rigidity, we shall assume here no prior familiarity
with bounded cohomology, and present in
Section~\ref{sec_background_hb} below a friendly and brief
introduction to this theory, including the main results that we
need from Burger-Monod's work. Suffice it to say at this point
that we define  (second) bounded cohomology  similarly to usual
(second) group cohomology, but using bounded cochains. 
As a by-product of our proofs, we get some new
cohomological invariants of measure equivalence, and consequently
some additional ``softer'' rigidity results,  as in the following (see Corollary~\ref{cor_cc_ME_invariants}):

\vglue-20pt
\phantom{up}

\begin{thm}
\label{thm_l2}%
The vanishing of the second bounded cohomology
with coefficients in the regular representation is an {\rm ME} invariant.
\end{thm}

\phantom{up}
\vglue-46pt
\phantom{up}
\begin{cor}
\label{cor_hyperbolic}%
A countable group containing an infinite normal\break amenable
subgroup is not {\rm ME} to any group in~$\creg$.
\end{cor}
\phantom{up}
\vglue-22pt

It follows for instance that such a group cannot be ME to any (nonelementary) Gromov-hyperbolic group; the latter statement was established for the particular case of infinite \emph{center} by S.~Adams~\cite{Adams95}.

\vskip4pt {\it Related results}.
In the framework of reducibility of Borel relations,\break G.~Hjorth and A.~Kechris~\cite{Hjorth-Kechris} 
established rigidity results for certain types of products in independent work carried out at about the same time. 

\vskip4pt {\it Acknowledgments}. It is our pleasure to thank Alex
Furman for many illuminating and helpful discussions on the
material of this paper. His approach to the subject,
particularly attacking orbit equivalence rigidity through the
notion of measure equivalence~\cite{Furman2} and the beautiful
idea of how to deduce superrigidity-type results from strong
rigidity-type results~\cite{Furman1}, substantially influenced this
work. We also use the opportunity to thank again the Mathematical
Institute at Oberwolfach, the FIM at the ETH-Zurich, and the
Mathematics Institute at the Hebrew University in Jerusalem, for
supporting and hosting mutual visits. The second author's travel to 
Oberwolfach was supported by the Edmund Landau Center for Research in Mathematical Analysis and Related 
Areas, sponsored by the Minerva Foundation. He also acknowledges the 
ISF support made through grant 50-01/10.0.

\vskip4pt
    {\it Added in Proof\/}: Since the acceptance of this paper for 
      publication, many new
      results in the emerging measurable group theory appeared, particularly 
      with the ground-breaking work of S. Popa. We refer the reader to the
      accounts [Po] [Sh2] for further details and references.

\section{Discussion and applications of the main results}
\label{sec_examples}%
\vglue-12pt

\Subsec{Weak orbit equivalence and measure equivalence} \label{sec_ME_WOE} 
In this subsection,  we recall some basic facts about the relation
between orbit and measure equivalence, which will enable us to
reformulate a number of our main results in the stronger form in
which they will be proved. The material of this subsection
follows~\cite[\S\S 2--3]{Furman2} wherein the reader can find more
details and proofs. As a matter of notation, we shall use only
\emph{left} actions and cocycles.

We  recall our standing convention that $(X,\mu)$ is called a
\emph{probability\break $\Gamma$-space} if it is a standard
probability space with an essentially free measurable $\Gamma$-action preserving $\mu$. Thus all corresponding measurable
equivalence relations will be of type $\mathrm{II}_1$.

\begin{defi}[Weak Orbit Equivalence]
\label{defi_WOE}%
 Let $\Gamma$  and $\Lambda$ be
countable\break groups and $(X,\mu)$, $(Y,\nu)$ be probability
$\Gamma$\ti\ and $\Lambda$-spaces respectively. The two actions
are said to be \emph{weakly orbit equivalent} (WOE) or stably
orbit equivalent, if either one of the following two equivalent
conditions holds:

\begin{itemize}
\item[(i)] The two equivalence relations induced by the $\Gamma$- and $\Lambda$-actions 
are \emph{weakly isomorphic}, i.e.,
there exist nonnull  measurable subsets $A\subseteq X$,
$B\subseteq Y$ on which the restrictions of the relations are
isomorphic. More precisely, for some $A,B$ as above, a measurable
isomorphism $F:A\to B$ and all $x_1,x_2\in A$, one has $\Gamma x_1
\cap A= \Gamma x_2 \cap A$ if and only if $\Lambda F(x_1) \cap B
=\Lambda F(x_2) \cap B$.

\item[(ii)] There exist measurable maps $p:X\to Y$, $q:Y\to X$ such
that:

\begin{itemize}
\item[1.] $p_* \mu\prec \nu$, $q_*\nu \prec\mu$ (where $\prec$
denotes absolute continuity of measures).

\item[2.] $p(\Gamma x)\subseteq \Lambda p(x)$ and $q(\Lambda y)\subseteq \Gamma q(y)$ for a.e.\ $x\in X$, $y\in Y$.

\item[3.] $q\circ p(x)\in\Gamma x$ and $p\circ q (y) \in \Lambda y$ for a.e.\ $x\in X$, $y\in Y$.
\end{itemize}
\end{itemize}
\end{defi}

Orbit  equivalence as defined in the introduction is of course a
special case of WOE with $A,B$ of full measure in~(i) or with
$p,q$ inverse measurable isomorphisms in~(ii). As we shall see, WOE is
a useful notion even if one is interested in OE only.

\begin{defi}[Compression Constant] With assumptions and notation  as
in Definition~\ref{defi_WOE}, one defines the \emph{compression
constant}
$$C(X, Y) = \nu(B)/\mu(A),$$
where $A,B$ are as in point~(i) of Definition~\ref{defi_WOE}. The compression constant depends on the given WOE but not on the choice of $A,B$.
\end{defi}

\begin{prop}
With notation  as above{\rm ,} assume that the $\Gamma$\ti\ and\break $\Lambda$-actions on $X,Y$ are
ergodic. Then $C(X, Y)=1$ if and only if the actions are {\rm OE}.\hfill\qed
\end{prop}

\begin{defi}[WOE Cocycles]
\label{defi_WOE_cocycles}%
 Retain the notation  of point~(ii) of
Definition~\ref{defi_WOE}. Due to essential freeness, one can
define measurable cocycles $\alpha:\Gamma\times X\to \Lambda$ and
$\beta:\Lambda\times Y\to \Gamma$ by the a.e.\ requirements
$\alpha(\gamma,x)p(x) = p(\gamma x)$ and $\beta(\lambda,y)q(y) =
q(\lambda y)$. Recall that the cocycle identity reads here
$\alpha(\gamma\gamma',x) =
\alpha(\gamma,\gamma'x)\alpha(\gamma',x)$.
\end{defi}

When two actions are  WOE~-- or even OE~-- the maps which send
orbits into orbits are of course far from being unique. Supposing
for simplicity that the actions on $X,Y$ are OE, one can perturb
an orbit equivalence $F:X\to Y$ by any measurable assignment
$\fhi:X\to \Lambda$, thereby defining $\widetilde F (x) = \fhi(x)
F(x)$, which induces the same OE. It is easy to see that any
isomorphism $\widetilde F$ inducing the same OE is actually
obtained in this way, and that this yields a cohomologous (or
\emph{equivalent}) cocycle $\widetilde \alpha\sim \alpha$. For
later reference we record the following elementary result.

\begin{lemma}
With the above  notation{\rm ,} suppose that the $\Gamma$\ti\ and
$\Lambda$-actions are {\rm OE} and that the associated cocycle
$\alpha:\Gamma\times X\to \Lambda$ is equivalent to a cocycle
$\widetilde\alpha$ which does not depend on $x \in X$. Then the
essential value map $f:\Gamma\to \Lambda$ determined by
$\widetilde \alpha$ is a group isomorphism and the {\rm OE} is induced
by an isomorphism $\widetilde F:X\to Y$ which intertwines the
actions relatively to $f$ \/{\rm (}\/i.e.\ the actions are isomorphic with
respect to $f${\rm )}.\hfill\qed
\end{lemma}

We finally observe  that even if one perturbs an OE map $F:X\to Y$
to obtain $\widetilde F$ as above, the latter will in general
\emph{not} be a bijection and hence {\it a priori\/} not describe an OE.
However it will induce a WOE, and the WOE context is stable under
this operation; hence this setting is more natural and convenient
to work with. The viewpoint of measure equivalence, which we now
turn to, enables us to remove completely the arbitrary choice of
the map $F$ inducing the (weak) orbit equivalence.

Recall from the introduction  (Definition~\ref{defi_ME}) the
definition of an ME coupling $(\Sigma,m)$ between two countable
groups $\Gamma,\Lambda$. We shall say that the ME coupling
$\Sigma$ is \emph{ergodic} if the $\Gamma\times \Lambda$-action
on $\Sigma$ is ergodic; this is equivalent to the ergodicity of
$\Gamma$ on $\Lambda\bsl\Sigma$, or to the ergodicity of $\Lambda$
on $\Gamma\bsl\Sigma$.

Recall that the $\Gamma$-action on $\Sigma$ admits by
definition a measurable fundamental domain $Y\subseteq \Sigma$
with $0<m(Y) <\infty$. Likewise, let $X$ be such a fundamental
domain for $\Lambda$. We shall always endow $\Gamma\bsl\Sigma$
with the measure \emph{restricted} from $m$ {via} the identification
$\Gamma\bsl\Sigma\cong Y$, and likewise for $\Lambda\bsl
\Sigma\cong X$. In order to distinguish from the original
$\Gamma$-action on $\Sigma$, we denote by $\gamma\cdot x$ the
measurable measure-preserving $\Gamma$-action on $X$ obtained
by $\Lambda\bsl \Sigma\cong X$ from the commutativity of the
$\Gamma$\ti\ and $\Lambda$-actions. Likewise, we have also a
``dot'' $\Lambda$-action $\lambda\cdot y$ on $Y$.

\begin{defi}[Retractions, ME Cocycles]
\label{defi_ME_cocycles}%
 Let $\chi:\Sigma\to \Gamma$ be
the measurable $\Gamma$-equivariant map defined by:
$\chi(x)^{-1}x\in Y$ for all $x\in \Sigma$. Then we call $\chi$
the \emph{retraction} associated to $Y$. Likewise, there is a
$\Lambda$-equivariant retraction $\kappa:\Sigma\to \Lambda$
associated to $X$. We obtain thus cocycles $\alpha:\Gamma\times
X\to \Lambda$ and $\beta:\Lambda\times Y\to \Gamma$ (with respect
to the ``dot'' actions) by setting $\alpha(\gamma,x) =
\kappa(\gamma x)^{-1}$ and $\beta(\lambda,y)=\chi(\lambda
y)^{-1}$.
\end{defi}

Thus we have for all $x\in X$ and $\gamma\in \Gamma$ the formula
$$\gamma\cdot x = \alpha(\gamma,x)\gamma x$$
and likewise for $\lambda\cdot y$. Observe also that one can define maps $p_\chi:X\to Y$ and $q_\kappa:Y\to X$ by $p_\chi(x)=\chi(x)^{-1}x$ and $q_\kappa(y)=\kappa^{-1}(y)y$.

\begin{exo}[Trivial Coupling]
\label{exo_trivial}%
 Let $(\Sigma,m)$ be an ME coupling of
$\Gamma$ with $\Lambda$ and assume that both actions on $\Sigma$
are simply transitive (with $m$   purely atomic). Then the choice
of any base point $x\in \Sigma$ defines an isomorphism $f:
\Gamma\to \Lambda$ by taking for $f(\gamma)$ the only
$\lambda\in\Lambda$ such that $\lambda\gamma x=x$. We call
$\Sigma$ a \emph{trivial} coupling and denote it by
$\mathbf{T\!}_f$.

Observe that another choice of $x$ gives a conjugated isomorphism. Observe also that upon identifying $\Lambda$ with $\Sigma$ as the orbit of $x$, the action becomes $(\gamma, \lambda)\,\eta = \lambda\eta f(\gamma)^{-1}$ for $\eta\in\Lambda$.
\end{exo}

A less trivial (but still very straightforward) source of examples is the following:

\begin{exo}[Lattices] 
\label{exo_lattices}%
Let $G$ be a locally  compact, second countable
group and $\Gamma,\Lambda$ two lattices in $G$. The existence of
lattices implies that any Haar measure $m$ is left and right
invariant; therefore, we obtain an ME coupling $\Sigma=(G,m)$ of
$\Gamma$ with $\Lambda$ by considering the $\Gamma\times
\Lambda$-action given by $(\gamma,\lambda)g = \gamma
g\lambda^{-1}$.
A very special case occurs when $G=\Gamma$ and $\Lambda$ is a finite index subgroup of $\Gamma$.
\end{exo}

Given an ME coupling  $(\Sigma,M)$ of $\Gamma$ with $\Lambda$, we
shall need the following concept which may seem pedantic at first
sight, but will turn out to be extremely useful: Since $\Sigma$ is
technically a $\Gamma\times \Lambda$-space, we may define the
\emph{opposite coupling} $\check\Sigma$ of $\Lambda$ with $\Gamma$
to be the $\Lambda\times \Gamma$-space obtained  {via} the
canonical isomorphism $\Lambda\times \Gamma\cong\Gamma\times
\Lambda$. As this will be particularly relevant in situations
where $\Lambda=\Gamma$, we will (though rarely!) have to
distinguish the $\Gamma$-actions on $\Sigma$ by writing
$(\gamma,x)\mapsto A^1_\gamma x$ and $A^2_\gamma x$ respectively
(then $\check\Sigma$ is obtained by switching $A^1$ and $A^2$).

\begin{defi}[Coupling Composition]
\label{defi_composition}%
 Assume we are furthermore\break given an ME coupling $(\Omega,n)$ of $\Lambda$
with a third (countable) group $\Delta$. Define the \emph{composed coupling}
$\Sigma\times_\Lambda \Omega$ to be the quotient space of
$\Sigma\times \Omega$ by the product $\Lambda$-action. By
commutativity, this is still a $\Gamma\times \Delta$-space, and
we turn it into an ME coupling of $\Gamma$ with $\Delta$ by
endowing it with the measure obtained by restricting $m\otimes n$
to an (infinite measure) fundamental domain for $\Lambda$ in
$\Sigma\times \Omega$.
\end{defi}

\begin{defi}[Coupling Index] 
Given an ME coupling $(\Sigma,m)$ of
$\Gamma$ with $\Lambda$, define its \emph{coupling index} to be
the following positive number:
$$[\Gamma:\Lambda]_\Sigma = \frac{m(\Lambda\bsl\Sigma)}{m(\Gamma\bsl\Sigma)}.$$
\end{defi}

The notation reflects the fact that in the particular case where $\Lambda$ is a finite index subgroup of $\Gamma$ (Example~\ref{exo_lattices}) we recover indeed the index $[\Gamma:\Lambda] = \big|\Gamma/\Lambda\big|$. More generally, the coupling index corresponds to the ratio of co-volumes if $\Gamma,\Lambda$ are lattices in one given locally compact group. It is straightforward to verify the formulae
\begin{equation}
\label{eq_constant_functorial}%
[\Gamma:\Lambda]_\Sigma = 1/[\Lambda:\Gamma]_{\check\Sigma},\kern.75cm [\Gamma:\Delta]_{\Sigma\times_\Lambda \Omega} = [\Gamma:\Lambda]_\Sigma \cdot [\Lambda:\Delta]_\Omega.
\end{equation}
We need one more

\begin{exo}[Standard Coupling]
\label{exo_standard}%
 Let $\Gamma$ be a countable group and\break
$(X,\mu)$ a probability $\Gamma$-space. We define an ME coupling
of $\Gamma$ with itself as follows: Endow $\Sigma = X\times
\Gamma$ with the product measure and define the\break $\Gamma$-actions
$A^1,A^2$ by $A^1_\gamma(x,\gamma_0) = (\gamma x,\gamma\gamma_0)$
and $A^2_\gamma(x,\gamma_0) = (x,\gamma_0\gamma^{-1})$; we call
this the \emph{standard} coupling associated to $X$. The two
resulting $\Gamma$-actions on $A^1(\Gamma)\bsl\Sigma$ and
$A^2(\Gamma)\bsl\Sigma$ are both isomorphic to the $\Gamma$-action on $X$. The subset $X\times\{e\}\subseteq \Sigma$ is a
common fundamental domain for both actions on $\Sigma$, the
associated cocycles are the identity isomorphism and there is a
natural quotient map $\Sigma\to \mathbf{T\!}_{\id}$ to the trivial
coupling whose fibres can be identified with $X$.

Conversely,  it is easy to verify that every ME coupling
satisfying the properties listed above is measurably isomorphic to
a standard coupling $\Sigma$ as above.
\end{exo}

We now state the fundamental observation concerning the relation
between ME and WOE. The following is proved by
A.~Furman~\cite{Furman2} (who gives credit also to M.~Gromov and
R.~Zimmer).

\begin{thm}[ME-WOE] 
\label{thm_OE_WOE}%
Let $\Gamma,\Lambda$ be countable groups and
$(X_0,\mu)${\rm ,} $(Y_0,\nu)$ be probability $\Gamma$\ti\ and
$\Lambda$-spaces respectively. To any {\rm WOE} given with $p,q$ as
in Definition~{\rm \ref{defi_WOE}} point~{\rm (ii)} corresponds an {\rm ME} coupling
$\Sigma$ of $\Gamma$ with $\Lambda${\rm ,} together with a choice of
$\Gamma$\ti\ and $\Lambda$-fundamental domains $Y,X$ resp.{\rm ,}
such that\/{\rm :}\/

\begin{itemize}
\item[{\rm (i)}] Modulo renormalisation of measures{\rm ,} one has isomorphisms
 of $\Gamma$-spaces $X_0\cong \Lambda\bsl\Sigma\cong X$ and of
  $\Lambda$-spaces $Y_0\cong \Gamma\bsl\Sigma\cong Y$.

\item[{\rm (ii)}] Under these identifications{\rm ,}
  $p_\chi=p$ and $q_\kappa=q$. Moreover{\rm ,} the {\rm WOE} cocycles $\alpha,\beta$ of
Definition~{\rm \ref{defi_WOE_cocycles}} coincide with the {\rm ME}
 cocycles of Definition~{\rm \ref{defi_ME_cocycles}.}

\item[{\rm (iii)}] $C(X_0,Y_0) = [\Gamma:\Lambda]_\Sigma$.
\end{itemize}
  Moreover{\rm ,} in the ergodic case{\rm ,}  
$[\Gamma:\Lambda]_\Sigma=1$ if and only if the {\rm WOE} is \/{\rm (}\/induced
by\/{\rm )}\/ an {\rm OE,} and then one can choose in $\Sigma$ a common
$\Gamma$\ti\ and $\Lambda$-fundamental domain.

Conversely{\rm ,} the above procedure produces {\rm WOE} probability
$\Gamma$-respectively $\Lambda$-spaces out of any {\rm ME}
coupling $\Sigma$ and the above three properties hold. \/{\rm (}\/Yet{\rm ,} in
contrast to our standing assumption these spaces need not be
essentially free {\rm --} see Remark~{\rm \ref{rem_essentially_free}} below\/{\rm .)}\/
\end{thm}

{\it On the proof}.
See~3.2 and~3.3 in Furman~\cite{Furman2} where ergodicity is assumed. However one can reduce to this case 
by~\cite[2.2]{Furman1}.
\Endproof

We can now see what \emph{commensurability} for actions should be.

\begin{exo}[Stability Properties]
\label{exo_commensurability}%
 Here are two constructions that
appear naturally and will be useful in the sequel; they are in a
sense mutually dual.

\smallskip

 {\bf(i).} Consider a countable  group $\Gamma$ and a
probability $\Gamma$-space $(X,\mu)$. Let $N\lhd \Gamma$ be a
finite normal subgroup and set $\Lambda = \Gamma/N$. Consider the
quotient $N\bsl X$ (with quotient measure) as a probability
$\Lambda$-space. Then the $\Gamma$-action on $X$ is WOE to
the $\Lambda$-action on $N\bsl X$ since one can take for
Definition~\ref{defi_WOE}~(ii) $p:X\to N\bsl X$ to be the quotient
map and $q:N\bsl X\to X$ any measurable cross-section.
Alternatively, one meets the other condition of that definition by
taking $A=q(N\bsl X)$ for $q$ as before and $B=N\bsl X$ (thus the
compression constant is $C(X, N\bsl X)=|N|$). The ME coupling
associated to this WOE is the following: First let $\Sigma$ be the
standard coupling of $\Gamma$ with itself associated to $X$ as in
Example~\ref{exo_standard}, and then consider the coupling $N\bsl
\Sigma = A^2(N)\bsl\Sigma$ obtained from $\Sigma$ by dividing out,
say, the second $N$-action. Of course we have
$[\Gamma:\Lambda]_{N\bsl\Sigma} = |N| = C(X, N\bsl X)$.

\smallskip

 {\bf(ii).} This time we consider a finite  index
subgroup $\Lambda$ of a countable group $\Gamma$ and a probability
$\Lambda$-space $(Y,\nu)$. We write
$Y\!\!\uparrow_{\Lambda}^\Gamma$ for the $\Gamma$-space which
is the \emph{suspension} (or \emph{induction}) of the
$\Lambda$-action on $Y$; this space is obtained (after dividing
the measure by $[\Gamma:\Lambda]$) by considering the quotient of
$\Gamma\times Y$ by the $\Lambda$-action
$\lambda(\gamma,y)=(\lambda \gamma,\lambda y)$ endowed with the
$\Gamma$-action descending from
$\gamma_1(\gamma_2,y)=(\gamma_2\gamma_1^{-1},y)$. Then the
$\Gamma$-action on $Y\!\!\uparrow_{\Lambda}^\Gamma$ is WOE to
the $\Lambda$-action on $Y$. Indeed, the first equivalent
characterisation in Definition~\ref{defi_WOE} is met by setting
$B=Y$ and letting $A\subseteq Y\!\!\uparrow_{\Lambda}^\Gamma$ be
the image of $\{e\}\times Y$. Alternatively, for the second
characterisation, let $p: Y\!\!\uparrow_{\Lambda}^\Gamma\to Y$ be
the natural quotient map and $q$ be the section obtained by
$y\mapsto (e,y)$. In particular, we have the compression constant
$C(Y\!\!\uparrow_{\Lambda}^\Gamma, Y)=[\Gamma:\Lambda]$. To
describe the ME coupling associated to this WOE, one considers
again the standard coupling $\Sigma$ (of $\Lambda$ this time)
associated to $Y$; then, either one takes the suspension of, say,
the first $\Lambda$-action on $\Sigma$, or~-- equivalently~--
one composes the coupling $\Sigma$ with the coupling arising from
the inclusion $\Lambda<\Gamma$ as in the end of
Example~\ref{exo_lattices}.

We conclude the example by remarking that if we have a probability\break
$\Gamma$-space $(X,\mu)$ and a finite index subgroup
$\Lambda<\Gamma$, then in general the restricted $\Lambda$-action on $X$ will {\it not} be WOE to the
original $\Gamma$-action. This can be seen for instance as follows: Suppose
$\Gamma=\Gamma_1\times \Gamma_2$, where the $\Gamma_i$'s are
torsion-free and in~$\creg$ (e.g.\ non-Abelian free groups). Let
$\Lambda<\Gamma$ be a finite index subgroup not isomorphic to
$\Gamma$. Now if $X$ is mildly mixing $\Gamma$-space, then the
$\Lambda$-action on $X$ cannot be WOE to the $\Gamma$-action, since that would contradict the generalisation of
Theorem~\ref{thm_superrigidity0} given below as
Theorem~\ref{thm_factors_irred}~(ii). Notice however, that this stands in contrast to,
but does not contradicts
the fact that any ME coupling of $\Gamma$ with some other
countable group $\Delta$ also forms an ME coupling of the finite index
subgroup $\Lambda<\Gamma$ with $\Delta$.
\end{exo}

\begin{rem}
\label{rem_essentially_free}%
There is some lack of symmetry in the relation between ME and WOE,
because the  WOE actions on probability spaces obtained as
quotients of an ME coupling can be far from being free (consider
e.g.\ Example~\ref{exo_lattices} with  Abelian $G$, or the trivial
coupling for which the quotients reduce to a point). As far as
proofs are concerned, this is not a difficulty for us, as we
establish all our proofs in the setting of ME couplings and then
deduce the WOE or OE statements, thereby using only the WOE $\lra$
ME direction. However, since the opposite direction will be useful
to us when constructing some examples, we observe that the
technicality arising in the inverse construction can easily be
circumvented. This is achieved by composing a given ME coupling
$\Sigma$ (with potentially nonfree $\Gamma$\ti\ or $\Lambda$-quotients) with a standard self-coupling of $\Gamma$
(Example~\ref{exo_standard}) associated to any free probability
$\Gamma$-space~$X$. The relevant properties of $\Sigma$ will be
preserved in the composed coupling; ergodicity properties, such as
irreducibility, are preserved if one chooses $X$ to be
``sufficiently ergodic'' (e.g.\ mixing) $\Gamma$-space.
\end{rem}
 
\vglue-8pt
\Subsec{Reformulation and discussion of the main results}  \label{sec_reformulation} 
The relation between OE and ME, as discussed in the preceding
subsection, enables us to reformulate our results in terms of the
latter notion. In doing so we shall also generalise the main
results to the framework of weak orbit equivalence.

As mentioned in the introduction, we will consider a family of
groups more general than the class~$\creg$. The property relevant
to our approach is described by the following:

\vglue-18pt
\phantom{up}
\begin{defi}[Class~$\cc$]
\label{defi_C}%
Denote by~$\cc$ the class of  groups admitting a
mixing unitary representation $\pi$ on a separable Hilbert space,
such that $\hb^2(\Gamma,\pi)\neq 0$.
\end{defi}

Recall that a unitary representation is called \textit{mixing} if
all its matrix coefficients vanish at infinity; the outstanding
example, and the one we shall actually use, is the regular
representation. It follows from Theorem~\ref{thm_exos_c} that the Examples~\ref{exos_c} introduced for the
sake of concreteness are all contained in~$\cc$ (see
Section~\ref{sec_in_c} which has more on~$\cc$).

We next extend Definition~\ref{defi_Irreducible} above in order
to cover products of any number of groups:

\begin{defi}
Let $\Gamma_1, \ldots, \Gamma_n$ be groups and  set
$\Gamma=\Gamma_1\times \cdots\times \Gamma_n$. A $\Gamma$-space
$(X,\mu)$ is called \emph{irreducible} if for every $1 \le j \le
n$ the subproduct $\Gamma'_j=\prod_{i\neq j}\Gamma_i$ acts
ergodically on $X$.
\end{defi}

Notice that this definition forces $n>1$ (unless $X$ is trivial).

 \vglue4pt

We begin reformulating our main results by considering
Theorem~\ref{thm_strong_rigidity1}. We shall in fact prove the
following more general version of it:

\begin{thm}
\label{thm_factors_irred}%
Let $\Gamma_1,\ldots, \Gamma_n$ be torsion-free 
groups in~$\cc$ and $(\Sigma, m)$ be an {\rm ME} coupling of
$\Gamma=\Gamma_1\times \cdots\times\Gamma_n$ with a product
$\Lambda=\Lambda_1\times\cdots \times\Lambda_n$ of any
torsion-free countable groups such that the
$\Lambda$-action on $\Gamma\bsl\Sigma$ is irreducible. Assume
that either
\begin{mitemize}
\item[{\rm (i)}] $[\Lambda:\Gamma]_\Sigma\geq 1${\rm ,} or 
\item[{\rm (ii)}] the $\Gamma$-action on $\Lambda\bsl\Sigma$ is irreducible.
\end{mitemize}
 Then{\rm ,} upon  permuting indices{\rm ,} there are
isomorphisms $f_i:\Lambda_i\xrightarrow{\ \cong\ } \Gamma_i$ such
that identifying $\Gamma$ with $\Lambda$ through $f=\prod
f_i:\Lambda\cong\Gamma${\rm ,} the coupling $\Sigma$ is a standard
coupling.
  Equivalently{\rm ,} by reference to Example~{\rm \ref{exo_standard}}
for the notion of standard couplings{\rm ,}  
$[\Lambda:\Gamma]_\Sigma= 1$ and there is a common fundamental
domain $Y\subseteq \Sigma$ for both actions such that $\lambda Y =
f(\lambda)Y$ for all $\lambda\in\Lambda$.
\end{thm}

Thus, at the level of (W)OE, Theorem~\ref{thm_factors_irred}
implies that, under the assumptions corresponding to the above,
\itshape any {\rm WOE} of the actions is in fact an {\rm OE} induced by an
isomorphism of the actions with respect to an isomorphism  of the
groups\upshape.

Therefore, Theorem~\ref{thm_strong_rigidity1} follows from Theorem~\ref{thm_factors_irred} in the particular case $\Lambda_i=\Gamma_i$, $n=2$, $\creg$ instead of~$\cc$ and essentially free quotients.

We give now an illustration of the necessity of the assumptions in Theorem~\ref{thm_factors_irred}:

\begin{exo}[Coupling Index Condition]
\label{exo_free_free}%
 Let $F_n$ denote the  free
group on $n$ generators. Realise $F_3$ and $F_5$ as index-two
subgroups of $F_2$ and $F_3$ respectively, and view $F_2\times F_5
<F_2\times F_3$ and $F_3\times F_3 <F_2\times F_3$ as index-two
subgroups. Thus $F_2\times F_3$ is an ME coupling of
$\Gamma=F_2\times F_5$ with $\Lambda=F_3\times F_3$ and the
coupling index is one (observe that there is indeed a common
$\Gamma$\ti\ and $\Lambda$-fundamental domain $\{(e,e),(x,y)\}$
in $F_2\times F_3$ given by representatives $x,y$ of the
nontrivial cosets in $F_2,F_3$). Thus we see that the coupling
index condition alone is not sufficient to derive the conclusion
of Theorem~\ref{thm_factors_irred}. On the other hand, if we
replace $F_5$ by $F_3$ then $F_2\times F_3$ becomes a coupling of
$F_2\times F_3$ with $\Lambda=F_3\times F_3$ for which the latter
acts irreducibly (indeed, it acts on a point). However this time
the inequality for the coupling constant is not satisfied, which
accounts for the failure of the conclusion of the theorem. (Recall
from Remark~\ref{rem_essentially_free} that one can also build OE
and WOE counter-examples out of the ME examples given here upon
making the quotient actions essentially free by composition with,
say, the standard coupling associated to a Bernoulli shift.)
\end{exo}

Here is now an example showing how the statement breaks down for groups not in~$\cc$.

\begin{exo}
[Class $\cc$ Condition, I] Let $G$ be a connected
noncompact simple Lie group with trivial center and consider four
copies of $G$ labeled $G_i$, $1\leq i\leq 4$. For each pair $1\leq
i\neq j\leq 4$ let $\Gamma_{ij}$ be an irreducible lattice in the
product $G_i\times G_j$. Of course, one may choose none, some, or
all $\Gamma_{ij}$ to be nonisomorphic as abstract groups. Now
$\Gamma=\Gamma_{12}\times\Gamma_{34}$ as well as
$\widetilde\Gamma=\Gamma_{13}\times\Gamma_{24}$ can both be
realised naturally as lattices in $\prod_{1\leq i\leq 4}G_i$, thus
producing an ME coupling of $\Gamma$ with $\widetilde \Gamma$.
Using Howe-Moore's theorem, it is easy to check that this coupling
is irreducible~-- namely each $\Gamma_{ij}$ acts ergodically on
the quotient of $G^4$ by the``other'' product. Moreover, the
conclusion of the theorem fails even if we take all $\Gamma_{ij}$
isomorphic, as this ME coupling is not a standard one
(Example~\ref{exo_standard}). This also shows that a nontrivial
assumption on the groups is needed in
Theorem~\ref{thm_strong_rigidity1} from the introduction.
\end{exo}

Next, we reformulate and generalise
Theorem~\ref{thm_superrigidity1}  (and thus
Theorem~\ref{thm_superrigidity0}) in the ME setting and discuss
the assumptions made there. For the sake of clarity, we separate the
statements for the groups~(\ref{thm_factors_mystery_new}) and for
the actions~(\ref{thm_factors_mystery_new}*):

\begin{thm}
\label{thm_factors_mystery_new}%
Let $\Gamma_1, \ldots, \Gamma_n$ be
torsion-free groups in~$\cc$ and let $\Lambda$ be any countable
group admitting an {\rm ME} coupling $(\Sigma,m)$ to
$\Gamma=\Gamma_1\times \cdots\times \Gamma_n$.

If the $\Gamma$-action on $\Lambda\bsl\Sigma$ is irreducible 
and the $\Lambda$-action on $\Gamma\bsl \Sigma$ is mildly mixing{\rm ,} then $\Lambda$ fits in an extension
\begin{equation}
\label{eq_extension}%
1\lra N\lra \Lambda \xrightarrow{\ \pi\ } \Gamma' \lra 1
\end{equation}
where $N$ is finite and $\Gamma'<\Gamma$ 
is a finite index subgroup whose projections to each $\Gamma_i$ are onto. Moreover{\rm ,}  
\begin{equation}
\label{eq_index_formula}%
[\Gamma:\Gamma']\ =\ |N|\cdot [\Gamma:\Lambda]_\Sigma.
\end{equation}
\end{thm}

Let us call an exact sequence~(\ref{eq_extension}) with $N$ finite and $\Gamma'$ of finite index in $\Gamma$ a \emph{virtual isomorphism} of the groups $\Lambda,\Gamma$. We have seen above (Example~\ref{exo_commensurability}) that in this setting the natural generalisation of isomorphic actions is a sort of commensurability of actions; we show in the proof of Theorem~\ref{thm_factors_mystery_new} that any WOE of actions as in the setting of that theorem are in fact a virtual isomorphism; more precisely:

\demo{\scshape Theorem {\rm  \ref{thm_factors_mystery_new}}$^\ast$} 
{\it Let $\Gamma_1,\ldots, \Gamma_n$ 
be torsion-free groups in~$\cc$ with an irreducible essentially free $\Gamma=\Gamma_1\times \cdots\times
\Gamma_n$-action on a probability space $Y$.

If this action is {\rm WOE}
 to a mildly  mixing{\rm ,} essentially free{\rm ,} action of any countable group $\Lambda$ on a probability space
$X${\rm ,} then there is  a virtual isomorphism as in~{\rm (\ref{eq_extension})}
 and the corresponding $\Gamma$-action on 
$(N\bsl X)\!\uparrow_{\Gamma'}^\Gamma$ is isomorphic to $Y$.}
\Enddemo

In the OE case we can deduce a stronger statement upon assuming aperiodicity of $\Gamma$:

\begin{cor}
\label{cor_OE_mystery}%
Let $\Gamma_1,\ldots, \Gamma_n$ be torsion-free
groups  in~$\cc${\rm ,} and let $Y$ be  an aperiodic irreducible
essentially free $\Gamma=\Gamma_1\times \cdots\times \Gamma_n$-space.

\nobreak If this action is {\rm OE} to a mildly mixing{\rm ,} essentially free
action of any countable group $\Lambda$ on  a probability space
$X${\rm ,} then there exists an isomorphism of $\Lambda$ and  $\Gamma$
with respect to which  the actions on $X,Y$ are isomorphic.
\end{cor}

(Observe that aperiodicity and irreducibility both hold if e.g.\ the $\Gamma$-action is mildly mixing.)

In the light of the discussion of Section~\ref{sec_ME_WOE}, the above result imply indeed Theorems~\ref{thm_superrigidity0} and~\ref{thm_superrigidity1} stated in
the introduction. In fact, we see that if we assume only that the actions  in Theorem~\ref{thm_superrigidity1}
are WOE, we still obtain both conclusions~(i) and~(ii), with the modified formula $[\Gamma:\Gamma']\ =\
|N|\cdot C(X,Y)$.

Note that  the OE assumption is equivalent to $C(X,Y)=1$; that $N$
is trivial as soon as $\Lambda$ is torsion-free; and that on the
other hand aperiodicity forces $[\Gamma:\Gamma']=1$ because in
that case the action cannot be a suspension of an action of a
proper finite index subgroup. This accounts for
Theorem~\ref{thm_superrigidity0} and
Corollary~\ref{cor_OE_mystery}.

\begin{exo}[Mild Mixing Condition]
\label{exo_mild_mixing}
 In order to put the mild mixing
assumption in a better perspective, consider the following
situation. Let $G$ be a connected, rank one, simple Lie group with
trivial center (e.g.\ $\mathrm{PSL}_2({\bf R})$) and choose two
lattices $\Gamma_1, \Gamma_2<G$ (in particular, the $\Gamma_i$'s
are in~$\creg$). Let $\Lambda<G\times G$ be an irreducible
lattice. Then, as in Example~\ref{exo_lattices}, $G \times G$ is an ME coupling of $\Gamma=\Gamma_1\times \Gamma_2$ with $\Lambda$, or
equivalently, the $\Gamma$-action on $G^2/\Lambda$ is WOE to
the $\Lambda$-action on $G^2/\Gamma$ (the essential freeness of
these actions can be deduced from the center freeness of $G$). One
can arrange to have the same co-volumes, so that then the actions
are in fact OE. Furthermore, the irreducibility of the lattice
$\Lambda$ ensures that the $\Gamma$-action is irreducible.
Thus, we have here a situation where the conclusion of
Theorem~\ref{thm_factors_mystery_new} (and
Theorem~\ref{thm_superrigidity0}) fails because the\break $\Lambda$-action on $G^2/\Gamma$ is not
mildly mixing, even though it does have very strong ergodicity properties: It is weakly mixing, and
moreover one can find $\Lambda$ for which every nontrivial
element acts ergodically (or weakly mixing). In fact, one can
detect precisely how the mild mixing property
in Definition~\ref{defi_mild_mixing} fails: By Howe-Moore's
theorem, it can be shown that the only nontrivial recurrent sets
are of the form $A\times G/\Gamma_2$ or $G/\Gamma_1\times B$, and
the associated recurrent sequences $(\lambda_n)$ of $\Lambda$ must
satisfy $\mathrm{pr}_1(\lambda_n)\to e_1$ or
$\mathrm{pr}_2(\lambda_n)\to e_2$, respectively, where
$\mathrm{pr}_i$ is the $i^{\rm th}$ quotient map $G \times G\to G$ and
$e_i$ the trivial element in the $i^{\rm th}$ factor of the two.
\end{exo}

Analogous to the case of products, we restate and generalise the rigidity results for groups with amenable radicals through the notion of measure equivalence:

\begin{thm}
\label{thm_kernel_irred}%
Let $\Gamma,\Lambda$ be countable groups and $M\lhd \Gamma, N\lhd
\Lambda$  amenable normal subgroups such that
$\overline{\Gamma}=\Gamma/M$ and $\overline{\Lambda}=\Lambda/N$
are in~$\cc$ and torsion-free. Let $(\Sigma,
m)$ be an {\rm ME} coupling of $\Gamma$ with $\Lambda$.

\nobreak If $N$ is ergodic on $\Gamma\bsl\Sigma$ and $M$ on
$\Lambda\bsl\Sigma${\rm ,}  then there is an isomorphism $f:
\overline{\Gamma}\xrightarrow{\ \cong\ } \overline{\Lambda}$.
 Moreover{\rm ,} $\Sigma$ admits a $\Gamma\times \Lambda$-equivariant factor $\Phi:\Sigma\to
\mathbf{T\!}_f${\rm ,}  where the latter is the trivial coupling of $\overline{\Gamma}$ with
$\overline{\Lambda}$ inducing $f$.
\end{thm}

The superrigidity-type statement goes as follows:

\begin{thm}
\label{thm_kernel_mystery_new}%
Let $\Gamma$ be a countable group with an amenable normal
subgroup $M\lhd \Gamma$ such that $\overline{\Gamma}=\Gamma/M$ is
in~$\cc$ and torsion-free{\rm ,} and let $\Lambda$ be
any countable group with an ME coupling $(\Sigma,m)$ to $\Gamma$.
 If the $M$-action on $\Lambda\bsl\Sigma$ is ergodic and the $\Lambda$-action on $\Gamma\bsl \Sigma$
is mildly mixing{\rm ,} then there is an amenable normal subgroup $N\lhd\Lambda$ such that
$\overline{\Lambda}=\Lambda/N$ is isomorphic to $\overline{\Gamma}$.
\end{thm}

To verify that these results indeed imply Theorems~\ref{thm_strong_rigidity2} and~\ref{thm_superrigidity2} stated in
the introduction, one appeals again to Theorem~\ref{thm_OE_WOE} above.
 
We now discuss some situations related to Theorem~\ref{thm_factors}.

\begin{exo}
[Class $\cc$ Condition, II] Let $G$ be any discrete group
with Kazhdan's property~(T) and $H$ be a group without
property~(T). Set $\Gamma_1 = G\times G$, $\Gamma_2 = H\times H$,
$\Lambda_1 = G\times H$ and $\Lambda_2 = H\times G$. Then $\Gamma
= \Gamma_1\times \Gamma_2$ is ME (indeed isomorphic) to $\Lambda =
\Lambda_1\times \Lambda_2$; however, $\Gamma_1$ is not ME to any
$\Lambda_i$ since property~(T) is an ME
invariant~\cite[1.4]{Furman1}. Thus, some nontrivial assumption
on the groups in Theorem~\ref{thm_factors} is necessary. In fact,
we do not have any natural candidate for a more general class of
groups than~$\cc$ for which a similar result should hold.

The fact that groups in~$\cc$ cannot have infinite  direct factors
(see Section~\ref{sec_in_c}) is illustrated in a patent way in the
above example. Indeed, we may arrange for both $G$ and $H$ to be
in~$\cc$ or even in~$\creg$: Take for instance for $G$ a lattice
in $\mathrm{Sp}(n,1)$ with $n\geq 2$ and for $H$ a free group on
two generators. Then, as above, the conclusion of
Theorem~\ref{thm_factors} fails for the products
$\Gamma=\Gamma_1\times \Gamma_2$, $\Lambda=\Lambda_1\times
\Lambda_2$,
 but of course after further splitting of the factors
one can shuffle the groups to get the (trivial) self-couplings of
$G$ and of $H$ respectively, in accordance with the theorem.
\end{exo}

As another example,  consider the construction described in
Example~\ref{exo_free_free}, namely the ME coupling of
$\Gamma=F_2\times F_5$ and $\Lambda=F_3\times F_3$ with coupling
index one~-- so that these groups admit actions (which can be made
free) that are indeed OE and not just WOE. By the recent result of
D.~Gaboriau~\cite{GaboriauL2}, the $\ell^2$-Betti numbers are OE
invariants, so that neither the couple $F_2$ and $F_3$, nor the
couple  $F_5$ and $F_3$, admit OE actions. Thus, even by assuming
that two products $\Gamma=\Gamma_1\times \Gamma_2$ and
$\Lambda=\Lambda_1\times \Lambda_2$ admit OE actions, one cannot
arrive at a stronger conclusion in Theorem~\ref{thm_factors}. The
reader is invited to examine the proof of
Theorem~\ref{thm_factors} in this very simple and concrete example
to see how the equality of co-volumes can be lost in passing from
the original ME coupling to couplings of the individual factors.

\begin{rem}
\label{rem_Betti}%
Suppose that a product $\Gamma = \Gamma_1 \times \dots \times
\Gamma_n$ of groups in~$\creg$ is ME to a
torsion-free product $\Lambda = \Lambda_1 \times \dots \times
\Lambda_n$ of any (countable) groups $\Lambda_i$. Gaboriau's
results~\cite{GaboriauCRAS}, \cite{GaboriauL2} imply that the
$\ell^2$-Betti numbers $\beta_{(2)}^i$ of $\Gamma$ are
proportional to those of $\Lambda$. But now
Theorem~\ref{thm_factors} tells us that (after permutation of
indices) we can also apply this to each pair, giving of course
more restrictions on the possible values of the $\ell^2$-Betti numbers of the factors.

For instance, suppose we have a group $\Gamma_1$ in~$\creg$ with the nonzero $\ell^2$-Betti numbers $\beta_{(2)}^2=2$, $\beta_{(2)}^3=3$, $\beta_{(2)}^4=1$ and set $\Gamma_2=\Gamma_1$. Choose now any torsion-free countable groups $\Lambda_i$ such that:
{\small\begin{center}
\noindent\begin{tabular}{c|c|c|c|c|c|}
 & $\beta_{(2)}^1$ & $\beta_{(2)}^2$ & $\beta_{(2)}^3$ & $\beta_{(2)}^4$ & $\beta_{(2)}^{\geq 5}$\\
\hline
$\Lambda_1$ & $0$ & $1$ & $2$ & $1$ & $0$\\
\hline
$\Lambda_2$ & $0$ & $4$ & $4$ & $1$ & $0$\\
\hline
\end{tabular}
\end{center}}
\noindent
Then $\Gamma = \Gamma_1\times \Gamma_2$ has the same $\ell^2$-Betti numbers as $\Lambda=\Lambda_1\times \Lambda_2$, so that Gaboriau's result does not exclude an ME coupling of these two groups. However, such a coupling is impossible in view of Theorem~\ref{thm_factors} since then we would have individual couplings, and \emph{that} would now contradict Gaboriau's 
proportionality.
\end{rem}

Finally, we make some  concluding remarks on the irreducibility
property in Theorem~\ref{thm_factors_irred}. Suppose that we have
an ME coupling $\Sigma$ between two groups $\Gamma=\Gamma_1\times
\Gamma_2$ and $\Lambda_1\times \Lambda_2$, where all four factors
are torsion-free and in~$\creg$ but the $\Gamma_i$ are not
isomorphic to the $\Lambda_j$. (For instance, the $\Gamma_i$ are
countable non-Abelian free groups and the $\Lambda_j$ surface
groups of genus $\geq 2$.)
Then Theorem~\ref{thm_factors} tells us that (upon permuting indices) $\Gamma_1$ is ME to $\Lambda_1$ and $\Gamma_2$ to $\Lambda_2$; but on the other hand, the coupling $\Sigma$ cannot be irreducible for both $\Gamma$ and $\Lambda$ because of Theorem~\ref{thm_factors_irred}~(ii). Can one deduce in certain situations that $\Sigma$ is actually a \emph{product coupling}? Likewise, if $[\Gamma:\Lambda]_\Sigma=1$, the coupling cannot even be irreducible for one side in view of Theorem~\ref{thm_factors_irred}~(i); so, again, must it be a product coupling?

\Subsec{Groups with many actions} 
\label{sec_continuum} 
We begin by proving an observation made in the introduction:

\begin{thm}
\label{thm_free_continuum}%
Let $\Gamma$ be a \/{\rm (}\/countable\/{\rm )}\/ free group. Then any given
probability $\Gamma$-space $(X,\mu)$ is {\rm OE} to actions of
uncountably many nonisomorphic groups.
\end{thm}

{\it Proof}.
We may assume that $\Gamma$ has rank at least two in view of the 
Ornstein-Weiss result~\cite{Ornstein-Weiss} for amenable groups.
In fact, for simplicity of notation only we shall take $\Gamma$ of
rank two. What we shall actually show is that for \emph{every}
pair of countable amenable groups $A,B$ the $\Gamma$-action on
$X$ is OE to an action of $\Lambda=A*B$ on the same space $X$. Let
$u,v$ be free generators of $\Gamma$; to avoid technical issues,
assume that both $u$ and $v$ are ergodic transformations of $X$
(it is not difficult to remove this assumption, keeping the same
strategy of proof). Consider the infinite (cyclic) amenable groups
$\langle u \rangle$ and $\langle v \rangle$ and note that by the result of Ornstein-Weiss,
there exist measure-preserving, essentially free actions of $A$ and
$B$ on $(X,\mu)$, each of which has a.e.\ the same orbits as
$\langle u \rangle$ and $\langle v \rangle$ respectively. These actions of $\langle u \rangle$ and $\langle v \rangle$
define by universality an action of their free product $\Lambda$,
which has the same orbits as $\Gamma$. The whole point of the
argument is to show that this action is essentially free. Indeed,
otherwise there is a nontrivial element $a_1b_1\cdots a_nb_n$ of
$A*B$ which fixes pointwise a measurable set $Y\subseteq X$ with
$\mu(Y)>0$. Now for a.e.\ $y\in Y$ there are integers
$p_1,\ldots,p_n$ and $q_1, \ldots, q_n$ such that
$u^{p_1}v^{q_1}\cdots u^{p_n}v^{q_n}y =a_1b_1\cdots a_nb_ny = y$.
Since there are countably many $n$-tuples $(p_i,q_i)$, this
contradicts the essential freeness of the $\Gamma$-action.
\hfill\qed

\begin{rem}
More generally,  it seems that whenever $G,H,A,B$ are countable
groups such that $G$ and $A$ admit OE actions, and likewise for $H$
and $B$, then $G*H$ admits an action OE to an action of $A*B$ (and
in particular $G*H$ is ME to $A*B$). This should follow from a
similar idea, realising the OE for $G$ and $A$ on a common space
$X$ and the OE for $H$ and $B$ on a space $Y$, only that now one
has to choose an isomorphism of standard probability spaces
$X\cong Y$ such that the resulting actions of the free products
are essentially free~-- e.g.\ by applying the Baire category
theorem to the Polish space of such isomorphisms. Observe that
this line of reasoning does not pass to WOE (this is not possible
in general, as the example of finite groups of different order
shows). The situation is reminiscent of the known difference between
bi-Lipschitz and quasi-isometric equivalence for free products of
finitely generated groups.
\end{rem}
\phantom{up}
\vskip-20pt

The above result stands in strong contrast to our Theorem~\ref{thm_continuum} from the introduction; let
us turn to the latter.

\vskip4pt
{\it Proof  of Theorem~{\rm \ref{thm_continuum}}}.
The idea is to apply  our Theorem~\ref{thm_strong_rigidity1} in
order to get many actions of a given group that are mutually not
OE; but actually, we shall rather use the stronger statement of
Theorem~\ref{thm_factors_irred} in order to be able to vary the
groups as well. That way, we shall construct a family of actions
as claimed in Theorem~\ref{thm_continuum}, but furthermore
\itshape no two of them will even be\upshape WOE.

Let $\mathcal  F$ be the continuum of isomorphism classes of all
groups $\Gamma=\Gamma_1\times \Gamma_2$, where $\Gamma_i=A*B$
range over all free products of any two
torsion-free countable groups. In view of Theorem~\ref{thm_factors_irred},
all we have to do is to find for each such $\Gamma$ in $\mathcal
F$ a continuum of \emph{nonisomorphic} irreducible probability
$\Gamma$-spaces. In order to produce the latter, we use the
well known Gaussian measure construction that associates to any
continuous unitary representation $\pi$ of a locally compact,
second countable group $\Gamma$, a measure-preserving $\Gamma$-action. 
As explained in~\cite{Zimmer84} (see~5.2.13 and p.~111), one obtains a
continuum of nonisomorphic $\Gamma$-actions once $\Gamma$ has a
continuum of nonequivalent irreducible unitary representations
$\pi$, and furthermore, for any closed subgroup $H<\Gamma$, the
following holds: If $\pi|_H$ is weakly mixing, then $H$ acts
ergodically on the measure space constructed in this manner.

On the other hand,  it is a well known fact that any discrete infinite
 group $\Gamma$ admits a continuum of irreducible
unitary representations $\pi$ that are weakly contained
in $L^2(\Gamma)$ (this follows from Corollaire~1 in
J.~Dixmier~\cite{Dixmier64}, a remark for which we thank Bachir
Bekka).  But then, for any nonamenable closed subgroup $H<G$, the
restriction $\pi|_H$ must be weakly mixing, since otherwise we
would have (using $\prec$ to denote weak containment):
$$
\bone_H\subseteq (\pi\otimes \overline\pi)|_H\prec (L^2(G)\otimes L^2(G))|_H \cong \bigoplus_{n=1}^\infty
L^2(H),$$
contradicting  nonamenability of $H$ in view of the (generalised)
Hulanicki criterion. Applying this discussion to
$H=\Gamma_i<\Gamma$, one constructs a continuum of
irreducible nonisomorphic probability $\Gamma$-spaces, thereby
finishing the proof.
\phantom{over}\hfill\qed

\Subsec{Outer automorphisms of certain type $\mathrm{II}_1$ relations}
 The goal of this subsection is to present the

\vskip5pt
{\it Proof of Theorem~{\rm \ref{thm_out}}}.
We shall use the following notation of\break A.~Furman~\cite{Furman3}: If $(X,\mu)$ is any
probability
 $\Gamma$-space for a \pagebreak countable group~$\Gamma$,    let
 $$\aut^*(X,\Gamma)=\big\{F\in\aut(\R_{\Gamma,X})\ :\ \exists\,f\in\aut(\Gamma)\
\forall\,\gamma\in\Gamma:\ \ F(\gamma x)=f(\gamma)F(x) \big\}$$
and write $A^*(X,\Gamma)$ for the image of $\aut^*(X,\Gamma)$ in $\out(\R_{\Gamma,X})$. With this notation one deduces immediately the following from Theorem~\ref{thm_strong_rigidity1}:

\smallskip

\itshape Let $\Gamma = \Gamma_1 \times \Gamma_2$ be a torsion-free group with both $\Gamma_i$ in~$\creg${\rm ,} and let $(X,\mu)$ be an irreducible probability $\Gamma$-space. Then $\out(\R_{\Gamma,X}) = A^*(X,\Gamma)$\upshape.

\smallskip

Now let $K$ be a (second countable) compact group, and $\mu$ be its
normalised Haar measure. Let us fix
$K=\mathrm{SO}(n)$ with $n$ odd, which enjoys the property of having
both trivial center and no nontrivial outer automorphisms. Let $\Lambda$
be a Kazhdan group which admits a dense embedding into $K$ and
such that every injective homomorphism $\Lambda\to \Lambda$ is an
inner automorphism. We note that for every $n\geq 5$ one can
indeed find such a group $\Lambda$ which is a lattice in an appropriate higher
rank simple Lie group; indeed, the dense embedding into $K$ is provided
by a standard Galois twist argument, while for the condition on injective
homomorphisms $\Lambda\to \Lambda$ we refer to~\cite{Prasad76}.
Let $F_p$ and $F_q$ be non-Abelian free
groups with $p\neq q$ and consider the free products
$\Gamma_1=\Lambda*F_p$ and $\Gamma_2=\Lambda*F_q$. Suppose for the
time being that we are given injective homomorphisms of $F_p,F_q$
into $K$ such that the induced maps $\Gamma_i\to K$ are still
injective; then we can view $K$ as  a probability
$\Gamma=\Gamma_1\times\Gamma_2$-space by letting $\Gamma_1$ and
$\Gamma_2$ act by right and left multiplication respectively (it
is easily verified that essential freeness here is satisfied once
every open subgroup of $K$ is center free).

\begin{prop}
\label{prop_out_trivial}%
The group $\out(\R_{\Gamma,K})$ is trivial.
\end{prop}

\Proof
By the above reformulation of  Theorem~\ref{thm_strong_rigidity1},
it is enough to show that $A^*(K,\Gamma)$ is trivial. Since we
chose $F_p$ and $F_q$ nonisomorphic, every element of
$\aut^*(K,\Gamma)$ induces a (perhaps twisted) isomorphism of both
$\Gamma_1$\ti\ and\break $\Gamma_2$-actions individually. We shall
see that $\aut^*(K,\Gamma_1)\cap \aut^*(K,\Gamma_2)$ is trivial
(even though each of these two groups is large).

A direct argument of A.~Furman~\cite[7.2]{Furman3} enables one to describe\break
$\aut^*(K,\Gamma_1)$ as
$$\aut^*(K,\Gamma_1) = \big\{ a_{\sigma,t}(k) = t\sigma(k)\ :\ \sigma\in\aut(K),\, t\in K,\, \sigma(\Gamma_1)=t^{-1}\Gamma_1 t \big\}.$$
Now, since $\out(K)$ is trivial, we can write  $\sigma(k) = c^{-1}
k c$ for some $c\in K$; hence $a_{\sigma,t}(k) = t c^{-1} k c$
with the condition $c^{-1}\Gamma_1 c = t^{-1} \Gamma_1 t$; i.e.\ $(c
t^{-1})^{-1} \Gamma_1 (c t^{-1})\break =\Gamma_1$. Recall now that
$\Gamma_1 = \Lambda * F_p$. We claim that up to a conjugation in
$\Gamma_1$ every $f\in\aut(\Gamma_1)$ is trivial on $\Lambda$.
Indeed, since any action of the Kazhdan group $\Lambda$ on the
Bass-Serre tree associated with the free product $\Lambda
* F_p$ has a fixed vertex, it follows that $f(\Lambda)$ is
contained in a conjugate of $\Lambda$ or of $F_p$, the latter
being of course impossible. Hence after conjugation every
$f\in\aut(\Gamma_1)$ satisfies $f(\Lambda)\subseteq \Lambda$, and
by the choice of $\Lambda$ we deduce that after further
conjugation $f$ is trivial on $\Lambda$, proving the claim.

If we apply this and the claim above to the automorphism $f$ given
by conjugation by $c t^{-1}$, recalling that by density of
$\Lambda$ in $K$ every continuous automorphism of the latter which
is trivial on the former must be trivial, we find $t
c^{-1}\in\Gamma_1$ and hence conclude
$$\aut^*(K,\Gamma_1) = \big\{k\mapsto \gamma_1 k c_1\ :\ \gamma_1\in\Gamma_1,\,c_1\in K \big\}.$$
The analogous argument for $\Gamma_2$ yields
$$\aut^*(K,\Gamma_2) = \big\{k\mapsto c_2 k \gamma_2 \ :\ \gamma_2\in\Gamma_2,\,c_2\in K \big\}$$
since $\Gamma_2$ acts from the right. Thus, for
$$F\in \aut^*(K,\Gamma) \subseteq \aut^*(K,\Gamma_1)\cap \aut^*(K,\Gamma_2)$$
we have $F(k) = c_2 k \gamma_2 = \gamma_1 k c_1$ and therefore
$\gamma_1^{-1} c_2 k = k c_1 \gamma_2^{-1}$. Taking\break $k=e$ (or
rather $k$ sufficiently close to $e$ since these equalities hold
only almost everywhere), we deduce $\gamma_1^{-1} c_2 = c_1
\gamma_2^{-1}$. Since $K$ has trivial center this forces
$\gamma_1^{-1} c_2 = c_1 \gamma_2^{-1} = e$, i.e.\ $c_2 = \gamma_1$
and $c_1 = \gamma_2$ so that $\aut^*(K,\Gamma)$ (and more
generally the above intersection) consists of maps $k\mapsto
\gamma_1 k \gamma_2$ which are of course trivial in
$A^*(K,\Gamma)$. This concludes the proof of
Proposition~\ref{prop_out_trivial}.
\Endproof\vskip4pt 

We return now to the proof of  Theorem~\ref{thm_out}. Fix once and
for all one dense embedding of $\Lambda$ into $K$ as above. We
considered for the statement of Proposition~\ref{prop_out_trivial}
injective homomorphisms of $F_p,F_q$ into $K$ such that the
induced maps $\Gamma_i\to K$ are still injective. However,
such embeddings not only exist, but are generic
with respect to the Haar measure~\cite{Epstein71}; in particular there is a continuum
of nonconjugate such homomorphisms. Now
Proposition~\ref{prop_out_trivial} shows that each member of the corresponding
family of $\Gamma$-actions determines a relation with
trivial $\out(\R_{\Gamma,K})$. By a similar argument,
Theorem~\ref{thm_strong_rigidity1} implies that no two distinct
actions in this family can be WOE, since they are nonconjugate;
thus the relations are not weakly isomorphic, as required.
\hfill\qed

\Subsec{Some examples with linear groups} 
Using Howe-Moore's theorem, one can easily deduce from our OE
rigidity results applications to rigidity for linear groups acting
on homogeneous spaces. We bring here two examples.

\begin{exo}
Let $\overline{\Gamma}$ be a torsion-free
group in~$\cc$ with an injective homomorphism
$\ro:\overline{\Gamma}\to \mathrm{SL}_n({\bf Z})$. Form the
semi-direct product $\Gamma = {\bf Z}^n\rtimes_\ro
\overline{\Gamma}$, where $\overline{\Gamma}$ acts linearly on
${\bf Z}^n$ via $\ro$. This realises $\Gamma$ as a subgroup of
$\mathrm{SL}_{n+1}({\bf Z})<G =\mathrm{SL}_{n+1}({\bf R})$. Let
now $\Delta,\Sigma$ be any two lattices in $G$.

\smallskip

\itshape If the translation $\Gamma$-action on $G/\Delta$ is
{\rm WOE} to a $\Lambda$-translation action on $G/\Sigma${\rm ,} where
$\Lambda<G$ is any discrete subgroup{\rm ,} then there is an infinite
normal amenable subgroup $N\lhd \Lambda$ such that $\Lambda/N$ is
isomorphic to $\Gamma$.  \/{\rm (}\/In particular{\rm ,} the Zariski closure of
$\Lambda$ is not semi-simple.\/{\rm )}\/\upshape 

\smallskip

This statement is a straightforward application of
Theorem~\ref{thm_kernel_mystery_new} together with Howe-Moore's
theorem.
\end{exo}

In our last example we consider linear embeddings that are not necessarily discrete. In fact, the following is of interest precisely in the nondiscrete cases:

\begin{exo}
Let $F=F_p\times F_q$ be a product of non-Abelian free groups (or
of any torsion-free groups in~$\cc$). Let
$i_1,i_2:F\to G=\mathrm{SL}_n({\bf R})$ be two embeddings such
that the image of each free group under both embeddings is
unbounded. For simplicity, assume $n$ is odd so that $G$ has
trivial center. Let $\Delta,\Sigma$ be two lattices in $G$.

\medskip

\itshape If the $F$-translation action on $G/\Delta$ through
$i_1$ is {\rm WOE} to the $F$-translation action on $G/\Sigma$
through $i_2${\rm ,} then there is an automorphism $f$ of $F$ such that
the embeddings $i_1$ and $i_2$ are topologically equivalent modulo
$f$.

\medskip

More precisely{\rm ,} $i_2\circ f\circ i_1^{-1}$ extends to an
isomorphism between the closures of $i_1(F)$ and $i_2(F)$ in $G$.
In particular{\rm ,} $\overline{i_1(F)}$ is isomorphic to
$\overline{i_2(F)}$\upshape.
\end{exo}

{\it Proof}.
By Howe-Moore's theorem, our assumption on the embeddings ensures
that both $F$-actions are irreducible. By
Theorem~\ref{thm_factors_irred}, it follows that the actions are
isomorphic with respect to an automorphism $f$; we may assume
$f=\id$ upon composing one of the embeddings with $f$. If $(g_k)$
is a sequence of $F$ such that $i_1(g_k)$ tends to $e\in G$, then
by Howe-Moore $i_2(g_k)$ is bounded since our two actions are
isomorphic. Thus $i_2(g_k)$ has a limit point $g\in G$; now, since
$G$ has trivial center, $g=e$ because otherwise $g$ would not act
trivially on $G/\Sigma$. By symmetry of that argument we deduce
that $i_1(g_k)\to e$ if and only if $i_2(g_k)\to e$, as required.
\hfill\qed

\section{Background in bounded cohomology}
\label{sec_background_hb}
 \vglue-5pt

The purpose of this section is to offer the most elementary possible account of the bounded cohomology tools 
that we shall need. In the setting of this paper, it is possible to derive most relevant statements from two fundamental principles (Theorems~\ref{thm_amenable_resolution} and~\ref{thm_boundary_exists} below). Thus, we shall indicate some proofs for the reader's convenience. For a more detailed introduction, we refer to~\cite{Burger-Monod3}, \cite{Monod}.

We consider throughout the paper bounded cohomology for countable
(discrete) groups $\Gamma$. The coefficients will be taken almost
always in unitary $\Gamma$-representations on separable Hilbert
spaces. However, for the purpose of induction of such modules (see
Section~\ref{sec_induction}), it will be essential to allow the
following more general setting:
\vglue-20pt 
\phantom{up}
\begin{defi}
\label{defi_coefficients}%
A \emph{coefficient} $\Gamma$-module $(\pi,E)$ is an isometric linear\break $\Gamma$-representation
$\pi$
 on a Banach space $E$ such that: (i)~$E$ is the dual of some separable Banach space, (ii)~$\pi$ consists of
adjoint operators (in this duality).
\end{defi}
\phantom{up}
\vglue-20pt

The \emph{bounded cohomology} of $\Gamma$ with coefficient module $(\pi,E)$ is defined to be the cohomology of the complex
\begin{equation}
\label{eq_bar_complex}%
0\lra \ell^\infty(\Gamma,E)^\Gamma \lra \ell^\infty(\Gamma^2,E)^\Gamma \lra \ell^\infty(\Gamma^3,E)^\Gamma \lra \cdots
\end{equation}
of bounded invariant functions and is denoted by $\hb^\bu(\Gamma,E)$ or $\hb^\bu(\Gamma,\pi)$. This complex is just the subcomplex of bounded functions in the standard (homogeneous) bar complex for the Eilenberg-MacLane cohomology; in other word, invariance is understood with respect to the regular representation
$$\big(\lambda_\pi(\gamma)f\big)(\gamma_0, \ldots,\gamma_n) = \pi(\gamma)\Big(f(\gamma^{-1}\gamma_0, \ldots,\gamma^{-1}\gamma_n)\Big)$$
and the maps in~(\ref{eq_bar_complex}) are the usual (Alexander-Spanier) coboundary maps.

The usual cohomological methods do not apply to bounded cohomology, which has proved difficult to compute. It is therefore essential to have at least some replacement for the intractable complex~(\ref{eq_bar_complex}). More useful complexes arise in connection with standard Borel $\Gamma$-spaces with a finite quasi-invariant measure
and are such that the $\Gamma$-action is amenable in R.~Zimmer's~\cite{Zimmer84} sense. For short, we
call such a space an \emph{amenable $\Gamma$-space}.

\begin{thm}[\cite{Burger-Monod3}, \cite{Monod}]
\label{thm_amenable_resolution}%
Let $E$ be a coefficient $\Gamma$-module and $S$ an amenable $\Gamma$-space. Then the complex
\begin{equation}
\label{eq_resolution_amenable}%
0\lra \linfty(S,E)^\Gamma \lra \linfty(S^2,E)^\Gamma \lra \linfty(S^3,E)^\Gamma \lra \cdots
\end{equation}
realises canonically $\hb^\bu(\Gamma,E)$. The corresponding statement holds for the subcomplex of alternating cochains.\hfill\qed
\end{thm}

We do not make the meaning of \emph{canonically} more precise here, but its importance will be obvious in certain arguments below. In the above, $\linfty$ denotes the space of essentially bounded \weak measurable functions. Below, we will often deal with cases where $E$ is separable, in which case \weak and strong measurability coincide; hence the simpler notation $L^\infty$.

The point of Theorem~\ref{thm_amenable_resolution} is that there are indeed examples of amenable spaces with very strong ergodicity properties: The following result was established in~\cite{Burger-Monod3}, \cite{Monod} for finitely (or compactly) generated groups; the general version was then provided by V.~Kaimanovich~\cite{Kaimanovich03}.

\begin{thm}
\label{thm_boundary_exists}%
For every countable group $\Gamma${\rm ,}
  there is an amenable $\Gamma$-space $S$ such that for every separable coefficient $\Gamma$-module
$E${\rm ,} the space $L^\infty(S^2,E)^\Gamma$ is reduced to constant functions. In particular{\rm ,}
 there is a canonical
isomorphism
$$\hb^2(\Gamma,E)\cong Z\la (S,E)^\Gamma,$$
where $Z\la$ denotes the space of alternating cocycles.\hfill\qed
\end{thm}

(We point out that the conditions on $E$ are not merely technical, and that there are counter-examples if one drops either the separability assumption or the duality of $E$.)

A first immediate application of this fact is the following.

\begin{cor}
\label{cor_integrate}%
Let $\Gamma$ be a countable group and 
$(\pi_n,\H_n)_{n=1}^\infty$ a family of unitary $\Gamma$-representations in separable Hilbert spaces
$\H_n$. Then
$$\hb^2\left(\Gamma,{\textstyle\bigoplus_{n=1}^\infty}\H_n\right)=0\ \ \Longleftrightarrow\ \ \forall\ n\geq 1:\ \hb^2(\Gamma,\H_n)=0.$$
An analogous statement holds for direct integrals of unitary representations.\hfill\qed
\end{cor}

\vglue-6pt Here is another immediate consequence taken from~\cite{Burger-Monod3}, \cite{Monod}:

\vglue-20pt
\phantom{up}
\begin{cor}
\label{cor_Banach_space}%
Let $\Gamma$ be a countable group and $\alpha: E\to F$ an adjoint $\Gamma$-map of coefficient $\Gamma$-modules. If $F$ is separable{\rm ,} then the induced map $\alpha_*:\hb^2(\Gamma, E)\to \hb^2(\Gamma, F)$ is injective.
\end{cor}
 
{\it Proof}.
By Theorem~\ref{thm_amenable_resolution} and functoriality of~(\ref{eq_resolution_amenable}), the map $\alpha_*$ is realised at the level of $\hb^n(\Gamma,-)$ by the corresponding map
$$\alpha_n:\ \linftya(S^{n+1},E)\lra \la(S^{n+1},F)$$
for any amenable $\Gamma$-space $S$; observe that this map ranges in measurable functions because $\alpha$ is adjoint. On the other hand, $\alpha_n$ is injective at the level of \emph{cocycles}, so that the cohomological statement with $n=2$ follows from $\la(S^2,F)=0$ with $S$ as in Theorem~\ref{thm_boundary_exists}.
\Endproof 

We can also derive readily the following special case of a general exact sequence~\cite[\No12.0.2]{Monod}:

\begin{cor}
\label{cor_inf}%
Let $\Gamma$ be a countable group{\rm ,}
 $N\lhd \Gamma$ a normal subgroup and $Q=\Gamma/N$ the quotient. If $E$ is a separable coefficient
$Q$-module{\rm ,} then the inflation map
$$\mathrm{inf}:\ \hb^2(Q,E) \lra \hb^2(\Gamma,E)$$
is injective.
\end{cor}

{\it Proof}.
Let $S$ be an amenable $\Gamma$-space as in Theorem~\ref{thm_boundary_exists} and let $S'$ be the Mackey realisation of $L^\infty(S)^N$. Then $S'$ is an amenable $Q$-space satisfying the condition of Theorem~\ref{thm_boundary_exists}. Thus we have canonical embeddings
$$\la((S')^{n+1}, E)^Q \lra \la(S^{n+1}, E)^\Gamma$$
which induce the inflation $\hb^n(Q,E) \to \hb^n(\Gamma,E)$. Since $\la((S')^2, E)^Q$ vanishes, we deduce that at the level of $\hb^2$ the inflation is
 still injective.
\Endproof 

A key ingredient for our use of bounded cohomology in this paper is the following product formula whose proof relies also on
Theorem~\ref{thm_boundary_exists}.

\begin{thm}[\cite{Burger-Monod3}, \cite{Monod}]
\label{thm_Kunneth}%
Let $\Gamma=\Gamma_1\times \cdots\times \Gamma_n$ be a product of countable groups $\Gamma_i$ and let $(\pi,E)$ be a separable coefficient $\Gamma$-module. Then
$$\hb^2(\Gamma,E)\cong \bigoplus_{i=0}^n\hb^2\big(\Gamma_i,E^{\Gamma'_i}\big),$$
where $E^{\Gamma'_i}$ denotes the subspace of vectors fixed by $\Gamma'_i = \prod_{j\neq i}\Gamma_j$.\hfill\qed
\end{thm}

We will mostly use the consequence that $\hb^2(\Gamma,E)\neq 0$ implies $E^{\Gamma'_i}\neq 0$ for some $i$; we emphasize that the above formula
does not follow formally like some K\"unneth formula, and indeed may fail when $E$ is not separable (or
not dual). Since~\cite{Burger-Monod3}, \cite{Monod} deal with the general case of group extensions, we
indicated the simpler proof of the product case in~\cite{Monod-ShalomCRAS}.

We end with a simple fact that is well known in this setting (\cite{Johnson}, \cite{Gromov}, \cite{Ivanov}, \cite{Noskov}) but can also be seen as an application of Theorem~\ref{thm_amenable_resolution}:

\begin{prop}
\label{prop_kernel}%
Let $N\lhd \Gamma$ be an amenable normal subgroup of the countable group $\Gamma$ and let $E$ be a coefficient $\Gamma$-module. Then the inflation
$$\mathrm{inf}:\ \hb^n(\Gamma/N,E^N)\lra \hb^n(\Gamma,E)$$
is an isomorphism. In particular{\rm ,} $E^N\neq 0$ if $\hb^n(\Gamma,E)\neq 0$.
\end{prop}

\Proof
Let $S$ be $\Gamma/N$ endowed with some probability measure of full support. Then it is both an amenable $\Gamma$\ti\ and $\Gamma/N$-space. The statement now follows from Theorem~\ref{thm_amenable_resolution} by realising the inflation by the isomorphisms
\vskip12pt
\hfill $\displaystyle{\linfty(S^{n+1},E^N)^{(\Gamma/N)} = \linfty(S^{n+1},E^N)^\Gamma\cong
\linfty(S^{n+1},E)^\Gamma.}$ 
\hfill\qed

\begin{rem}
\label{rem_amenable_group}%
A degenerate case of Proposition~\ref{prop_kernel} occurs when $\Gamma$ itself is amenable: one deduces
then that $\hb^n(\Gamma,E)$ vanishes for all $n\geq 1$ and every coefficient $\Gamma$-module $E$.
\end{rem}

More advanced tools from~\cite{Burger-Monod3}, \cite{Monod} include a low degree exact sequence for group extensions (used below for Proposition~\ref{prop_subgroup_in_C}).

\section{Cohomological induction through couplings}
\label{sec_induction}%

Before considering cohomological induction, we start with some properties of the operation of inducing \emph{representations}. This is well known in the OE setting but becomes more transparent for ME.

\Subsec{Induced representations}
\label{sec_ind_rep}
Let $(\Sigma,m)$ be an ME coupling of two countable groups $\Lambda,\Gamma$ and let $(\pi,E)$ be
a unitary $\Gamma$-representation in a separable Hilbert space $E$, or more generally a separable
coefficient module (Definition~\ref{defi_coefficients}).

There are two equivalent ways to define the \emph{induced representation} $\Mind\Sigma\Gamma\Lambda\pi$. First, one can define the Banach space
\begin{equation}
\label{eq_defi_induced_space}%
\Mind\Sigma\Gamma\Lambda\pi = L^{[2]}(\Sigma, E)^\Gamma
\end{equation}
of $\Gamma$-equivariant measurable maps $f:\Sigma\to E$, wherein the notation $L^{[2]}$ means that the $\Gamma$-invariant function $\|f\|_E$ is to be in $L^2(\Gamma\bsl\Sigma)$. The $\Lambda$-action in this model is simply given by translation on $\Sigma$. Equivalently, one can define
$$\Mind\Sigma\Gamma\Lambda\pi = L^2(\Gamma\bsl\Sigma, E)$$
and endow it with the twisted $\Lambda$-action defined a.e.\ by
$$(\lambda f)(\Gamma x) = \pi(\chi(x)^{-1}\chi(\lambda^{-1}x))f(\lambda^{-1} \Gamma x),$$
where $\chi$ is a retraction as in Definition~\ref{defi_ME_cocycles}. Although this viewpoint is useful too, one should remember that the isomorphism between the latter and the more natural former depends on the choice of $\chi$. We also mention that upon identifying $\Gamma\bsl\Sigma$ with a fundamental domain $Y$, the action on $f\in L^2(Y,E)$ becomes the
well-known twisted action
$$(\lambda f)(y) = \pi(\beta(\lambda^{-1}, y)^{-1}) f(\lambda^{-1}\cdot y)$$
for the associated cocycle $\beta:\Lambda\times Y\to \Gamma$. This model is relevant when thinking of OE or WOE.

At any rate, $\Mind\Sigma\Gamma\Lambda\pi$ is a separable coefficient $\Lambda$-module, and a unitary representation if $E$ was unitary. The definition~(\ref{eq_defi_induced_space}) implies that the construction is natural and that one has the following \emph{transitivity} property: If $\Sigma'$ is an ME coupling of $\Lambda$ with a further group $\Delta$, then
\begin{equation}
\label{eq_Mind_transitive}%
\Mind{\Sigma'}\Lambda\Delta{(\Mind\Sigma\Gamma\Lambda\pi)} \cong \Mind{(\Sigma'\times_\Lambda \Sigma)}\Gamma\Delta\pi.
\end{equation}

Here are some elementary properties of the induction operation:

\begin{lemma}
\label{lemma_induce_properties}%
Let $\Lambda,\Gamma$ be two countable groups with commuting{\rm ,} measure-preserving actions on a measure
space $(\Sigma,m)$. Suppose that the $\Gamma$-action admits a finite measure fundamental domain and
that the $\Lambda$-action admits some fundamental domain. Then{\rm ,} for every unitary
$\Gamma$-representation  $\pi,\sigma$ one has\/{\rm :}\/
\begin{mitemize}
\item[{\rm (i)}] If $\pi$ is mixing for $\Gamma$ then 
$\Mind\Sigma\Gamma\Lambda\pi$ is mixing for $\Lambda$.

\item[{\rm (ii)}] If $\pi\prec\sigma$ then
 $\Mind\Sigma\Gamma\Lambda\pi\prec\Mind\Sigma\Gamma\Lambda\sigma$.

\item[{\rm (iii)}] $\Mind\Sigma\Gamma\Lambda\ell^2(\Gamma)\cong L^2(\Sigma)$ as $\Lambda$-representations.
\end{mitemize}
\end{lemma}

\Proof
This follows e.g.\ from Lemma~6.2 in~\cite{ShalomANNALS}, where for~(i) we use the fact that the
$\Lambda$-representation on $L^2(\Sigma)$ is mixing since the $\Lambda$-action on $\Sigma$ admits a
fundamental domain.
\hfill\qed

\begin{cor}
\label{cor_induce_amenablility}%
Let $M,N$ be two countable groups with commuting{\rm ,}\break measure-preserving actions on a
$\sigma$-finite measure space $(\Sigma,m)$. Suppose that $N$ has a fundamental domain in $\Sigma$ and
that $M$ is amenable and has a finite measure{\rm ,} fundamental domain. Then $N$ is amenable{\rm ,} too.
\end{cor}

\Proof
Recall that by A.~Hulanicki's criterion~\cite{Hulanicki66}, a group $\Lambda$ is amenable if and only if the regular representation $\ell^2(\Lambda)$ contains weakly the trivial representation $\bone_\Lambda$. Thus by assumption $\bone_M \prec \ell^2(M)$. The first and third points of Lemma~\ref{lemma_induce_properties} imply $\Mind\Sigma{M}{N}\bone_M \prec L^2(\Sigma)$. On the other hand, $\Mind\Sigma{M}{N}\bone_M=L^2(M\bsl \Sigma)$ contains $\bone_N$ since $M\bsl \Sigma$ has finite measure. Since $L^2(\Sigma)\cong \bigoplus_{n=1}^\infty\ell^2(N)$ we conclude that $\bone_N \prec \bigoplus_{n=1}^\infty\ell^2(N)$ and thus $N$ is amenable by a generalised version of A.~Hulanicki's criterion, see~\cite[7.3.6]{Zimmer84}.
\hfill\qed

\Subsec{Cohomological induction} \label{sec_ind_coh}
Given an ME coupling $(\Sigma,m)$ of two countable groups $\Lambda,\Gamma$ and separable coefficient $\Gamma$-module $(\pi,E)$, there is a natural way to induce the bounded cohomology. More specifically, we define a map
\begin{equation}
\label{eq_def_MMind}%
\MMind\Sigma\Gamma\Lambda:\ \hb^n(\Gamma,\pi) \lra \hb^n(\Lambda, \Mind\Sigma\Gamma\Lambda\pi)
\end{equation}
as follows. Let $\chi:\Sigma\to\Gamma$ be a retraction for the $\Gamma$-action on $\Sigma$. Realise the left hand side of~(\ref{eq_def_MMind}) by the complex~(\ref{eq_bar_complex}) whose $n^{\rm th}$ term is $\ell^\infty(\Gamma^{n+1}, E)^\Gamma$. Likewise, the $n^{\rm th}$ term for the right hand side is $\ell^\infty\big(\Lambda^{n+1}, L^{[2]}(\Sigma,E)^\Gamma\big)^\Lambda$. For every $f$ in the former, define $\MMind\Sigma\Gamma\Lambda f$ in the latter by
\begin{equation}
\label{eq_def_MMind_cocycle}%
\MMind\Sigma\Gamma\Lambda f(\lambda_0, \ldots, \lambda_n)(x) = f(\chi(\lambda_0^{-1}x), \ldots, \chi(\lambda_n^{-1}x)).
\end{equation}
It is straightforward to verify the statement:

\begin{lemma}
This map $\MMind\Sigma\Gamma\Lambda$ is a
 well defined linear map ranging in\break $\ell^\infty\big(\Lambda^{n+1},
L^{[2]}(\Sigma,E)^\Gamma\big)^\Lambda$. Moreover{\rm ,} it is continuous 
  and when $n$ varies one obtains a morphism of complexes.\hfill\qed
\end{lemma}

Observe that the very fact that $\MMind\Sigma\Gamma\Lambda$ yields cochains with square-summable coefficients would {\it a priori\/} not be true if we tried to induce general (unbounded) cocycles.

The main result of this section is:

\begin{thm}[Induction]
\label{thm_mind}%
Let $(\Sigma, m)$ be an {\rm ME} coupling of countable groups $\Lambda,\Gamma$. Then the induction map
$$\MMind\Sigma\Gamma\Lambda:\ \hb^2(\Gamma,\pi) \lra \hb^2(\Lambda, \Mind\Sigma\Gamma\Lambda\pi)$$
is injective for every separable coefficient $\Gamma$-module $(\pi,E)$. Moreover{\rm ,} it does not depend of the choice of $\chi$.
\end{thm}

The tools provided by Theorems~\ref{thm_amenable_resolution} and~\ref{thm_boundary_exists}~-- more specifically, Corollary~\ref{cor_Banach_space}~-- enable us to deduce Theorem~\ref{thm_mind} from the following:

\begin{prop}
\label{prop_mind_reduced}%
Let $(\Sigma, m)$ be an {\rm ME} coupling of countable groups $\Lambda,\Gamma${\rm ,}
 let $\chi:\Sigma\to\Gamma$ be a retraction for the $\Gamma$-action on $\Sigma$ and let $(\pi,E)$ be a
coefficient $\Gamma$-module. Then the map
\begin{eqnarray*}
\chi^*:\ \ell^\infty\big(\Gamma^{n+1}, \linfty(\Lambda\bsl\Sigma,E)\big)^\Gamma &\lra &\ell^\infty\big(\Lambda^{n+1}, \linfty(\Sigma,E)^\Gamma\big)^\Lambda\\
\chi^*f(\lambda_0, \ldots, \lambda_n)(s) &= &f(\chi(\lambda_0^{-1}s), \ldots, \chi(\lambda_n^{-1}s))(s)
\end{eqnarray*}
induces an injection
$$\mathfrak s_\Sigma:\ \hb^n\big(\Gamma, \linfty(\Lambda\bsl \Sigma,E)\big) \lra \hb^n\big(\Lambda, \linfty(\Sigma,E)^\Gamma\big).$$
Moreover{\rm ,} $\mathfrak s_\Sigma$ does not depend on the choice of $\chi$.
\end{prop}

\phantom{up}
\vglue-18pt

{\it Proof that Proposition~{\rm \ref{prop_mind_reduced}} implies Theorem~{\rm \ref{thm_mind}}}.
Let $E$ be separable. In view of   formula~(\ref{eq_def_MMind_cocycle}), the induction map
$\MMind\Sigma\Gamma\Lambda$ factors as
\begin{multline}
\label{eq_mind_factors}%
\hb^n(\Gamma,E) \xrightarrow{\ \varepsilon_*\ }\hb^n\big(\Gamma, L^\infty(\Lambda\bsl\Sigma,E)\big) \xrightarrow{\ \mathfrak s_\Sigma\ }\\
\hb^n\big(\Lambda, L^\infty(\Sigma,E)^\Gamma\big) \xrightarrow{\ \iota_*\ } \hb^n\big(\Lambda, L^{[2]}(\Sigma,E)^\Gamma\big),
\end{multline}
where $\varepsilon: E \to L^\infty(\Sigma,E)^\Lambda$ is the inclusion of constant functions and $\iota$ is the inclusion of $L^\infty(\Sigma,E)^\Gamma$ into $L^{[2]}(\Sigma,E)^\Gamma= \Mind\Sigma\Gamma\Lambda\pi$. The map $\varepsilon$ admits a $\Gamma$-equivariant right inverse given by integration over the finite measure space $\Lambda\bsl\Sigma$. Therefore, the first map in~(\ref{eq_mind_factors}) is injective. The second map is injective and independent of $\chi$ by Proposition~\ref{prop_mind_reduced}. Since $\iota$ is the dual of the inclusion of $L^{[2]}(\Sigma,E)^\Gamma$ into $L^{[1]}(\Sigma,E)^\Gamma$ and $L^{[2]}(\Sigma,E)^\Gamma$ is separable, we are in situation to apply Corollary~\ref{cor_Banach_space} and $\iota_*$ is also injective for $n=2$. This is the only point where we use the separability of $E$ and $n=2$.
\Endproof  

Thus we are left to prove Proposition~\ref{prop_mind_reduced}. Given the importance of Theorem~\ref{thm_mind} for the paper, we shall present two proofs: first, we give a functorial proof that actually shows more. Then, for the convenience of a reader who would want to avoid the use of cohomological machinery, we outline an independent down-to-earth proof.

\vglue-4pt
\Subsec{Functorial proof} We prove a more general ``reciprocity'' statement:

\vglue-18pt
\phantom{up}

\begin{prop}
\label{prop_reciprocity}%
Let $\Lambda,\Gamma$ be countable groups{\rm ,}
 $S$ a standard measure space with measure class preserving $\Lambda\times\Gamma$-action and let
$(\pi,E)$ be a coefficient $\Lambda\times\Gamma$-module.
If both the $\Lambda$\ti\ and $\Gamma$-actions on $S$ are amenable{\rm ,}
 then there is a canonical isomorphism
\begin{equation}
\label{eq_mfraks}%
\mathfrak s:\ \hb^n\big(\Gamma, \linfty(S,E)^\Lambda\big) \cong \hb^n\big(\Lambda, \linfty(S,E)^\Gamma\big)
\end{equation}
for all $n\geq0$. Moreover{\rm ,}
 if the $\Gamma$-action on $S$ admits a retraction $\chi:S\to\Gamma${\rm ,} then the map
\begin{eqnarray*}
\chi^*:\ \ell^\infty\big(\Gamma^{n+1}, \linfty(S,E)^\Lambda\big)^\Gamma &\lra &\ell^\infty\big(\Lambda^{n+1}, \linfty(S,E)^\Gamma\big)^\Lambda\\
\chi^*f(\lambda_0, \ldots, \lambda_n)(s) &= &f(\chi(\lambda_0^{-1}s), \ldots, \chi(\lambda_n^{-1}s))(s)
\end{eqnarray*}
induces $\mathfrak s$ on cohomology.
\end{prop}

This implies indeed Proposition~\ref{prop_mind_reduced} since an action with fundamental domain is amenable (and the $\Lambda$-representation on $E$ is taken to be trivial).

\demo{Proof of Proposition~{\rm \ref{prop_reciprocity}}}
The $\Lambda\times\Gamma$-action on $\Gamma^{n+1}\times S$ given by diagonal $\Gamma$-action on $\Gamma^{n+1}\times S$ and $\Lambda$-action on $S$ is amenable because it is isomorphic to the $\Lambda\times\Gamma$-action given by diagonal $\Gamma$-action on $\Gamma^{n+1}$ and $\Lambda$-action on $S$. Therefore the coefficient $\Lambda\times\Gamma$-module
$$\ell^\infty\big(\Gamma^{n+1}, \linfty(S,E)\big) \cong \linfty(\Gamma^{n+1}\times S, E)$$
is relatively injective; see~\cite[\No5.7.1]{Monod}. Likewise, $\ell^\infty(\Lambda^{n+1},
\linfty(S,E)\big)$ is relatively injective for $\Lambda\times\Gamma$. Therefore, there is a
$\Lambda\times\Gamma$-morphism of complexes 
$$\ell^\infty\big(\Gamma^{\bu+1}, \linfty(S,E)\big) \to \ell^\infty\big(\Lambda^{\bu+1}, \linfty(S,E)\big)$$
and any two such morphisms are $\Lambda\times\Gamma$-homotopic. (This follows
from~\cite[\S 7]{Monod}; the complexes above have indeed a contracting homotopy by evaluation of the
first variable.) In particular -- and due to the symmetry between $\Lambda$ and $\Gamma$ -- there is a
canonical isomorphism between the cohomology of the associated nonaugmented complexes of
$\Lambda\times \Gamma$-invariants. But those identify to the complexes with $n^{\rm th}$ term
$\ell^\infty(\Gamma^{n+1}, \linfty(S,E)^\Lambda)^\Gamma$, respectively $\ell^\infty(\Lambda^{n+1},
\linfty(S,E)^\Gamma)^\Lambda$, which compute canonically both sides of~(\ref{eq_mfraks}), whence the
first part of the proposition.

If now we have a retraction $\chi$, then   the formula for $\chi^*$ also yields an example of such a
morphism of complexes, thus inducing the same map since all morphisms are
$\Lambda\times\Gamma$-homotopic.
\Endproof \vskip4pt 

Observe that this proof shows that both sides in~(\ref{eq_mfraks}) are canonically isomorphic to $\hb^n\big(\Gamma\times\Lambda, \linfty(S,E)\big)$.
(Note that this situation contrasts sharply with Theorem~\ref{thm_Kunneth} and illustrates the necessity of the separability assumption in the latter, an assumption not satisfied by $\linfty(S,E)$ above.)

\Subsec{Another proof} We briefly  indicate here another way to deduce
Proposition~\ref{prop_mind_reduced}, starting with an elementary
proof for a special case (Lemma~\ref{lemma_induction_inclusion})
and then dealing with increasing levels of generality. First we
state the following fact, skipping the tedious verification of the
computation (recall that it follows anyway from the previous
functorial proof).

\begin{lemma}
\label{lemma_s_independent}%
The map $\mathfrak s_\Sigma$ in Proposition~{\rm \ref{prop_mind_reduced}}
 does not depend on the choice of $\chi$.\hfill\qed
\end{lemma}

\begin{lemma}
\label{lemma_induction_inclusion}%
Proposition~{\rm \ref{prop_mind_reduced}}
 holds if the $\Lambda$-action on $\Sigma$ admits a fundamental domain contained in some fundamental
domain for the $\Gamma$-action.
\end{lemma}

\Proof
In that case, we can choose retractions $\chi:\Sigma\to \Gamma$ and $\kappa:\Sigma\to\Lambda$ such that $\kappa^{-1}(e_\Lambda)\subseteq \chi^{-1}(e_\Gamma)$. If we define now $\kappa^*$ by a formula analogous to that of Proposition~\ref{prop_mind_reduced}, we have
\begin{eqnarray*}
\kappa^*\chi^* f(\gamma_0, \ldots, \gamma_n)(s) &=&
 \chi^* f(\kappa(\gamma_0^{-1}s),\ldots,\kappa(\gamma_n^{-1}s))(s)   \\
&=&f(\chi(\kappa(\gamma_0^{-1}s)^{-1}s),\ldots,\chi(\kappa(\gamma_n^{-1}s)^{-1}s))(s).
\end{eqnarray*}
Since
$$\kappa(\gamma_i^{-1}s)^{-1}s = \gamma_i \kappa(\gamma_i^{-1}s)^{-1}\gamma_i^{-1}s \in \gamma_i \kappa^{-1}(e_\Lambda)\subseteq \gamma_i \chi^{-1}(e_\Gamma),$$
we deduce $\kappa^*\chi^* f=f$ so that $\chi^*$ has a left inverse as morphism of complexes. This implies the injectivity of $\mathfrak s_\Sigma$ in view of Lemma~\ref{lemma_s_independent}.
\hfill\qed

\begin{cor}
\label{cor_induction_inequality}%
Proposition~{\rm \ref{prop_mind_reduced}} holds if the coupling is ergodic and\break
$[\Gamma:\Lambda]_\Sigma\leq 1$.
\end{cor}

{\it Proof}. 
In the ergodic situation, the full group of automorphisms of the $\Gamma\times \Lambda$-action on
 $\Sigma$ acts transitively on sets of equal measure (up to null-sets). Therefore, we can find a
$\Lambda$-fundamental domain contained in some fundamental domain for the $\Gamma$-action. Now we
can apply Lemmas~\ref{lemma_induction_inclusion} and~\ref{lemma_s_independent}.
\hfill\qed

\begin{prop}
\label{prop_induction_inequality}%
Proposition~{\rm \ref{prop_mind_reduced}} holds if
$[\Gamma:\Lambda]_\Sigma\leq 1$ \/{\rm (}\/without ergodicity assumption\/{\rm ).}\/
\end{prop}

{\it Proof}. 
By Lemma~2.2 in~\cite{Furman1}, any ME coupling $(\Sigma,m)$ can be disintegrated into ergodic couplings $(\Sigma_t,m_t)$ with $m=\int_T m_t\,d\eta(t)$ for some probability space $(T,\eta)$. The set
$$T_+ = \big\{t\in T\ :\ [\Gamma:\Lambda]_{\Sigma_t} \leq 1\big\}$$
has positive $\eta$-measure since $[\Gamma:\Lambda]_\Sigma\leq 1$. Let $\Sigma_+$ be the coupling obtained by integrating $m_t$ over $T_+$; using Corollary~\ref{cor_integrate}, one checks that the injectivity of $\MMind{\Sigma_+}\Gamma\Lambda$ follows from the injectivity of $\eta$-a.e.\ $\MMind{\Sigma_t}\Gamma\Lambda$, which is granted by Corollary~\ref{cor_induction_inequality}. One checks similarly that the injectivity of $\MMind\Sigma\Gamma\Lambda$ follows from the injectivity of $\MMind{\Sigma_+}\Gamma\Lambda$.
\Endproof \vskip4pt 

Using Lemma~\ref{lemma_s_independent}, one verifies:

\begin{lemma}
\label{lemma_ind_transitive}%
Let $\Sigma$ be an {\rm ME}
 coupling of countable groups $\Lambda,\Gamma$ and $\Sigma'$ an {\rm ME} coupling of countable groups
$\Delta,\Lambda$. Then{\rm ,} for every coefficient $\Gamma$-module $(\pi, E)${\rm ,} $\mathfrak
s_{\Sigma'}\mathfrak s_\Sigma =\mathfrak s_{\Sigma'\times_\Lambda \Sigma}$ and hence $\
\MMind{\Sigma'}\Lambda\Delta{(\MMind\Sigma\Gamma\Lambda\pi)} =
\MMind{(\Sigma'\times_\Lambda \Sigma)}\Gamma\Delta\pi$.\hfill\qed
\end{lemma}

{\it End of second proof of Proposition~{\rm \ref{prop_mind_reduced}}}.
Consider the composed ME coupling of $\Gamma$ with itself $\Omega = \check\Sigma\times_\Lambda \Sigma$. In view of Lemma~\ref{lemma_ind_transitive}, it is enough to show that the map
$$\mathfrak s_\Omega:\ \hb^n\big(\Gamma, \linfty(\Gamma\bsl \Omega,E)\big) \lra \hb^n\big(\Gamma, \linfty(\Omega,E)^\Gamma\big)$$
defined as in Proposition~\ref{prop_mind_reduced} is injective. By~(\ref{eq_constant_functorial}) 
$$[\Gamma:\Gamma]_\Omega = [\Gamma:\Lambda]_\Sigma\cdot[\Lambda:\Gamma]_{\check\Sigma}=1,$$
so that we may conclude by applying Proposition~\ref{prop_induction_inequality} to $\Omega$.
\hfill\qed

\section{Strong rigidity} 
\label{sec_strong}%
\vglue-12pt
\Subsec{Strong rigidity for products}  
Our first goal is to prove Theorem~\ref{thm_factors} from the introduction, for the more general class of
groups~$\cc$ defined in~\ref{defi_C}. The use of bounded cohomology in the proof of
Theorem~\ref{thm_factors}  is detailed in the following result
which we isolate for further reference:

\begin{prop}
\label{prop_fund_domain}%
Let $(\Sigma,m)$ be an {\rm ME} coupling of $\Gamma = \Gamma_1\times\cdots\times \Gamma_n$ 
with $\Lambda = \Lambda_1\times \cdots\times \Lambda_{n'}${\rm ,}
 where $\Gamma_i$ are torsion-free groups in~$\cc$ and $\Lambda_j$ are any countable groups.
Then there are a surjective map $t:\{1,\ldots,n\}\to \{1,\ldots,n'\}$ and a fundamental domain $Y\subseteq
\Sigma$ for the $\Gamma$-action such that $\Lambda'_{t(i)} Y\subseteq \Gamma'_i Y$ for all $i$.
Moreover{\rm ,}
 if   $n = n'${\rm ,} then $\Lambda_{t(i)} Y\subseteq \Gamma_i Y$ for all $i$ and the groups
$\Lambda_j$ are also in~$\cc$.
\end{prop}

{\it Proof}.
First we perform an induction on $1\leq k\leq n$ to construct a map $t$ and fundamental domains $Y_k\subseteq \Sigma$ for the $\Gamma$-action such that
\begin{equation}
\label{eq_inclusion_inuction}%
\Lambda'_{t(i)} Y_k\ \subseteq\ \Gamma'_i Y_k\kern1cm\forall\ 1\leq i\leq k.
\end{equation}
Let $\pi_1$ be a mixing unitary representation of $\Gamma_1$ such that $\hb^2(\Gamma_1,\pi_1)\neq 0$. We also denote by $\pi_1$ the corresponding representation of $\Gamma$ (factoring through $\Gamma_1$), and recall that by Corollary~\ref{cor_inf} we have $\hb^2(\Gamma,\pi_1)\neq 0$. Therefore, applying Theorem~\ref{thm_mind}, we deduce that $\hb^2(\Lambda, \Mind\Sigma\Gamma\Lambda{\pi_1})$ is 
nonx-trivial. We apply now the product formula of Theorem~\ref{thm_Kunneth} to find some $1\leq t(1)\leq n'$ such that
\begin{equation}
\label{eq_factors_coh}%
\hb^2\left(\Lambda_{t(1)}, (\Mind\Sigma\Gamma\Lambda{\pi_1})^{\Lambda'_{t(1)}}\right)\ \neq 0.
\end{equation}
Let $Y_0$ be any fundamental domain for $\Gamma$, and $\beta: \Lambda\times Y_0\to \Gamma$ be the associated cocycle. Take for $\Mind\Sigma\Gamma\Lambda{\pi_1}$ the model $L^2(Y_0,\H_{\pi_1})$ with $\beta$-twisted action. Since~(\ref{eq_factors_coh}) forces $(\Mind\Sigma\Gamma\Lambda{\pi_1})^{\Lambda'_{t(1)}}\neq 0$, this amounts to the existence of a nonzero measurable function $f: Y_0\to \H_{\pi_1}$ such that $f(\lambda'\cdot x) = \pi_1(\beta(\lambda', x))f(x)$ for every $\lambda'\in \Lambda'_{t(1)}$ and almost every $x\in Y_0$. Since the $\Gamma_1$-representation on $\H_{\pi_1}$ is mixing, the $\Gamma$-representation is tame and thus the cocycle reduction lemma (Lemma~5.2.11 in~\cite{Zimmer84}) can be applied to the restriction $\beta: \Lambda'_{t(1)}\times Y_0\to \Gamma$. Mind that the $\Lambda'_{t(1)}$-action on $Y_0$ is not assumed ergodic, but since $\Gamma_1$ is torsion-free and $\pi_1(\Gamma_1)$ mixing,
 the only possible stabiliser in $\Gamma$ of
 nonzero
  elements of $\H_{\pi_1}$ is $\Gamma'_1$. Thus cocycle reduction applied to every ergodic component yields a measurable $\fhi:Y_0\to \Gamma$ such that
$$\beta'(\lambda, x) = \fhi(\lambda\cdot x)\beta(\lambda,x)\fhi(x)^{-1}$$
ranges in $\Gamma'_1$ whenever $\lambda\in\Lambda'_{t(1)}$. Observe that if we replace $\fhi$ by its composition with the projection $\Gamma\to\Gamma_1<\Gamma$ and take for $\beta'$ the resulting cocycle $\beta': \Lambda\times Y_0\to \Gamma$, we have still $\beta'(\Lambda'_{t(1)}\times Y_0)\subseteq \Gamma'_1$ almost everywhere; we make this change. We consider the new fundamental domain
$$Y_1 = \big\{\fhi(x)x\ :\ x\in Y_0\big\}$$
for the $\Gamma$-action on $\Sigma$. For every $\lambda'\in\Lambda'_{t(1)}$ and almost every $x\in
Y_0$ 
\begin{eqnarray*}
\lambda'\fhi(x) x& =& \fhi(x)\lambda' x = \fhi(x)\beta(\lambda', x)^{-1}\lambda'\cdot x \\
&=& \big(\fhi(\lambda'\cdot x)\beta(\lambda', x)\fhi(x)^{-1}\big)^{-1}\fhi(\lambda'\cdot x)\lambda'\cdot x\
\in\ \Gamma'_1\fhi(\lambda'\cdot x)\lambda'\cdot x  \subseteq  \Gamma'_1 Y_1.
\end{eqnarray*}
This shows
\begin{equation}
\label{eq_fund_inclusion}%
\Lambda'_{t(1)} Y_1\ \subseteq\ \Gamma'_1 Y_1.
\end{equation}
Let now $k\geq 2$ and assume that  $Y_{k-1}$ and $t:\{1,\ldots,k-1\}\to \{1,\ldots,n'\}$ such
that~(\ref{eq_inclusion_inuction}) holds for $k-1$. We take now a mixing unitary representation $\pi_k$ of
$\Gamma_k$ such that $\hb^2(\Gamma_k,\pi_k)\neq 0$. Arguing as above, there is $1\leq t(k) \leq n'$ with
\begin{equation}
\label{eq_factors_coh2}%
\hb^2\left(\Lambda_{t(k)}, (\Mind\Sigma\Gamma\Lambda{\pi_k})^{\Lambda'_{t(k)}}\right)\ \neq 0 
\end{equation}
and thus $(\Mind\Sigma\Gamma\Lambda{\pi_k})^{\Lambda'_{t(k)}}\neq 0$. We perform again a reduction of cocycle and obtain as above a map $\psi: Y_{k-1}\to \Gamma_k$ such that
$$Y_k = \big\{\psi(x)x\ :\ x\in Y_{k-1}\big\}$$
is a fundamental domain for the $\Gamma$-action on $\Sigma$ satisfying
\begin{equation}
\label{eq_fund_inclusion2}%
\Lambda'_{t(k)} Y_k\ \subseteq\ \Gamma'_k Y_k.
\end{equation}
We claim that $Y_k$ still satisfies
\begin{equation}
\label{eq_fund_inclusion3}%
\Lambda'_{t(i)} Y_k\ \subseteq\ \Gamma'_i Y_k\kern1cm\forall\ 1\leq i\leq k-1.
\end{equation}
Indeed, for all $\lambda'\in\Lambda'_{t(i)}$ and almost every $x\in Y_{k-1}$ we have $\lambda'\psi(x)x = \psi(x)\lambda'x$ which by~(\ref{eq_fund_inclusion}) is $\psi(x) \gamma' y$ for some $\gamma'\in\Gamma'_i$ and $y\in Y_{k-1}$. Thus
$$\lambda'\psi(x)x = \big(\psi(x)\gamma'\psi(y)^{-1}\big)\psi(y)y\ \in\ \Gamma'_i Y_k,$$
as claimed.

Thus~(\ref{eq_fund_inclusion2}) and~(\ref{eq_fund_inclusion3})
complete the induction step to
prove~(\ref{eq_inclusion_inuction}). We let now $Y=Y_k$, so that
for the first claim of Proposition~5.1,  it remains only to show
that $t$ is surjective onto $\{1,\ldots,n'\}$. Suppose for a
contradiction that there is $1\leq j\leq n'$ not in the image of
$t$. Then $\Lambda_j$ is contained in $\bigcap_{i=1}^n
\Lambda'_{t(i)}$, so that in view of~(\ref{eq_inclusion_inuction})
we have
$$\Lambda_j Y \subseteq \bigcap_{i=1}^n \Gamma'_i Y = Y,$$
contradicting the properness of the $\Lambda$-action, since $\Lambda_j$ is 
nonx-trivial.

Under the additional assumption $n = n'$, the map $t$ must be bijective. Then, using~(\ref{eq_inclusion_inuction}), we have for all $1\leq i\leq n$ the inclusion
\begin{equation}
\label{eq_inclusion_individual}%
\Lambda_{t(i)} Y = \bigcap_{j\neq t(i)}\Lambda'_j Y \subseteq \bigcap_{\ell\neq i}\Gamma'_\ell Y = \Gamma_i Y,
\end{equation}
as claimed. It remains only to show that $\Lambda_{t(i)}$ is in~$\cc$ for all $i$. Let $\pi_i$ be the mixing $\Gamma_i$-representation that gave~(\ref{eq_factors_coh2}) for $k=i$ in the inductive argument above; in particular~(\ref{eq_factors_coh2}) implies
\begin{equation}
\label{eq_factors_final}%
\hb^2\left(\Lambda_{t(i)}, \Mind\Sigma\Gamma\Lambda{\pi_i}\right)\ \neq 0.
\end{equation}
The inclusion~(\ref{eq_inclusion_individual}) shows that the set $\Sigma_i = \Gamma_i Y\subseteq \Sigma$ is $\Lambda_{t(i)}\times \Gamma_i$-invariant. The $\Lambda_{t(i)}$-action on $\Sigma_i$ admits some fundamental domain because it is a subspace of $\Sigma$, and the $\Gamma_i$-action admits the finite measure fundamental domain $Y$. Therefore the $\Lambda_{t(i)}$-representation $\tau = \Mind{\Sigma_i}{\Gamma_i}{\Lambda_{t(i)}}{\pi_i}$ is mixing by point~(i) in Lemma~\ref{lemma_induce_properties}. On the other hand, since $\Sigma_i$ is a fundamental domain for the 
$\Gamma'_i$-action on $\Sigma$, and $\Gamma'_i$ acts trivially on $\H_{\pi_i}$, we have
$$\H_\tau = L^{[2]}(\Sigma_i, \H_{\pi_i})^{\Gamma_i} \cong L^{[2]}(\Sigma, \H_{\pi_i})^\Gamma.$$
Therefore, the $\Lambda_{t(i)}$-representation $\tau$ is isomorphic to the restriction to $\Lambda_{t(i)}$ of $\Mind\Sigma\Gamma\Lambda{\pi_i}$. Thus~(\ref{eq_factors_final}) shows indeed that $\Lambda_{t(i)}$ is in~$\cc$.
\hfill\qed

\demo{Proof of Theorem~{\rm \ref{thm_factors}} for the class~$\cc$ instead of~$\creg$}
Whenever two groups $\Gamma,\Lambda$ are measure
equivalent, one can assume by disintegration that there is an
\emph{ergodic} ME coupling $(\Sigma,m)$ of $\Gamma$ with
$\Lambda$; see~\cite[2.2]{Furman1}. First we apply
Proposition~\ref{prop_fund_domain} to obtain a bijection $t$ and a
$\Gamma$-fundamental domain $Y$ with $\Lambda_{t(i)} Y
\subseteq \Gamma_i Y$ for all $i$. Since by the proposition all
$\Lambda_j$ are in~$\cc$ and $n=n'$, the assumptions of the theorem
(with class~$\cc$) now also hold if we reverse
$\Gamma$ and $\Lambda$. Therefore, a second application of
Proposition~\ref{prop_fund_domain} gives us a bijection $s$ and a
fundamental domain $X$ for the $\Lambda$-action on $\Sigma$
such that $\Gamma_{s(j)} X \subseteq \Lambda_j X$ for all $j$.
Since all groups commute, the properties of our $\Gamma$-fundamental domain $Y$ are not altered if we translate it by an
element of $\Gamma$; therefore, we may and do assume that the
intersection $A=X\cap Y$ has positive measure. We claim that the
bijections $s$ and $t$ are inverse to each other.

Indeed, pick $1\leq j\leq n$, write $i=s(j)$ and let us show that $t(i)=j$. Write $A_\Gamma$ for the image of $A$ in $\Gamma\bsl\Sigma$ and apply Poincar\'e recurrence to the $\Lambda_{t(i)}$-action on $\Gamma\bsl\Sigma$. This yields a nontrivial element $\lambda\in\Lambda_{t(i)}$ such that $\lambda A_\Gamma\cap A_\Gamma$ has positive measure. Therefore, there is $\gamma\in\Gamma$ such that $B=\lambda A\cap \gamma A$ has positive measure. Since $A\subseteq Y$, the inclusion $\Lambda_{t(i)} Y \subseteq \Gamma_i Y$ implies $B\subseteq \gamma Y\cap \Gamma_i Y$ and hence $\gamma\in \Gamma_i$. But now $\gamma A\subseteq \Gamma_i X = \Gamma_{s(j)} X\subseteq \Lambda_j X$ so that $B\subseteq \Lambda_j X\cap\lambda X$. It follows that $\lambda\in\Lambda_j$, and since $\lambda$ is nontrivial we conclude $t(i)=j$, proving the claim.

In conclusion, we may permute the indices in such a way that
\begin{equation}
\label{eq_matching_domains}%
\Lambda_i Y \subseteq \Gamma_i Y \kern.5cm \text{and}\kern.5cm \Gamma_i X \subseteq \Lambda_i X \kern1cm\forall\ 1\leq i\leq n.
\end{equation}

Choose now some $i$ and let $\overline{\Sigma}$ be the space of ergodic components of the $\Gamma'_i\times \Lambda'_i$-action on $\Sigma$. We shall show that for an appropriate measure $\nu$ the space $(\overline{\Sigma},\nu)$ yields an ME coupling of $\Gamma_i$ with $\Lambda_i$ for the natural $\Gamma_i\times\Lambda_i$-action inherited from $\Sigma$. 
Note that we cannot take for $\nu$ the projection of $m$, since this could e.g.\ give infinite measure to every point in $\overline{\Sigma}$. Instead, we identify $\overline{\Sigma}$ with the space of ergodic components of the $\Lambda'_i$-action on $\Gamma'_i\bsl\Sigma$, the latter being measurably identified with $\Gamma_i Y$. We let $\nu$ be the projection of $m|_{\Gamma_i Y}$ under the corresponding map $\Gamma_i Y\to \overline{\Sigma}$.

We claim that $\Gamma_i$ has a $\nu$-finite fundamental domain in $\overline{\Sigma}$. Indeed, let $\overline{Y}$ be the image of $Y$ in $\overline{\Sigma}$; one has $\Gamma_i\overline{Y}=\overline{\Sigma}$ and $\nu(\overline{Y})$ is finite because in view of~(\ref{eq_matching_domains}) the pre-image of $\overline{Y}$ in $\Gamma'_i\bsl\Sigma$ is just $Y$ under the identification $\Gamma'_i\bsl\Sigma \cong \Gamma_i Y$. Therefore, to conclude that $\overline{Y}$ is a finite measure fundamental domain for $\Gamma_i$ in $(\overline{\Sigma},\nu)$, it remains to show that $\gamma_i\in\Gamma_i$ must be trivial whenever $\gamma_i \overline{Y} \cap \overline{Y}$ is nonnull for $\nu$. This is indeed the case, since then $\gamma_i\Gamma'_i\Lambda'_i Y\cap \Gamma'_i\Lambda'_i Y$ has positive measure in $\Sigma$; by~(\ref{eq_matching_domains}), this set is $\gamma_1\Gamma'_i Y\cap \Gamma'_i Y$ so that $\gamma_i\in\Gamma'_i$ and hence $\gamma_i$ is trivial.

Arguing in a symmetric manner, we consider the {\it a priori\/} different measure $\nu'$ on $\overline{\Sigma}$ obtained by projecting $m|_{\Lambda_i X}$ under $\Lambda_i X\to \overline{\Sigma}$ and deduce that $\Lambda_i$ has a $\nu'$-finite fundamental domain $\overline{X}$ in $\overline{\Sigma}$. In order to complete the proof that $(\overline{\Sigma},\nu)$ yields an ME coupling of $\Gamma_i$ with $\Lambda_i$, it is enough to show that $\nu'$ is a scalar multiple of $\nu$. The measures $\nu,\nu'$ are absolutely continuous with respect to each other since they are projected from $\Gamma_i Y$ and $\Lambda_i X$ respectively. Since we started by reducing to an ergodic ME coupling $(\Sigma,m)$, the measures $\nu,\nu'$ are both ergodic for the $\Gamma_i\times\Lambda_i$-action, so due to absolute continuity they are a scalar multiple one of the other. This completes the proof of Theorem~\ref{thm_factors}.
\phantom{over}\hfill\qed

\demo{Proof of Theorem~{\rm \ref{thm_factors_irred}}}
Proposition~\ref{prop_fund_domain} shows that the groups $\Lambda_j$ are also in~$\cc$. Therefore, under the assumption~(ii) the situation is symmetric and we may switch $\Gamma$ with $\Lambda$ if necessary so that in either case the assumption~(i) holds; in other words, we may assume
\begin{equation}
\label{eq_inequality_volumes}%
m(\Gamma\bsl\Sigma)\ \geq m(\Lambda\bsl\Sigma).
\end{equation}
Moreover,  Proposition~\ref{prop_fund_domain} gives us a
fundamental domain $Y\subseteq \Sigma$ for $\Gamma$ such that
after possibly permuting indices
\begin{equation}
\label{eq_fund_inclusion4}%
\Lambda_i Y\ \subseteq\ \Gamma_i Y \kern1cm (1\leq i\leq n).
\end{equation}
Fix some $i$. By the above, the subset $C=\Gamma'_i Y\subseteq \Sigma$ is preserved by $\Lambda'_i\times \Gamma'_i$. Since $C\cong \Gamma_i\bsl \Sigma$, the ergodicity of $\Lambda'_i$ on $\Gamma\bsl \Sigma$ implies that $C$ is an ergodic component of the $\Lambda'_i\times \Gamma'_i$-action on $\Sigma$. We know from the proof of Theorem~\ref{thm_factors} that the space $\overline{\Sigma}$ of ergodic components is an ME coupling of $\Lambda_i$ with $\Gamma_i$ for the measure projected from $m|_{\Gamma_i Y}$ (indeed, the ergodicity assumption of that proof is implied by the irreducibility of $\Lambda$ on $\Gamma\bsl \Sigma$). Now $C$ corresponds to an atom in $\overline{\Sigma}$ and thus the $\Gamma'_i$-invariant partition $\Sigma=\bigsqcup\limits_{\gamma_i\in\Gamma_i}\gamma_i C$ shows that $\Gamma_i$ acts simply transitively on $\overline{\Sigma}$.

Since $C$ is an ergodic component and all actions commute, $\lambda_i C\cap C$ can only be $C$ 
or $\varnothing$ (up to null-sets) when $\lambda_i\in \Lambda_i$. The first case implies in particular
$\lambda_i Y\subseteq C = \Gamma'_i Y$.
   But~(\ref{eq_fund_inclusion4}) gives $\lambda_i Y\subseteq \Gamma_i Y$ and thus
$\lambda_i Y\subseteq Y$. Since $\Lambda_i$ is torsion-free and acts properly on $\Sigma$, we conclude
that
$\lambda_i$ is trivial. We have shown that $\Lambda_i$ acts freely on $\overline{\Sigma}$, and also that
the sets $\lambda_i Y$ are mutually disjoint. 

A dual argument shows that all $\lambda_j Y$ are mutually disjoint when $\lambda_j$ ranges over $\Lambda_j$ 
and $1\leq j\leq n$.
We claim that this is also the case for all $\lambda Y$ with $\lambda\in \Lambda$. Indeed, assume that $\lambda Y\cap Y$ is nonnull and let us show that the component $\lambda_k$ of $\lambda$ in $\Lambda_k$ vanishes for any $k$. Write $\lambda = \lambda_k\lambda'$ with $\lambda'\in\Lambda'_k$. Then the set $\lambda' Y\cap \lambda_k^{-1} Y$ is nonnull; however, the inclusions~(\ref{eq_fund_inclusion4}) give
$$\lambda' Y\cap \lambda_k^{-1} Y\ \subseteq\ \Gamma'_k Y \cap \Gamma_k Y\ =\ Y$$
and therefore $\lambda_k Y\cap Y$ has positive measure. By what we already know, this makes $\lambda_k$ trivial as claimed.

The claim just proved forces $Y$ to be a fundamental domain for the\break $\Lambda$-action on $\Sigma$
because of the inequality~(\ref{eq_inequality_volumes}). In particular, we have\break
$[\Lambda:\Gamma]_\Sigma =1$ and
\begin{equation}
\label{eq_fund_inclusion6}%
\Lambda_j Y\ =\ \Gamma_j Y\kern1cm(1\leq j\leq n).
\end{equation}

Thus, returning now to fix some $i$, the arguments used to show that $\Gamma_i$ acts simply transitively on $\overline{\Sigma}$ can be reversed to show that $\Lambda_i$ does so too. Hence the coupling $\overline{\Sigma}$ is 
trivial (see Example~\ref{exo_trivial}) and yields an isomorphism $f_i:\Lambda_i\to \Gamma_i$. Observe that $\overline{\Sigma}$ can be viewed as the space of ergodic components of the $\Lambda'_i$-action induced on $\Gamma_i Y$ via $\Gamma_i Y\cong \Gamma'_i\bsl \Sigma$, and therefore by~(\ref{eq_fund_inclusion6}) it is the space of ergodic components of the corresponding action on $\Lambda_i Y$. These components are precisely the $\Lambda_i$-translates of $Y$ by irreducibility; we choose for $f_i$ the conjugate obtained by fixing as base point in $\overline{\Sigma}$ the component corresponding to $Y$ in the above identification. We apply the whole argument to each index $i$ and observe that with our choice of $f_i$ we have indeed
$$\Pi_{i=1}^n \big(\lambda_i f_i(\lambda_i)\big) Y\ =\ Y$$
for all $\lambda_i\in \Lambda_i$.
\hfill\qed

\Subsec{Strong rigidity with radicals} 
The main goal of this subsection is to prove
Theorem~\ref{thm_kernel_irred}, which will actually be  simpler
than in the product case (unlike the case of the superrigidity
theorem). We begin by proving a more general version of
Theorem~\ref{thm_kernel_intro}:

\demo{Proof of Theorem~{\rm \ref{thm_kernel_intro}} for the class~$\cc$ instead of~$\creg$}
Let $(\Sigma,m)$ be an ME coupling of $\Gamma$ with $\Lambda$; as
mentioned before, we may assume it is ergodic~\cite[2.2]{Furman1}.
The argument is similar to a part of the proof of
Theorem~\ref{thm_factors}, and we shall prove that the space
$(\overline{\Sigma}, \nu)$ or ergodic components of the $M\times
N$-action on $\Sigma$ is an ME coupling of $\overline{\Gamma}$
with $\overline{\Lambda}$ for an appropriate measure $\nu$. Let
$m_N$ be the measure $m$ restricted to $C=N\bsl \Sigma$ (e.g.\ via
a $\Lambda$-fundamental domain in $\Sigma$  and a section of
$\Lambda/N$) and take for $\nu$ the projection of $m_N$ under the
map $C\to \overline{\Sigma}$. As in
Theorem~\ref{thm_factors}, the ergodicity implies that one
obtains the same measure up to a scalar by proceeding with $M\bsl
\Sigma$ instead. Thus, due to the symmetry of the situation, it is
enough to show that $\overline{\Lambda}$ has a $\nu$-finite
measure  fundamental domain in $\overline{\Sigma}$.

Given a mixing unitary representation $\pi$ of $\overline{\Lambda}$ with nonvanishing $\hb^2$, we use Corollary~\ref{cor_inf} and Theorem~\ref{thm_mind} exactly as in the proof of Proposition~\ref{prop_fund_domain} to deduce that $\hb^2(\Gamma, \Mind\Sigma\Lambda\Gamma\pi)$ is nontrivial. By Proposition~\ref{prop_kernel}, the space $(\Mind\Sigma\Lambda\Gamma\pi)^M$ is nontrivial. Given a fundamental domain $X'\subseteq \Sigma$ for $\Lambda$ and the associated cocycle $\alpha: \Gamma\times X'\to \Lambda$, we obtain a nonzero measurable function $f: X'\to \H_\pi$ such that $f(\gamma\cdot x) = \pi(\alpha(\gamma, x))f(x)$ for every $\gamma\in M$ and almost every $x\in X'$. Since $\pi$ is tame as $\Lambda$-representation, we may apply cocycle reduction for the restriction $\alpha: M\times X'\to \Lambda$ in the nonergodic setting. We assumed $\overline{\Lambda}$ torsion-free and $\pi$ mixing,
and  so the only possible stabiliser in $\Lambda$ is $N$. Therefore, the reduction provides
  us with a measurable map $\fhi:X'\to \Lambda$ such that $\fhi(\gamma\cdot x)\alpha(\gamma, x)\fhi(x)^{-1}$ is in $N$ for every $\gamma\in M$ and almost every $x\in X'$. Now the new $\Lambda$-fundamental domain $X = \{\fhi(x)x : x\in X'\}$ satisfies
\begin{equation}
\label{eq_fund_kernel}%
M X\ \subseteq\ N X
\end{equation}
because for almost every $x\in X'$ and all $\gamma\in M$ we have
\begin{eqnarray*}
\gamma\fhi(x) x = \fhi(x) \gamma x &= &\fhi(x) \alpha(\gamma,x)^{-1}\gamma\cdot x\\
&= & \big(\fhi(\gamma_2\cdot x)\alpha(\gamma_2, x)\fhi(x)^{-1}\big)^{-1}(\fhi(\gamma\cdot
x)\gamma\cdot x)
\end{eqnarray*}
which is in $N X$.

We turn  back  now to $\overline{\Sigma}$. Write $\overline{X}$ for the image of $X$ in
$\overline{\Sigma}$. We claim that for every $\lambda\in \Lambda$ the sets $\lambda\overline{X}$ and
$\overline{X}$ are
$\nu$-essentially disjoint unless $\lambda\in N$. Indeed, otherwise the $N\times M$-invariant subset $N
M\lambda X \cap NM X$ of $\Sigma$ would have positive $m$-measure. Then~(\ref{eq_fund_kernel})
implies $m(N\lambda X \cap NX)>0$ and thus $\lambda\in N$, as claimed. On the other hand, the
$\Lambda$-translates of $\overline{X}$ cover $\overline{\Sigma}$.

It remains to see that $\nu(\overline{X})$ is finite. View $C$ as a $\Gamma\times \overline{\Lambda}$-space and let $X_N$ be the image of $X$ in $C$; we have $m_N(X_N) = m(X)$ since $X$ is a $\Lambda$-fundamental domain. Consider now $\overline{\Sigma}$ as the space of ergodic components of $M$ in $C$. By~(\ref{eq_fund_kernel}) $M$ preserves $X_N$ and so $\nu(\overline{X})$ is finite.
\hfill\qed

\demo{Proof of Theorem~{\rm \ref{thm_kernel_irred}}}
Let $\overline{\Sigma}$  and $X$ be as in the proof of
Theorem~\ref{thm_kernel_intro}. The subset $B=N X\subseteq \Sigma$
is $M\times N$-invariant by~(\ref{eq_fund_kernel}). Since $M$
is ergodic on $\Lambda\bsl \Sigma$, $B$ is an ergodic component
for the $M\times N$-action on $\Sigma$. It follows now from the
proof of Theorem~\ref{thm_kernel_intro} that $\overline{\Lambda}$
acts simply transitively on $\mathrm{T}=\overline{\Sigma}$. By
symmetry, we see that $\overline{\Gamma}$ acts also simply
transitively on $\mathrm{T}$ and the statement follows.
\hfill\qed

\section{Superrigidity}
\vglue-12pt

\Subsec{Preliminaries}

\begin{lemma}
\label{lemma_Z_dense}%
Let $(Z, \zeta)$ be an ME coupling of countable groups $\Gamma,\Lambda$. Let $E\subseteq Z$ be a fundamental domain for $\Gamma$ and $\teta: \Lambda\times E\to \Gamma$ the associated cocycle. If $\teta$ is equivalent to a cocycle ranging $\zeta$-essentially in a subgroup
 $\Gamma_0<\Gamma${\rm ,} then $\Gamma_0$ has finite index in $\Gamma$.
\end{lemma}

\Proof
Let $\chi:Z\to \Gamma$ be the $\Gamma$-equivariant retraction associated to $F$ and assume that there is $\fhi: F\to \Gamma$ such that the cocycle $\teta'$ defined by
$$\teta'(\lambda, x) = \fhi(\lambda\cdot x)\teta(\lambda, x)\fhi(x)^{-1}$$
ranges in $\Gamma_0$. Define $\widetilde{\fhi}:Z\to \Gamma/\Gamma_0$ by
$$\widetilde{\fhi}(x) = \chi(x)\fhi(\chi(x)^{-1}x)^{-1}\Gamma_0;$$
this is $\Gamma$-equivariant. For $x\in F$, $\lambda\in\Lambda$ and $\gamma\in\Gamma$ one has
\begin{align*}
\widetilde{\fhi}(\lambda\gamma x) &= \gamma\chi(\lambda x)\fhi(\chi(\lambda x)^{-1}\lambda x)^{-1}\Gamma_0 = \gamma\chi(\lambda x)\fhi(\lambda\cdot x)^{-1}\Gamma_0 \\
&= \gamma\chi(\lambda x)\fhi(\lambda\cdot x)^{-1}\teta'(\lambda, x)^{-1}\Gamma_0 = \gamma\chi(\lambda x)\teta(\lambda,x)\fhi(x)^{-1}\Gamma_0.
\end{align*}
Since $x\in F$, we have $\teta(\lambda, x) = \chi(\lambda x)^{-1}$ and the above reduces to $\gamma\fhi(x)\Gamma_0 =\widetilde{\fhi}(\gamma x)$. Therefore we obtain a $\Gamma$-equivariant map $\Lambda\bsl Z\to \Gamma/\Gamma_0$. Projecting the $\Gamma$-invariant finite measure to $\Gamma/\Gamma_0$, we deduce that the latter is finite.
\hfill\qed

\begin{lemma}[Furman's homomorphism]
\label{lemma_Furman}%
 Let $\Lambda$ be a countable
group acting on a measure space $(\Sigma,m)${\rm ,} preserving the
measure class. Let $F:\Sigma\times \Sigma\to \Gamma$ be a
measurable map to a countable group $\Gamma$ and assume that for
all $\lambda\in\Lambda$
\begin{itemize}
\item[{\rm (i)}] $F(x,y) = F(\lambda x,\lambda y)$\hfill for $m^2$-almost all $(x,y)\in \Sigma^2$.

\item[{\rm (ii)}] $F(\lambda x,y)F(x,y)^{-1} = F(\lambda x,z)F(x,z)^{-1}$ 
 \hfill for $m^3$-almost all $(x,y,z)\in \Sigma^3$.

\item[\/{\rm (iii)}\/] $F(x,y) = F(y,x)^{-1}$\hfill for $m^2$-almost all $(x,y)\in \Sigma^2$.

\item[{\rm (iv)}] $F(x,y)F(y,z) = F(x,z)$\hfill for $m^3$-almost all $(x,y,z)\in \Sigma^3$.
\end{itemize}
Then for $m$-almost every $x\in\Sigma$ the map
$$\ro_x: \Lambda\lra \Gamma,\ \ \ro_x(\lambda) = F(\lambda^{-1}x, y)F(x,y)^{-1}$$
is a homomorphism{\rm ,} and for $m^2$-almost all $(x,y)\in \Sigma^2$ the homomorphisms
 $\ro_x$ and $\ro_y$ are conjugated\/{\rm :}\/
$$\ro_y(\lambda) = F(x,y)^{-1}\ro_x(\lambda)F(x,y).\eqno{(\forall\,\lambda\in\Lambda)}$$
\end{lemma}

\Proof
This follows from Lemma~5.2 in~\cite{Furman1} upon modification of notation. \phantom{overthere}
\hfill\qed

\begin{defi}
\label{defi_ICC}%
A group is said to be {ICC} (for \emph{infinite conjugacy
classes}) if the conjugacy class of every nontrivial element is
infinite.
\end{defi}

We shall see (Proposition~\ref{prop_ICC_C}) that this property is automatically satisfied for torsion-free groups in~$\cc$. The following is a straightforward verification:

\begin{lemma}
\label{lemma_Dirac_inv}%
If $\Gamma$ is an {\rm ICC} group{\rm ,}
 then the Dirac mass $\delta_e$ at the neutral element is the only probability measure on $\Gamma$ which is
invariant for the $\Gamma$-action by conjugation.\hfill\qed
\end{lemma}

The following lemma isolates the use of the mild mixing assumption
in the results where it appears. We recall that we have introduced
in Section~\ref{sec_ME_WOE} the notation $A^1$, $A^2$ in order to
distinguish between two commuting $\Gamma$-actions.

\begin{lemma}
\label{lemma_mixing_coupling}%
Let $(\Sigma,m)$ be an ME coupling  of two countable groups
$\Gamma,\Lambda${\rm ,} let $\Delta<\Gamma$ be a subgroup and let
$\Omega = \Sigma\times_\Lambda \Lambda\times_\Lambda \check\Sigma$
be the composed self-$\Gamma$-coupling \/{\rm (}\/see Definition~{\rm \ref{defi_composition}}
and the paragraph preceding it\/{\rm ).}\/ If the $\Lambda$-action on $\Gamma\bsl \Sigma$
 is mildly mixing and the $\Delta$-action on $\Lambda\bsl \Sigma$ is ergodic{\rm ,} then the
$A^i(\Delta)$-action on $A^j(\Gamma)\bsl\Omega$ is ergodic for
all $i\neq j$.
\end{lemma}

\Proof
By~\cite{Schmidt-Walters} (a reference for which we thank Eli
Glasner), the assumption on the $\Lambda$-action on $\Gamma\bsl
\Sigma$ is equivalent to the following: for any standard Borel
$\Lambda$-space $(S,\nu)$ with a nonatomic invariant ergodic
$\sigma$-finite measure $\nu$, the diagonal $\Lambda$-action
on $\Gamma\bsl \Sigma\times S$ is ergodic. Assume now that we have
a real-valued $A^i(\Delta)$-invariant measurable function $F$
on $A^j(\Gamma)\bsl\Omega$. In view of $\Omega\cong
\Sigma\times_\Lambda \check\Sigma$, the function $F$ corresponds
to a measurable function $F'$ on $\Sigma\times \check\Sigma$ that
is $\Lambda$-invariant for the diagonal action and
$A^i(\Delta)\times A^j(\Gamma)$-invariant. View $F'$ as a
$\Lambda$-invariant function on $\Gamma\bsl \Sigma\times
\Delta\bsl\Sigma$; the ergodicity of the $\Delta$-action on
$\Lambda\bsl \Sigma$ implies that $S=\Delta\bsl\Sigma$ is an
ergodic $\Lambda$-space (notice that $\Delta\bsl\Sigma$ is not
atomic since otherwise $\Sigma$ and then $\Gamma\bsl\Sigma$ would
be, contradicting the mild mixing assumption). Therefore, the
above criterion for mild mixing shows that $F'$ is constant.
\hfill\qed

\Subsec{Proof of Theorems~{\rm \ref{thm_factors_mystery_new}}
and~{\rm \ref{thm_factors_mystery_new}*}} Let $(\Omega,\omega)$ be the composed ME coupling of
$\Gamma$ with itself defined by $\Omega =
\Sigma\times_\Lambda \Lambda\times_\Lambda \check\Sigma$. The class of $(x,\lambda,y)\in
\Sigma\times \Lambda\times \check\Sigma$ in $\Omega$ will be denoted by $[x,\lambda,y]$.

By Lemma~\ref{lemma_mixing_coupling} applied to $\Delta=\Gamma_j$ for each $1\leq j\leq n$, the coupling $\Omega$ of $\Gamma$ with itself satisfies the assumptions of Theorem~\ref{thm_factors_irred} point~(ii). Thus there
are an automorphism $f$ of $\Gamma$ and a factor map $\Phi:\Omega\to  \mathrm{T}$ to a trivial coupling
$\mathbf{T}$ such that for every $\gamma, \gamma'\in \Gamma$ and a.e.\ $x\in \Omega$ one has
$$\Phi(A^1_\gamma A^2_{\gamma'} x) = \gamma f(\gamma')\Phi(x).$$
We may identify $\Gamma$ with $\mathrm{T}$ as in Example~\ref{exo_trivial} so that we consider now
$\Phi$ as a map $\Omega\to \Gamma$ satisfying
\begin{equation}
\label{eq_Phi}%
\Phi(A^1_\gamma A^2_{\gamma'} x) = \gamma \Phi(x)f(\gamma')^{-1}.
\end{equation}
We prove now an analogue of A.~Furman's  crucial Lemma~5.5
in~\cite{Furman1} using a similar approach, but bypassing the
seemingly nontrivial question of tameness of the conjugation
action of a group on the space of probability measures on it.
We point out that our approach also applies to Furman's setting
and simplifies the argument therein.

\begin{lemma}
\label{lemma_not_depend}%
The map $\Psi: \Sigma\times\Sigma\times\Sigma \to \Gamma$ defined $m^3$-a.e.\ by
\begin{equation}
\label{eq_Psi}%
\Psi(x,y,z) = \Phi([x,e,z])\Phi([y,e,z])^{-1}
\end{equation}
does not depend on $z$ in the sense that $\Psi(x,y,z_1)=\Psi(x,y,z_2)$ holds for\break $m^4$-almost every
$(x,y,z_1,z_2)\in\Sigma^4$.
\end{lemma}

{\it Proof}.
Define the map ${\mathscr T}: \Sigma^4\to \Gamma$ by
\begin{equation}
\label{eq_def_T}%
{\mathscr T}(x,y,z_1, z_2) = \Psi(x,y,z_1)\Psi(x,y,z_2)^{-1}.
\end{equation}
The above definitions imply that for $m^4$-a.e.\ $(x,y,z_1,z_2)\in\Sigma^4$ and every $\lambda\in\Lambda$, $\gamma\in\Gamma$ one has
\begin{align*}
{\mathscr T}(x,y,z_1, z_2) &= {\mathscr T}(\lambda x,\lambda y,\lambda z_1, \lambda z_2)\\
&= {\mathscr T}(x,\gamma y,z_1, z_2)\\
&= {\mathscr T}(x,y,\gamma z_1, z_2)\\
&= {\mathscr T}(x,y,z_1, \gamma z_2).
\end{align*}
Consider the $\Gamma\times \Lambda$-space $Z= \Sigma\times(\Gamma\bsl\Sigma)^3$, where $\Gamma$ acts on the first factor and $\Lambda$ by diagonal (fourtuple) action. This is again a coupling and ${\mathscr T}$ induces a map ${\mathscr T}_0:\Lambda\bsl Z\to \Gamma$. Substituting~(\ref{eq_Phi}) in~(\ref{eq_Psi}) and then in~(\ref{eq_def_T}), one gets
$${\mathscr T}(\gamma x,y,z_1, z_2) = \gamma {\mathscr T}(x,y,z_1, z_2)\gamma^{-1},$$
so that ${\mathscr T}_0$ is $\Gamma$-equivariant for the conjugating action on $\Gamma$. If we project
now the $\Gamma$-invariant measure of $\Lambda\bsl Z$ through ${\mathscr T}_0$, we get a conjugating
invariant probability measure on $\Gamma$. Since the product of two ICC groups is again ICC, the
(independent) Proposition~\ref{prop_ICC_C} ensures that $\Gamma$ is ICC, and so
Lemma~\ref{lemma_Dirac_inv} shows that ${\mathscr T}_0$ is essentially constant with value
$e\in\Gamma$, proving Lemma~\ref{lemma_not_depend}.
\Endproof \vskip4pt 

We may now define $F:\Sigma\times \Sigma\to \Gamma$ as in Lemma~\ref{lemma_Furman} by $F(x,y)=\Psi(x,y,z)$. All properties listed in Lemma~\ref{lemma_Furman} are readily verified, and we obtain thus a family of conjugated homomorphisms $\ro_x: \Lambda\to \Gamma$. In particular, we obtain a
well-defined subgroup $N\lhd \Lambda$ as the kernel of almost all $\ro_x$. For later reference, we record
that the definition of $\ro_x$ boils down to
\begin{equation}
\label{eq_def_ro}%
\ro_x(\lambda) = \Phi([x,\lambda,z])\Phi([x,e,z])^{-1}\ \ \forall\,\lambda\in\Lambda,\mbox{\ $m^2$-a.e.\ $(x,z)\in \Sigma^2$.}
\end{equation}

Let $D\subseteq \Sigma$ be a fundamental domain for the $\Lambda$-action on $\Sigma$ and consider the measure space
$$\widetilde\Omega = D\times\Lambda\times D\ \subseteq\ \Sigma\times\Lambda\times\Sigma.$$
This inclusion yields an isomorphism of measure spaces $\widetilde\Omega\cong \Omega$ and through this identification we let $\Gamma$ act on $\widetilde\Omega$ by the second $\Gamma$-action $A^2$ on $\Omega$. We also endow $\widetilde\Omega$ with the $\Lambda$-action coming from left multiplication on itself, so that we obtain on $\widetilde\Omega$ a $\Gamma\times\Lambda$-structure given explicitly by
\begin{equation}
\label{eq_G_L_struct}%
(\gamma,\lambda)\,(x,\lambda_1,y) = (x,\lambda\lambda_1 \alpha(\gamma,y)^{-1},\gamma\cdot y),
\end{equation}
where $\alpha: \Gamma\times D\to \Lambda$ is the cocycle corresponding to $D$; we recall also that the $\Gamma$-action on $D$ is given by $\gamma\cdot y = \alpha(\gamma,y)\gamma y$. If we denote by $\widetilde\Phi$ the map $\widetilde\Omega\to \Gamma$ induced by $\Phi$ under the identification $\widetilde\Omega\cong \Omega$, then $E=\widetilde\Phi^{-1}(e)\subseteq \widetilde\Omega$ is a finite measure fundamental domain for the $\Gamma$-action since $\widetilde\Phi:\widetilde\Omega\to \Gamma$ is $\Gamma$-equivariant (with respect to right multiplication on $\Gamma$ through the automorphism $f$). Thus $\widetilde\Omega$ is an ME coupling of $\Gamma$ with $\Lambda$ because the latter has an obvious finite measure fundamental domain. Equation~(\ref{eq_def_ro}) gives
\begin{equation}
\label{eq_def_ro2}%
\ro_x(\lambda) = \widetilde\Phi(x,\lambda,z) \widetilde\Phi(x,e,z)^{-1}\ \ \forall\,\lambda\in\Lambda,\mbox{\ $m^2$-a.e.\ $(x,z)\in D^2$,}
\end{equation}
but it shows further that for $\lambda_0\in\Lambda$
\begin{equation}
\label{eq_kernel}%
\begin{array}{ll}
\lambda_0\in N \Longleftrightarrow &\forall\ \lambda_1,\lambda_2\in\Lambda, \mbox{\ for a.e.\ $(x,y)\in D^2$:}\\
&\widetilde\Phi(x,\lambda_1\lambda_0\lambda_2, y)= \widetilde\Phi(x,\lambda_1\lambda_2, y).
\end{array}
\end{equation}
This characterisation shows in particular (for $\lambda_1=e$) that $N$ preserves $E$, and thus, by properness of the $\Lambda$-action, the group $N$ is finite. This is the only point where we use
the fact   that $\widetilde\Phi$ has finite measure fibres.

By~(\ref{eq_def_ro2}), we have for a.e.\ $t=(x,\lambda_1,y)$ in $\widetilde\Omega$ and every 
$\lambda\in \Lambda$,
\begin{equation}
\label{eq_ro_phi}%
\begin{array}{ll}
\ro_x(\lambda) &= \ro_x(\lambda\lambda_1)\ro_x(\lambda_1)^{-1}\\
&= \widetilde\Phi(x,\lambda\lambda_1,z) \widetilde\Phi(x,e,z)^{-1} \Big(\widetilde\Phi(x,\lambda_1,z) \widetilde\Phi(x,e,z)^{-1}\Big)^{-1}\\
&= \widetilde\Phi(\lambda t)\widetilde\Phi(t)^{-1}.
\end{array}
\end{equation}
This is just $\widetilde\Phi(\lambda t)$ whenever $t\in E$. On the other hand, if $\teta: \Lambda\times E\to \Gamma$ is the cocycle associated to $E$, we have
$$e = \widetilde\Phi(\teta(\lambda,t)\lambda t) = \widetilde\Phi(\lambda t)f(\teta(\lambda,t))^{-1}.$$
We deduce
\begin{equation}
\label{eq_ro_teta}%
\ro_x(\lambda) = f(\teta(\lambda,t))\kern1cm \forall\,\lambda\in\Lambda,\ {\rm a.e.}\ t=(x,\lambda_1,y)\in E.
\end{equation}
Applying the Fubini-Lebesgue theorem to the conclusion of Lemma~\ref{lemma_Furman}, we have some $x_0\in D$ such that $\ro = \ro_{x_0}: \Lambda\to \Gamma$ is a homomorphism with kernel $N$ and such that for a.e.\ $x\in D$
\begin{equation}
\label{eq_ro_conj}%
\ro_x(\lambda) = F(x_0,x)^{-1}\ro(\lambda)F(x_0,x) \kern1cm \forall\,\lambda\in\Lambda.
\end{equation}
We now proceed to show that $\ro(\Lambda)$ has finite index in $\Gamma$. Define $\fhi: E\to \Gamma$ by
$$\fhi(t) = f^{-1}\big(F(x_0, x)\big)\ \mbox{\ for $t=(x,\lambda_1,y)\in E$}.$$
Observe that $\fhi(\lambda\cdot t) = \fhi(t)$ for $\lambda\in\Lambda$ because by~(\ref{eq_G_L_struct})
$$\lambda\cdot t = \lambda \teta(\lambda,t)t = \Big(x,\lambda\lambda_1\alpha(\teta(\lambda,t),y)^{-1}, \teta(\lambda,t)\cdot y\Big).$$
Consider the cohomologous cocycle
$$\teta'(\lambda, t) = \fhi(\lambda\cdot t)\teta(\lambda, t)\fhi(t)^{-1}.$$
Applying successively~(\ref{eq_ro_teta}) and~(\ref{eq_ro_conj}), we have a.e.\ %
$$
\begin{array}{ll}
f(\teta'(\lambda,t)) &= F(x_0, x)f(\teta(\lambda,t))F(x_0,x)^{-1}\\
&= F(x_0, x)\ro_x(\lambda)F(x_0,x)^{-1} = \ro(\lambda).
\end{array}
$$
In other words, $\teta'$ ranges essentially in $f^{-1}(\ro(\Lambda))$. Applying Lemma~\ref{lemma_Z_dense}, we deduce that $f^{-1}(\ro(\Lambda))$, and thus also $\ro(\Lambda)$, have finite index in $\Gamma$. Summing up, we obtain 
a virtual isomorphism of $\Lambda$ and $\Gamma$:
\begin{equation}
\label{eq_extension_ro}%
1\lra N\lra \Lambda \xrightarrow{\ \ro\ } \ro(\Lambda) \lra 1.
\end{equation}
We proceed now to construct the virtual isomorphism~(\ref{eq_extension}) with the additional properties stated in Theorem~\ref{thm_factors_mystery_new} by considering again the ME coupling $(\Sigma,m)$ of $\Gamma,\Lambda$ but now in view of~(\ref{eq_extension_ro}). Endow $N\bsl\Sigma$ with the restricted measure $m_N$ (alternatively, we could work with the quotient measure $|N|\cdot m_N$); we consider this as a $\ro(\Lambda)\times \Gamma$-space via $\ro$. Since $\Gamma$ is torsion-free, the $N$-action on $\Gamma\bsl \Sigma$ is essentially free and thus the $\Gamma$-action on $N\bsl\Sigma$ still admits a fundamental domain (of measure $m(\Gamma\bsl\Sigma)/|N|$). Therefore, $(N\bsl\Sigma,m_N)$ is an ME coupling of $\ro(\Lambda)$ with
$\Gamma$, with coupling index satisfying
\begin{equation}
\label{eq_coupling_N}%
[\Gamma:\ro(\Lambda)]_{N\bsl\Sigma}=|N|\cdot[\Gamma:\Lambda]_\Sigma.
\end{equation}
We write now $\Gamma^*_i$ for the image of the projection of $\ro(\Lambda)$ to $\Gamma_i$ and consider the finite index subgroup $\Gamma^*<\Gamma$ containing $\ro(\Lambda)$ defined by $\Gamma^*= \Gamma^*_1\times\cdots\times\Gamma^*_n$. Let $(\Sigma^*,m^*)$ be the $\Gamma^*\times \Gamma$-space obtained by suspension from $(N\bsl\Sigma,m_N)$; that is, we have
$$\Sigma^* = \ro(\Gamma)\bsl\Big(N\bsl\Sigma\times \Gamma^*\Big),$$
wherein $\ro(\Lambda)$ acts diagonally, while the $\Gamma^*$-action is given by (inverted) right multiplication on the second factor and $\Gamma$ acts on the first factor. Then $\Sigma^*$ yields an ME coupling of $\Gamma^*$ with $\Gamma$ and we have
\begin{equation}
\label{eq_coupling_induced}%
[\Gamma:\Gamma^*]_{\Sigma^*}\cdot[\Gamma^*:\ro(\Lambda)]=[\Gamma:\ro(\Lambda)]_{N\bsl\Sigma}.
\end{equation}
The $\Gamma$-action on $\Gamma^*\bsl\Sigma^*$ is still irreducible by construction. On the other hand, we claim that the $\Gamma^*$-action on $\Gamma\bsl\Sigma^*$ is irreducible. Indeed, the ergodicity of the $\Gamma^*_i$-action on $\Gamma\bsl\Sigma^*$ is equivalent to the ergodicity of $\ro(\Lambda)\cap\Gamma_i$ on $\Gamma\bsl(N\bsl\Sigma)$, which in turn is equivalent to the ergodicity of $\ro^{-1}(\Gamma_i)$ on $\Gamma\bsl \Sigma$. The latter follows from mild mixing since $\ro^{-1}(\Gamma_i)$ is infinite, proving the claim.

At this point we are in position to apply
Theorem~\ref{thm_factors_irred} point~(ii) to $\Sigma^*$ and
obtain isomorphisms $h_i:\Gamma^*_i\to \Gamma_i$ after possibly
permuting factors. If we set $h=\Pi_{i=1}^n h_i$, $\pi =
h\circ\ro$ and $\Gamma'=\pi(\Lambda)$, then we have a virtual
isomorphism~(\ref{eq_extension}) such that the projections of
$\Gamma'$ to each $\Gamma_i$ are onto. Moreover,
Theorem~\ref{thm_factors_irred} yields
$[\Gamma:\Gamma^*]_{\Sigma^*}=1$, which in view
of~(\ref{eq_coupling_N}) and~(\ref{eq_coupling_induced}) implies
indeed the formula~(\ref{eq_index_formula}) since
$[\Gamma:\Gamma']=[\Gamma^*:\ro(\Lambda)]$. This completes the
proof of Theorem~\ref{thm_factors_mystery_new}. 
Theorem~\ref{thm_factors_mystery_new}* then follows from the
properties of the fundamental domain in $\Sigma^*$ granted by
Theorem~\ref{thm_factors_irred} in view of the construction of
$\Sigma^*$.

\setcounter{Subsec}{3}

\vskip6pt 6.3.  {\it Proof of Theorem~{\rm \ref{thm_kernel_mystery_new}}}.
We will   start by constructing a homomorphism $\ro$ of
$\Lambda$, but in contrast to the proof of
Theorem~\ref{thm_factors_mystery_new}, $\ro$ will range in
$\overline{\Gamma}$ instead of $\Gamma$, and we shall show that
its kernel is amenable rather than finite. Moreover, in order to
show that it has finite index image, we shall need the following
variation of Lemma~\ref{lemma_Z_dense} (which will be complemented
by Lemma~\ref{lemma_pseudo_cohomologous} below):

\vglue-16pt
\phantom{up}

\begin{lemma}
\label{lemma_Z_dense2}%
Let $(Z, \zeta)$ be an {\rm ME}
 coupling of countable groups $\Gamma,\Lambda$ and let $M\lhd \Gamma$ be a normal subgroup. Suppose
that there is a measurable $M$-invariant set $E\subseteq Z$ such that $Z = \bigsqcup\limits_{\gamma M\in
\Gamma/M} \gamma E$. Pick a measurable map $\sigma: Z\to \Gamma$ such that $x\in\sigma(x)E$ for all
$x\in Z${\rm ,} so that in particular the induced map $\bar\sigma:Z\to\overline{\Gamma}=\Gamma/M$ is
$\Gamma$-equivariant and determined by $E$. If $\bar \sigma|_{\Lambda E}$ ranges $\zeta$-essentially in
a subgroup $\overline{\Gamma}_0<\overline{\Gamma}$, then $\overline{\Gamma}_0$ has finite index in
$\overline{\Gamma}$.
\end{lemma}

\phantom{up}
\vglue-16pt

{\it Proof}.
Denote by $\Gamma_0$ the pre-image of $\overline{\Gamma}_0$ in $\Gamma$, so that $\sigma$ ranges essentially in $\Gamma_0$. It is enough to show that $\Gamma_0$ has finite index in $\Gamma$. Let $F\subseteq Z$ be a fundamental domain for the $\Gamma$-action and define $F' = \{\sigma(x)^{-1}x : x\in F\}$. One checks that $F'$ is also a fundamental domain. Moreover, $F'\subseteq E$. Indeed, for $x\in F$ and $y=\sigma(x)^{-1}x$ we have $\bar\sigma(y) = \bar\sigma(x)^{-1}\bar\sigma(x)$ so that $\sigma(y)\in M$; now $y\in\sigma(y)E$ implies $y\in E$.

Let now $\chi, \chi'$ be the retractions associated to $F,F'$. One computes from the definition of $F'$ that $\chi'(x) = \chi(x)\sigma(\chi(x)^{-1}x)$ for a.e.\ $x\in Z$. Thus $\chi'(x)M = \bar\sigma(x)$ and in particular $\chi'(\lambda x)M\in\overline{\Gamma}_0$ for $\zeta$-a.e.\ $x\in E$ and all $\lambda\in\Lambda$. Therefore the cocycle associated to $F'$ ranges essentially in $\Gamma_0$. An application of Lemma~\ref{lemma_Z_dense} completes the proof.
\hfill\qed

\demo{Proof of Theorem~{\rm \ref{thm_kernel_mystery_new}}}
We write $\gamma\mapsto \bar\gamma = \gamma M$ for the natural map. We consider again the composed self-coupling $\Omega = \Sigma\times_\Lambda \Lambda\times_\Lambda \check\Sigma$ of $\Gamma$ and apply Lemma~\ref{lemma_mixing_coupling} to $\Delta = M$. This time we may apply Theorem~\ref{thm_kernel_irred} and deduce an automorphism $f$ of $\overline{\Gamma}$ and a factor map $\Phi:\Omega\to  \mathrm{T}$ to a trivial self-coupling $\mathrm{T}\cong \overline{\Gamma}$ of $\overline{\Gamma}$, and we have
\begin{equation*}
\Phi(A^1_\gamma A^2_{\gamma'} x) = \bar\gamma \Phi(x)f(\bar\gamma')^{-1}.
\end{equation*}
Since $\overline{\Gamma}$ is ICC by Proposition~\ref{prop_ICC_C}, we can argue exactly as for Theorem~\ref{thm_factors_mystery_new} and deduce that the map $\Psi: \Sigma\times\Sigma\times\Sigma \to \overline{\Gamma}$ defined as in Lemma~\ref{lemma_not_depend} does not depend on the last variable. Taking the corresponding definitions for $F$, $\ro_x:\Lambda\to\overline{\Gamma}$, the $\Gamma\times\Lambda$-space $\widetilde\Omega = D\times\Lambda\times D$ and $\widetilde\Phi$,
we need not change   the arguments to apply Furman's Lemma~\ref{lemma_Furman} and
get~(\ref{eq_def_ro2}) and~(\ref{eq_kernel}).

An important difference, though, is that $E=\widetilde\Phi^{-1}(e)$ is in general not a fundamental domain and may have infinite measure. However, $M$ preserves $E$ and we have
\begin{equation}
\label{eq_E}%
\widetilde\Omega = \bigsqcup_{\gamma M\in \Gamma/M} \gamma E.
\end{equation}
\begin{lemma}
\label{lema_N_amenable}%
The group $N$ is amenable.
\end{lemma}

{\it Proof}.
Equation~(\ref{eq_kernel}) shows  that $N$ preserves $E$, which is
thus a $M\times N$-space. Both actions are proper, measure-preserving and have some fundamental domain
since they are the restriction of the actions on $\widetilde \Omega$. Therefore we
can conclude by Corollary~\ref{cor_induce_amenablility}, provided
$M$ has a \emph{finite measure} fundamental domain in $E$. Take a
fundamental domain $F'\subseteq E$ for the $\Gamma$-action on
$\widetilde\Omega$ as constructed in the proof of
Lemma~\ref{lemma_Z_dense2}. It remains to show that the $M$-translates of $F'$ cover $E$. Since $\Gamma F' =\widetilde
\Omega$, this follows from the fact that the union~(\ref{eq_E}) is
disjoint.
\Endproof \vskip4pt 

We fix a map $\sigma: \widetilde \Omega\to \Gamma$ as in Lemma~\ref{lemma_Z_dense2} and we find,
as after~(\ref{eq_ro_phi}), that $\ro_x(\lambda)= \widetilde\Phi(\lambda t)$ whenever $t=(x,\lambda_1,y)\in
E$. We have
$$e = \widetilde\Phi(\sigma(\lambda t)^{-1}\lambda t) = \widetilde\Phi(\lambda t)f(\bar\sigma(\lambda t))$$
whence again
\begin{equation}
\label{eq_ro_teta2}%
\ro_x(\lambda) = f(\bar\sigma(\lambda t))^{-1}\kern1cm \forall\,\lambda\in\Lambda,\ {\rm a.e.}\ t=(x,\lambda_1,y)\in E.
\end{equation}
Fix $x_0\in D$ with~(\ref{eq_ro_conj}) and let $\ro = \ro_{x_0}$, $\overline{\Gamma}_0 = \ro(\Lambda)$ so that with~(\ref{eq_ro_conj}) and~(\ref{eq_ro_teta2}) we get for a.e.\ $t\in E$
\begin{equation}
\label{eq_image_sigme}%
\ro(\lambda)^{-1} = F(x_0,x)f(\bar\sigma(\lambda t))F(x_0,x)^{-1}.
\end{equation}
Pick an arbitrary section $\tau:\overline{\Gamma}\to \Gamma$ of the projection and define $\psi:\widetilde\Omega\to \Gamma$ by $\psi(x,\lambda_1,y) = \tau f^{-1} F(x_0,x)$. Note that $\psi$ is $\Gamma\times \Lambda$-invariant by~(\ref{eq_G_L_struct}). Define further $\sigma':\widetilde\Omega\to \Gamma$ by $\sigma'(t) = \psi(t)\sigma(t)\psi(t)^{-1}$ and $E'=\{\psi(t)t:t\in E\}$. Now~(\ref{eq_image_sigme}) implies
\begin{equation}
\label{eq_pseudo_cohomologous}%
\overline{\sigma'}(\lambda t)\in f^{-1}(\overline{\Gamma}_0)\kern1cm \forall\,\lambda\in\Lambda,\ {\rm a.e.}\ t \in E.
\end{equation}
Define $\sigma'':\widetilde\Omega\to \Gamma$ by $\sigma''(s) = \sigma'(\psi(s)^{-1}s)$. We replace cocycle equivalence by:

\begin{lemma}
\label{lemma_pseudo_cohomologous}%
$E'$ and $\sigma''$ satisfy the assumptions of Lemma~{\rm \ref{lemma_Z_dense2}}
 for $Z=\widetilde\Omega$ and $f^{-1}(\overline{\Gamma}_0)$ instead of $\overline{\Gamma}_0$.
\end{lemma}

{\it Proof}.
First we check $ME'=E'$. Let $m\in M$, $t\in E$, $s=\psi(t)t\in E'$ and $m'=\psi(t)^{-1}m\psi(t)\in M$. The $M$-invariance of $\psi$ gives
$$m s=m\psi(t)t = \psi(t)m't = \psi(m't)m't$$
which is in $E'$ since $ME = E$. One checks in a similar way that $\widetilde\Omega = \bigsqcup\limits_{\gamma M\in \Gamma/M} \gamma E'$. We check now $s\in \sigma''(s)E'$ for all $s\in\widetilde\Omega$. Setting $t=\psi(s)^{-1}s$ we have $s=\psi(t)t$ since $\psi(t)=\psi(s)$. Write
$$t = \sigma(t)\sigma(t)^{-1}t = \psi(t)^{-1}\sigma'(t)\psi(t)\sigma(t)^{-1}t,$$
so that
$$s = \sigma'(t)\psi(t)\sigma(t)^{-1}t = \sigma'(t)\psi(\sigma(t)^{-1}t)\sigma(t)^{-1}t.$$
By the choice of $\sigma$ we have $\sigma(t)^{-1}t\in E$ so the above shows that $s$ is in $\sigma'(t)E' = \sigma''(s)E'$ as claimed. Finally, for $t\in E$ and $s=\psi(t)t\in E'$ we have indeed
$$\overline{\sigma''}(\lambda s) = \overline{\sigma'}(\lambda\psi(s)^{-1}s) \overline{\sigma'}(\lambda t)$$
since $t=\psi(s)^{-1}s$ and thus $\overline{\sigma''}|_{\Lambda
E'}$  ranges essentially in $f^{-1}(\overline{\Gamma}_0)$
by~(\ref{eq_pseudo_cohomologous}).
\Endproof \vskip4pt 

In conclusion, Lemma~\ref{lemma_Z_dense2} forces
$f^{-1}(\overline{\Gamma}_0)$ and hence also $\overline{\Gamma}_0$
to have finite index in $\overline{\Gamma}$. Writing
$\overline{\Lambda}=\Lambda/N\cong \overline{\Gamma}_0$, we
observe that this group is torsion-free. It
is also in~$\cc$ since this property passes to finite index
subgroups (Lemma~\ref{lemma_cofoelner_in_C} below). Therefore, in
order to complete the proof of
Theorem~\ref{thm_kernel_mystery_new} by an application of
Theorem~\ref{thm_kernel_irred}, it remains only to check that $N$
is ergodic on $\Gamma\bsl \Sigma$. Since the $\Lambda$-action
on $\Gamma\bsl \Sigma$ is mildly mixing, this is guaranteed if $N$
is infinite. On the other hand, we claim that if $N$ were finite
then $M$ would be so too and therefore the ergodicity assumption
for the $M$-action on $\Lambda\bsl\Sigma$ would make $\Sigma$
purely atomic, in which case Theorem~\ref{thm_kernel_mystery_new}
holds trivially.

Finally, to verify the claim above,  observe that with $N$ finite we deduce that
$\Lambda$ is in~$\cc$ (Lemma~\ref{lemma_finite_kernel} below), so
that by Corollary~\ref{cor_cc_ME_invariants} below we deduce that
$\Gamma$ is in~$\cc$ as well. Since $M$ is amenable and normal, by Proposition~\ref{prop_kernel}
it must fix a nonzero vector in any $\pi$ with $H^2_b(\Gamma, \pi) \neq 0$.
In our case such $\pi$ which is mixing exists, hence $M$ is finite.
\hfill\qed

\section{Groups in the class~$\cc$ and ME invariants}
\label{sec_in_c}%

In this section, which is independent of the rest of the paper, we
collect some information about the class~$\cc$ of groups to which our
results apply.

Recall that we defined the class~$\cc$ as the collection of all countable groups $\Gamma$
admitting a mixing unitary representation $\pi$ such that $\hb^2(\Gamma,\pi)$ is nonzero,
and that~$\creg$ is the subclass of those for which one can take $\pi$ to be the regular
representation on $\ell^2(\Gamma)$. We do not know whether or not the inclusion $\creg\subseteq \cc$
 is strict; it so happens that all our proofs that certain geometrically defined classes of
groups belong to~$\cc$ actually show that they belong to~$\creg$.

We begin by justifying that the classes of groups listed in Examples~\ref{exos_c} do indeed belong to the class~$\creg$. This is established in the companion paper~\cite{Monod-Shalom1} using also~\cite{Mineyev-Monod-Shalom}; we refer to~\cite{Monod-Shalom1} for the geometric background and context. However, let us emphasize that we shall present at the end of this section an alternative argument to show that many groups in the list of Examples~\ref{exos_c} are in~$\cc$, without appealing to~\cite{Monod-Shalom1} (or~\cite{Mineyev-Monod-Shalom})

\demo{Proof of Theorem~{\rm \ref{thm_exos_c}}}
Let $\Gamma$ be a group as in Examples~\ref{exos_c}. In case~(i), Corollary~7.8 of~\cite{Monod-Shalom1} states that $\Gamma$ belongs to~$\creg$ (note that we may assume the tree is countable since $\Gamma$ is). In case~(ii), Corollary~7.6 of~\cite{Monod-Shalom1} ensures $\Gamma \in \creg$. Case~(iii) is Theorem~3 in~\cite{Mineyev-Monod-Shalom}.
\Endproof \vskip4pt 

Via basic Bass-Serre theory, it follows that any nonelementary free product of groups is in~$\creg$, and that this is more generally so if one amalgamates over a finite subgroup (we recall here that an amalgamated product $A*_CB$ is \emph{nonelementary} if $A\neq C$ and $[B\!:\!C]>2$, or vice-versa). We have observed in~7.10, 7.11 of~\cite{Monod-Shalom1} that even if $C$ is infinite, $A*_CB$ is in $\creg$ as soon as $C$ is malnormal (or almost malnormal) in one factor, and that there is in fact a sequence of weaker and weaker acylindricality conditions generalising this fact. Already the case of free products has the following immediate consequence:

\begin{cor}
\label{cor_continuum}%
 {\rm (i)}~Every countable group embeds into a group in $\creg$.

\nobreak
{\rm (ii)}~There is a continuum $\,2^{\textstyle\aleph_0}$ of nonisomorphic finitely generated torsion-free groups in $\creg$.\hfill\qed
\end{cor}

We also point out that the third part of Examples~\ref{exos_c} applies to the Cayley graph of Gromov-hyperbolic groups:

\begin{cor}
\label{cor_hyp_in_c}%
Every nonelementary Gromov-hyperbolic group is in~$\creg$. This holds more generally for nonelementary subgroups of Gromov-hyperbolic groups.\hfill\qed
\end{cor}

\vglue-12pt
\Subsec{Stability properties and {\rm ME} invariants}
We begin with a simple observation:

\begin{lemma}
\label{lemma_finite_kernel}%
Let $\Gamma$ be a group and $F\lhd \Gamma$ a finite normal
subgroup. If $\,\Gamma/F$ is in~$\creg$ or~$\cc${\rm ,} then
so is $\Gamma$.
\end{lemma}

\Proof
If $\Gamma/F$ is in~$\creg$, then
$\hb^2(\Gamma,\ell^2(\Gamma/F))\neq 0$ because an averaging
argument gives $\hb^2(\Gamma,\ell^2(\Gamma/F))\cong
\hb^2(\Gamma/F,\ell^2(\Gamma/F))$. Since the $\Gamma$-representation $\ell^2(\Gamma/F)$ is contained in
$\ell^2(\Gamma)$, we see that $\Gamma$ is in~$\creg$. If
$\Gamma/F$ is in~$\cc$ and $\pi$ is a mixing $\Gamma/F$-representation with $\hb^2(\Gamma/F,\pi)\neq 0$, then as before
$\hb^2(\Gamma,\pi)\neq 0$. But $\pi$ is still mixing as
$\Gamma$-representation, and  so $\Gamma$ is in~$\cc$.
\hfill\qed

\begin{prop}
\label{prop_subgroup_in_C}%
Let $\Gamma$ be a group in~$\cc$ \/{\rm (}\/respectively~$\creg${\rm )}
 and $N\lhd \Gamma$ an infinite normal subgroup. Then $N$ is in~$\cc$ \/{\rm (}\/respectively~$\creg${\rm )}.
\end{prop}

The normality assumption is necessary as can be seen readily by
taking a group which is not in~$\cc$ (such as a product of two
infinite groups) and applying Corollary~\ref{cor_continuum}.
Observe that on the other hand if $\Gamma$ is a group as in 
Examples~\ref{exos_c}, then \emph{any} subgroup that acts still
\emph{nonelementarily} on the corresponding space is again in
the list of Examples~\ref{exos_c}.

\demo{Proof of Proposition~{\rm \ref{prop_subgroup_in_C}}}
 Let $\pi$ be a mixing unitary $\Gamma$-representation with
nonvanishing $\hb^2(\Gamma,\pi)$. Since $N$ is infinite, it
cannot fix any nonzero vector for $\pi$. Therefore the second
exact sequence of Theorem~12.4.2 in~\cite{Monod} (wherein one must
read $\Delta$ for $N$) shows that the restriction
$$\mathrm{res}:\ \hbc^2(\Gamma,\pi) \lra \hbc^2(N,\pi)$$
is injective and thus $\hbc^2(N,\pi)\neq 0$, so that $N$ is in~$\cc$. In the case $\pi=\ell^2(\Gamma)$, we
observe that $\pi|_N$ is a multiple of $\ell^2(N)$ so that we are done by Corollary~\ref{cor_integrate}.
\phantom{over}\Endproof \vskip4pt 

Recall that a subgroup $\Lambda$ of a group $\Gamma$ is called
\emph{co-amenable} (or one says that the coset $\Gamma$-space
$\Gamma/\Lambda$ is \emph{amenable in Eymard\/{\rm '}\/s sense}) if there is
a\break $\Gamma$-invariant mean on $\ell^\infty(\Gamma/\Lambda)$. This is, for example, 
the case if $\Lambda$ has finite index in $\Gamma$, or if
$\Lambda$ is normal in $\Gamma$ and the quotient is an amenable
group.

\begin{lemma}
\label{lemma_cofoelner_in_C}%
Let $\Gamma$ be a group in~$\cc$ \/{\rm (}\/respectively~$\creg$\/{\rm )}\/
 and $\Lambda< \Gamma$ a co-amenable subgroup. Then $\Lambda$ is also in~$\cc$
\/{\rm (}\/respectively~$\creg${\rm ).}
\end{lemma}

\Proof
Let $\pi$ be a mixing $\Gamma$-representation with $\hb^2(\Gamma,\pi)\neq 0$. By~\cite[\No8.6.2]{Monod}, the restriction map $\hb^2(\Gamma,\pi)\to \hb^2(\Lambda,\pi)$ is injective, so that $\Lambda$ is in~$\cc$ because $\pi$ is also $\Lambda$-mixing. The case of~$\creg$ is handled with Corollary~\ref{cor_integrate} as in Proposition~\ref{prop_subgroup_in_C}.
\Endproof \vskip4pt 

As implicitly observed in Section~\ref{sec_strong}, Theorem~\ref{thm_mind} has the following consequence.

\begin{cor}
\label{cor_cc_ME_invariants}%
Let $\Lambda,\Gamma$ be {\rm ME} countable groups. Assume that $\Gamma$ is in~$\cc$
 \/{\rm (}\/respectively~$\creg${\rm ).} Then $\Lambda$ is also in~$\cc$ \/{\rm (}\/respectively~$\creg${\rm ).}
\end{cor}

In other words, being in the classes~$\cc$ or~$\creg$ are ME invariants;  in particular this proves Theorem~\ref{thm_l2} from
the introduction.

\demo{Proof of the corollary}
Let  $(\Sigma, m)$ be an ME coupling of $\Lambda$ with $\Gamma$ and
let $\pi$ be a mixing $\Gamma$-representation with
$\hb^2(\Gamma,\pi)\neq 0$. By Theorem~\ref{thm_mind}, the space
$\hb^2(\Lambda, \Mind\Sigma\Gamma\Lambda\pi)$ is nonzero. By
Lemma~\ref{lemma_induce_properties} point~(i), the $\Lambda$-representation $\Mind\Sigma\Gamma\Lambda\pi$ is mixing,
and so
$\Lambda$ is in~$\cc$. In the case $\pi=\ell^2(\Gamma)$,
Lemma~\ref{lemma_induce_properties} point~(iii) gives
$\Mind\Sigma\Gamma\Lambda\ell^2(\Gamma)\cong L^2(\Sigma)$. The
latter is a multiple of $\ell^2(\Lambda)$ since there is a
$\Lambda$-fundamental domain in $\Sigma$, so that
Corollary~\ref{cor_integrate} implies that
$\hb^2(\Lambda,\ell^2(\Lambda))$ is nonzero.
\Endproof  

We can in fact refine our cohomological ME-invariants to
distinguish between some groups which are not in these two
classes.

\begin{defi}
Denote by $w\creg$ the class of countable groups $\Gamma$ admitting a
unitary representation $\pi$ which is \textit{weakly contained} in
$\ell^2(\Gamma)$ and such that $\hb^2 (\Gamma, \pi) \neq 0$.
\end{defi}

Obviously, we have $\creg \subseteq w\creg$ and as we shall
see, this inclusion is strict. All the stability properties
established above for the classes~$\cc$ and~$\creg$ remain valid,
with similar proofs, also for this class. Natural examples of
groups in $w\creg$ are provided by the following:

\begin{prop}
\label{prop_wc-c} {\rm (i)} Suppose that $N \lhd \Gamma$ is a normal
amenable subgroup. If $\Gamma/N$ is in~$\creg${\rm ,} then $\Gamma$ is in
$w\creg$.

{\rm  (ii)} If a countable group $\Gamma$ splits nontrivially as a free product over an amalgamated amenable
subgroup then $\Gamma$ is in $w\creg$ \/{\rm (}\/unless the amalgamated group has index~$\leq 2$
 in both factors\/{\rm ).}\/
\end{prop}

\Proof
(i) Because $N$ is amenable we have $\bone \prec \ell^2(N)$, and inducing both sides to $\Gamma$ gives $\ell^2(\Gamma/N) \prec \ell^2(\Gamma)$. Thus it is enough to show $\hb^2(\Gamma,\pi)\neq 0$ for $\pi = \ell^2(\Gamma/N)$. Since we have $\hb^2(\Gamma/N,\pi)\neq 0$ by assumption, the result follows from Corollary~\ref{cor_inf}.

(ii) It is enough to find  an amenable subgroup $N <\Gamma$ with $\hb^2(\Gamma,
\ell^2(\Gamma/N))\break
\neq 0$ since in the above argument the normality of $N$ in $\Gamma$ was not used for $\ell^2(\Gamma/N)
\prec \ell^2(\Gamma)$. Here, all stabilisers for the $\Gamma$-action on the set $E$ of edges of the
Bass-Serre tree associated with the amalgamated decomposition are amenable; therefore, the action on
$E\times E $ also has amenable stabilisers. Now our Corollary~7.8 in~\cite{Monod-Shalom1} completes the
proof.
\Endproof  

An argument similar to the one given in the proof of
Corollary~\ref{cor_cc_ME_invariants}, by~(ii) in
Lemma~\ref{lemma_induce_properties}, shows:

\begin{cor}
Being in the class~$w\creg$ is an {\rm ME} invariant.\hfill\qed
\end{cor}

Now that we have large families of groups in the various classes defined above, we add a few observations to complete the picture:

\begin{prop}
\label{prop_not_C}%
Let $\Gamma$ be a countable group.

\begin{itemize}
\item[{\rm (i)}] If $\Gamma$ is amenable then $\Gamma$ is not in~$\cc$ or~$\creg${\rm ,} nor in~$w\creg$.

\item[{\rm (ii)}] If $\Gamma$ contains an infinite normal amenable \/{\rm (}\/e.g.\ central\/{\rm )}\/
 subgroup $N$ then $\Gamma$ is not in~$\cc$ or~$\creg${\rm ;}
 but if $\Gamma/N$ is in $\creg${\rm ,} then $\Gamma$
is in $w\creg$.

\item[{\rm (iii)}] If $\Gamma$ is the product of two infinite
groups then $\Gamma$ is not in $\cc$ or $\creg${\rm ,} and if both
groups are nonamenable then $\Gamma$ is also not in $w\creg$ \/{\rm (}\/though $\Gamma$ may
be in $w\creg$ if one of the factors is amenable  as in~{\rm (ii)} above\/{\rm )}\/.

\item[{\rm (iv)}] If $\Gamma$ has infinitely many ends{\rm ,} i.e.\ \/{\rm (}\/by
Stallings\/{\rm )}\/ if it is a nontrivial free product over a
\textrm{finite} amalgamated subgroup{\rm ,} then $\Gamma$ is in $\cc$
and $\creg$. If $\Gamma$ is a nontrivial free product over an
\textrm{amenable} subgroup then it is in $w\creg$.

\item[{\rm (v)}] $\Gamma$ is not in any of these classes if it is a lattice in a higher rank simple Lie group{\rm ,}
or in a higher rank simple algebraic group over a local field.

\item[{\rm (vi)}] $\Gamma$ is not in any of these classes if it is an irreducible lattice in a product 
of nonamenable compactly generated locally compact groups.
\end{itemize}
\end{prop}

\Proof
For~(i), see Remark~\ref{rem_amenable_group}. The first part of~(ii) follows from Proposition~\ref{prop_kernel} since an infinite subgroup cannot fix a nonzero vector in a mixing representation; the second part was noted in Proposition~\ref{prop_wc-c}. For point~(iii), assume that $\hb^2(\Gamma,\pi)\neq 0$ for $\Gamma = \Gamma_1\times \Gamma_2$. Theorem~\ref{thm_Kunneth} implies $\H_\pi^{\Gamma_i}\neq 0$ for some $i\in\{1,2\}$. This forces $\Gamma_i$ to be finite if $\pi$ is mixing and to be amenable if $\pi\prec\ell^2(\Gamma)$. Point~(iv) has been addressed above, see the above discussion after the proof of Theorem~\ref{thm_exos_c} and Proposition~\ref{prop_wc-c}. One gets~(v) by applying
 Theorem~1.4 from~\cite{Monod-Shalom1}. Theorem~16 in~\cite{Burger-Monod3} implies~(vi). 
\Endproof \vskip4pt 

At this point, we can complete the

\demo{Proof of Corollary~{\rm \ref{cor_hyperbolic}}}
The result follows from the juxtaposition of Corollary~\ref{cor_cc_ME_invariants} and point~(ii) of Proposition~\ref{prop_not_C}.
\Endproof \vskip4pt 

Next, recall that a group is said to be ICC if the conjugacy
class of every nontrivial element is infinite. We used the
following fact each time we needed Lemma~\ref{lemma_Dirac_inv}.

\begin{prop}
\label{prop_ICC_C}%
Any countable torsion-free group in~$\cc$ is {\rm ICC}.
\end{prop}

\Proof
Let $\Gamma$ be as in the statement and suppose for a
contradiction that it contains a nontrivial element $\gamma_0$
with finite conjugacy class. The centraliser $\Gamma_0$ of
$\gamma_0$ has finite index in $\Gamma$ so that it is also
in~$\cc$. The
subgroup $C_0$ generated by $\gamma_0$ is normal in $\Gamma_0$ and
amenable, so by Proposition~\ref{prop_kernel} the space
$\H_\pi^{C_0}$ is nonzero for every $\Gamma_0$-representation
such that $\hb^2(\Gamma_0, \pi)\neq 0$. However this can never
happen for $\pi$ mixing unless $C_0$ is finite, a contradiction
since $\Gamma$ is torsion-free and $\gamma_0$ is nontrivial.
\hfill\qed

\Subsec{A shortcut to~$\cc$}
\label{sec_shortcut}%
Finally, we describe a short alternative approach to establish that certain groups (including free groups) are in
the class~$\cc$, thereby avoiding completely the dependence on the companion
paper~\cite{Monod-Shalom1}. The more restricted family we cover here is still rich enough to provide many
interesting examples to which our foregoing results apply, as well as to establish
Theorem~\ref{thm_continuum} (see below) and Theorem~\ref{thm_out} (with some modification of the
proof).

\begin{prop}
\label{prop_rank_one}%
Every lattice in a simple \/{\rm (}\/connected{\rm ,} center-free\/{\rm )}\/ Lie group of rank one is in~$\cc$.
\end{prop}

Before proving the proposition we notice that it yields yet another proof (in the finitely generated case) that non-Abelian free groups are in~$\cc$. Using this fact, we can give an alternative argument (independent of~\cite{Monod-Shalom1}) to show that \itshape there is a continuum of nonisomorphic finitely generated torsion-free groups in~$\cc$\upshape\ (Corollary~\ref{cor_continuum}):

Indeed, if $A,B$ are any two finitely generated amenable groups, then the argument given in the proof of Theorem~\ref{thm_free_continuum} above shows that the free product $A*B$ is ME to a free group on two generators. Thus $A*B$ is in~$\cc$ by Corollary~\ref{cor_cc_ME_invariants}. In other words, it remains only to justify that there is a continuum of nonisomorphic finitely generated torsion-free amenable groups. This is true even for \emph{soluble} groups; P.~Hall proves in~\cite{Hall} that there are uncountably many (in fact, a continuum of) nonisomorphic groups $G$ on two generators with $[G'',G]=1$. A close examination of his proof (Theorem~6 and pages~433--435 in~\cite{Hall}) shows that $G$ can be chosen torsion-free.

\demo{Proof of Proposition~{\rm \ref{prop_rank_one}}}
Let $G$ be a rank-one Lie group. Since all lattices in $G$ are ME to each other, it is enough by
Corollary~\ref{cor_cc_ME_invariants} to prove the statement for \emph{some} lattice $\Gamma<G$. We
are therefore free to choose $\Gamma$ cocompact in $G$, which implies that $\Gamma$ is hyperbolic in the
sense of Gromov and therefore the bounded cohomology with \emph{trivial} coefficients $\hb^2(\Gamma)$
is infinite-dimensional~\cite{Epstein-Fujiwara}.

We now consider the quasi-regular $G$-representation $\pi=L^2(G/\Gamma)$ which splits as $\bone \oplus \pi_0$, where $\pi_0$ denotes the kernel of the orthogonal projection $p:\pi\to\bone$ to the constants. Observe that the induced $\Gamma$-representation $\Mind G\Gamma\Gamma\bone$ as defined in Section~\ref{sec_ind_rep} through the ME self-coupling of $\Gamma$ given by the right and left $\Gamma$-actions on $G$ is the restriction $\pi|_\Gamma$ of $\pi$ to $\Gamma$. In fact, in the present case
every definition of cohomological induction
$\MMind{G}\Gamma\Gamma:\hb^\bu(\Gamma)\to\hb^\bu(\Gamma, \Mind G\Gamma\Gamma\bone)$
through this self-coupling (Section~\ref{sec_ind_coh}) coincides with the map
$$\hb^\bu(\Gamma)\xrightarrow{\ \mathbf{i}_\Gamma^G\ } \hbc^\bu(G, \pi) \xrightarrow{\ \mathrm{res}\ } \hb^\bu (\Gamma, \pi|_\Gamma)$$
wherein $\mathbf{i}_\Gamma^G$ is the induction to continuous bounded cohomology $\hbc^\bu$ as 
defined in~\cite{Burger-Monod3}, \cite{Monod} and $\mathrm{res}$ is  the restriction map. Thus, we have
a commutative diagram:
$$\xymatrix{
& \hb^\bu(\Gamma,\pi|_\Gamma) \ar[r]^{p_*} & \hb^\bu(\Gamma)\\
\hb^\bu(\Gamma) \ar[ur]^{\MMind{G}\Gamma\Gamma} \ar[dr]_{\mathbf{i}_\Gamma^G} & & \\
& \hbc^\bu(G,\pi) \ar[r]^{p_*} & \hbc^\bu(G) \ar[uu]_{\mathrm{res}}\\
}$$
However, in degree two we have the following additional information: (i)~the space $\hbc^2(G)$ has dimension at most one
 (Lemma~6.1 in~\cite{Burger-Monod1}; compare~\cite{Burger-Monod3});
(ii)~$\MMind{G}\Gamma\Gamma$ is injective (Theorem~\ref{thm_mind}). It follows that
$\hb^2(\Gamma,\pi_0|_\Gamma)$ is infinite-dimensional (bounded cohomology is additive in the
coefficients~\cite[8.2.10]{Monod}). This finishes the proof because $\pi_0$ is a mixing representation by
the Howe-Moore theorem.
\hfill\qed

%

\references {CFW2}

\bibitem[A1]{Adams94a}
\name{S.\ Adams},  Indecomposability of equivalence relations generated by word
  hyperbolic groups, {\it Topology\/} {\bf 33} (1994),  785--798.

\bibitem[A2]{Adams95}
\bibline, Some new rigidity results for stable orbit equivalence,
  {\it Ergodic Theory Dynam.\ Systems\/} {\bf 15} (1995), 209--219.


\bibitem[BG]{BG81}
\name{S.\ Bezugly\u {\i}} and \name{V.\ Golodets},
 Hyperfinite and II$_1$ actions for nonamenable groups,
{\it  J. Funct.\  Anal\/}.\ {\bf 40} (1981), 30--44.


\bibitem[BM1]{Burger-Monod1}
\name{M.\ Burger} and \name{N.\ Monod},  Bounded cohomology of lattices in higher
  rank {L}ie groups, {\it J.\ European Math.\ Soc\/}.\ ({\it JEMS\/}) 
  {\bf 1} (1999),
  199--235.

\bibitem[BM2]{Burger-Monod3}
\bibline, Continuous bounded cohomology and
  applications to rigidity theory (with an appendix by {M}.~{B}urger and
  {A}.~{I}ozzi), {\it Geom.\ Funct.\ Anal\/}.\ {\bf 12} (2002), 219--280.

\bibitem[BMz]{Burger-Mozes2}
\name{M.\ Burger} and \name{S.\ Mozes},  Lattices in product of trees, {\it Inst. Hautes
  \'Etudes Sci.\ Publ.\ Math\/}.\ (2000),  151--194 (2001).

\bibitem[CFW]{Connes-Feldman-Weiss}
\name{A.\ Connes, J.\ Feldman}, and \name{Benjamin Weiss},  An amenable equivalence
  relation is generated by a single transformation, {\it Ergodic Theory 
  Dynam.\ Systems\/} {\bf 1} (1981),  431--450.

\bibitem[Dx]{Dixmier64}
\name{J.\ Dixmier},  Utilisation des facteurs hyperfinis dans la th\'eorie des
  ${C}\sp{\ast} $-alg\`ebres, {\it C. R. Acad.\ Sci.\ Paris\/} {\bf 258} (1964),
  4184--4187.

\bibitem[Dy]{Dye}
\name{H.\ A. Dye},  On groups of measure-preserving transformations. {I}, {\it Amer.\
  J. Math\/}.\ {\bf 81} (1959), 119--159.

\bibitem[EpF]{Epstein-Fujiwara}
\name{D.\ B.~A. Epstein} and \name{K.\ Fujiwara},  The second bounded cohomology of
  word-hyperbolic groups, {\it Topology\/} {\bf 36} (1997),  1275--1289.

\bibitem[EF1]{Eskin-Farb1}
\name{A.\ Eskin} and \name{B.\ Farb},  Quasi-flats and rigidity in higher rank
  symmetric spaces, {\it J. Amer.\ Math.\ Soc\/}.\ {\bf 10} (1997),  653--692.

\bibitem[EF2]{Eskin-Farb2}
\bibline,  Quasi-flats in ${H}\sp 2\times {H}\sp 2$, in
{\it Lie
  Groups and Ergodic Theory\/} (Mumbai, 1996) (Bombay), {\it 
  Tata Inst.\ Fund.\ Res\/}.\ {\bf 14}  
  (1998), 75--103.

\bibitem[E]{Epstein71}
\name{D.\ B.~A. Epstein},  Almost all subgroups of a {L}ie group are free, {\it J.\
  Algebra\/} {\bf 19} (1971), 261--262.

\bibitem[F1]{Furman1}
\name{A.\ Furman},  Gromov's measure equivalence and rigidity of higher rank
  lattices, {\it Ann.\ of Math.\/}.\  {\bf 150} (1999),  1059--1081.

\bibitem[F2]{Furman2}
\bibline,   Orbit equivalence rigidity, {\it Ann.\ of Math\/}.\ {\bf 150}
  (1999),  1083--1108.

\bibitem[F3]{Furman3}
\bibline,   Outer automorphism groups of some ergodic equivalence
  relations, {\it Comment.\ Math.\ Helv\/}.\ {\bf 80} (2005), 157--196.

\bibitem[Ga1]{Gaboriau00}
\name{D.\ Gaboriau},  Co\^ut des relations d'\'equivalence et des groupes,
  {\it Invent.\ Math\/}.\ {\bf 139} (2000), no.~1, 41--98.

\bibitem[Ga2]{GaboriauCRAS}
\bibline,   Sur la (co-)homologie ${L}\sp 2$ des actions pr\'eservant
  une mesure, {\it C. R. Acad.\ Sci.\ Paris S{\hskip.5pt\rm \'{\hskip-5pt\it e}}r.\ I Math\/}.\ {\bf 330} 
  (2000), 
  365--370.

\bibitem[Ga3]{GaboriauL2}
\bibline,  Invariants {$l\sp 2$} de relations d'\'equivalence et de
  groupes, {\it Publ.\ Math.\ Inst.\ Hautes \'Etudes Sci\/}.\ {\bf 95} (2002), 
  93--150.

\bibitem[Ge1]{Gefter93}
\name{S.\ L.\ Gefter},  Ergodic equivalence relation without outer
  automorphisms, {\it Dopov.\ Dokl.\ Akad.\ Nauk Ukra\"\i ni\/} {\bf 11} 
  (1993),  25--27.

\bibitem[Ge2]{Gefter96}
\bibline,  Outer automorphism group of the ergodic equivalence
  relation generated by translations of dense subgroup of compact group on its
  homogeneous space, {\it Publ.\ Res.\ Inst.\ Math.\ Sci\/}.\ {\bf 32} (1996), 
  517--538.

\bibitem[GP]{Gaboriau-Popa}
\name{D.\ Gaboriau} and \name{S.\ Popa}, An uncountable family of nonorbit
  equivalent actions of $\mathbf{F}_n$, {\it J.\ Amer.\ Math.\ Soc\/}.\
  {\bf 18} (2005), 547--599.

\bibitem[Gr]{Gromov}
\name{M.\ Gromov}, Volume and bounded cohomology, {\it Inst.\ Hautes \'Etudes
  Sci.\ Publ.\ Math\/}.\ {\bf 56} (1982),  5--99 (1983).

\bibitem[Ha]{Hall}
\name{P.\ Hall},  Finiteness conditions for soluble groups, {\it Proc.\ London Math.\
  Soc\/}.\ {\bf 4} (1954), 419--436.

\bibitem[Hj]{Hjorth_Dye}
\name{G.\ Hjorth},  A converse to {D}ye's theorem,  {\it Trans.\ Amer.\ Math.\ Soc\/}.\ {\bf 357} (2005), 3083--3103.

\bibitem[HK]{Hjorth-Kechris}
\name{G.\ Hjorth} and \name{A.\ ~S. Kechris},  Rigidity theorems for actions of
  product groups and countable Borel equivalence relations, {\it Memoirs
  Amer.\  Math.\ Soc\/}.\ {\bf 177}, A.\ M.\ S., Providence, RI
  (2005).

\bibitem[Hu]{Hulanicki66}
\name{A.\ Hulanicki},  Means and {F}\o lner condition on locally compact
  groups, {\it Studia Math\/}.\ {\bf 27} (1966), 87--104.

\bibitem[I]{Ivanov}
\name{N.\ V.\ Ivanov},  Foundations of the theory of bounded cohomology, {\it
J.\ Soviet  Math\/}.\ {\bf 37}
   (1987), 1090--1115.

\bibitem[J]{Johnson}
\name{B.\ E. Johnson},  {\it Cohomology in Banach Algebras}, 
{\it Memoirs Amer.\ Math.\ Soc\/}.\ 
{\bf 127}, A.\ M.\ S., Providence, RI (1972).

\bibitem[K]{Kaimanovich03}
\name{V.\ A. Kaimanovich},  Double ergodicity of the {P}oisson boundary and
  applications to bounded cohomology, {\it Geom.\ Funct.\ Anal\/}.\ 
  {\bf 13} (2003),
  852--861.

\bibitem[KL]{Kleiner-Leeb}
\name{B.\ Kleiner} and \name{B.\ Leeb}, Rigidity of quasi-isometries for
  symmetric spaces and {E}uclidean buildings, {\it Inst.\ 
  Hautes \'Etudes Sci.\ Publ.\
  Math\/}.\ {\bf 86} (1997), 115--197 (1998).

\bibitem[MMS]{Mineyev-Monod-Shalom}
\name{I.\ Mineyev, N.\ Monod}, and \name{Y.\ Shalom},  Ideal bicombings for
  hyperbolic groups and applications, {\it Topology\/} {\bf 43} (2004),
  1319--1344.

\bibitem[M]{Monod}
\name{N.\ Monod}, {\it Continuous Bounded Cohomology of Locally Compact Groups},
  {\it Lecture Notes in Math\/}.\ {\bf 1758}, Springer-Verlag, New York,
  2001.

\bibitem[MS1]{Monod-ShalomCRAS}
\name{N.\ Monod} and \name{Y.\ Shalom},  Negative curvature from a cohomological
  viewpoint and cocycle superrigidity, {\it C.\ R.\  Acad.\ Sci.\ Paris S{\hskip.5pt\rm \'{\hskip-5pt\it e}}r.\ 
  I Math\/}.\
  {\bf 337} (2003),  635--638.

\bibitem[MS2]{Monod-Shalom1}
\bibline,   Cocycle superrigidity and bounded
  cohomology for negatively curved spaces, {\it J.\ Differential Geom\/}.\ 
   {\bf 67}
  (2004), 395--455.

\bibitem[N]{Noskov}
\name{G.\ A.\ Noskov},  Bounded cohomology of discrete groups with
  coefficients, {\it Leningr. Math. J\/}.\ {\bf 2} (1991),  1067--1084.

\bibitem[OP]{Ozawa-Popa}
\name{N.\ Ozawa} and \name{S.\ Popa},  Some prime factorization results for type
  ${\rm II}\sb 1$ factors, {\it Invent.\ Math\/}.\ {\bf 156} (2004), 
  223--234.

\bibitem[OW]{Ornstein-Weiss}
\name{D.\ S.\ Ornstein} and \name{B.\ Weiss},  Ergodic theory of amenable group
  actions.\ {I}.\ {T}he {R}ohlin lemma, {\it Bull.\ Amer.\ Math.\ Soc\/}.\ {\bf 2}
  (1980),  161--164.

\bibitem[Po]{Po} \name{S. Popa}, 
Deformation and rigidity for group actions and von Neumann
algebras, Plenary lecture at ICM Madrid 2006,
{\it Proc.\  of the International Congress of Mathematicians\/} (Madrid,
Spain, 2006), European Math.\ Soc., to appear.



\bibitem[Pr]{Prasad76}
\name{G.\ Prasad},  Discrete subgroups isomorphic to lattices in semisimple
  {L}ie groups, {\it Amer.\ J. Math\/}.\ {\bf 98} (1976),  241--261.

\bibitem[Sh1]{ShalomANNALS}
\name{Y.\ Shalom},  Rigidity, unitary representations of semisimple groups, and
  fundamental groups of manifolds with rank one transformation group, 
  {\it Ann.\ of
  Math\/}.\ {\bf 152} (2000),  113--182.

\bibitem[Sh2]{Sh2}
\bibline, Measurable group theory,
{\it European Congress of Mathematics}, 391--423,
European Math.\ Soc., Z\"urich, 2005.


\bibitem[SW]{Schmidt-Walters}
\name{K.\ Schmidt} and \name{P.\ Walters},  Mildly mixing actions of locally compact
  groups, {\it Proc.\ London Math.\ Soc\/}.\ {\bf 45} (1982),  506--518.

\bibitem[Z1]{Zimmer83}
\name{Robert~J. Zimmer},  Ergodic actions of semisimple groups and product
  relations, {\it Ann.\ of Math\/}.\ {\bf 118} (1983),  9--19.

\bibitem[Z2]{Zimmer84}
\bibline,   {\it Ergodic Theory and Semisimple Groups\/}, Birkh\"auser
  Verlag, Basel, 1984.

\Endrefs
\def\cprime{$'$}
\ifx\undefined\bysame
\newcommand{\bysame}{\leavevmode\hbox to3em{\hrulefill}\,}
\fi

\end{document}